\documentclass[12pt]{article}
\usepackage{arxiv}
\usepackage[noblocks]{authblk}

\usepackage{amsmath}
\usepackage{amsthm}
\usepackage{amssymb}
\usepackage{amsfonts}
\usepackage[T1]{fontenc}
\usepackage[symbol]{footmisc}
\usepackage[style=alphabetic, backend=biber, natbib=true,url=false, doi=false,isbn=false,eprint=true]{biblatex}
\usepackage{dsfont}
\usepackage{mathtools,mathrsfs,bm}
\usepackage{import,booktabs,times,fancyhdr}
\usepackage{nicefrac,microtype,xcolor,xspace}
\usepackage{graphicx,tikz,tcolorbox}
\usetikzlibrary{decorations.pathreplacing}

\usepackage{algorithm,algpseudocode}
\usepackage{booktabs}
\usepackage{tabularx}
\usepackage{array} 
\usepackage[inline]{enumitem}
\usepackage{subfigure}
\usepackage[parfill]{parskip}
\usepackage{hyperref,cleveref}

\usepackage{stmaryrd}
\usepackage{attention_macros}
\usepackage[normalem]{ulem}
\usepackage{tcolorbox}
\usepackage{standalone}
\usepackage{caption}
\usepackage{floatpag}

\usepackage{mathMacros}

\definecolor{darkorange}{RGB}{215, 155, 0}
\definecolor{lightorange}{RGB}{255, 230, 204}
\definecolor{darkgreen}{RGB}{130, 179, 102}
\definecolor{lightgreen}{RGB}{213, 232, 212}
\definecolor{darkyellow}{RGB}{214, 182, 86}
\definecolor{lightyellow}{RGB}{255, 242, 204}
\definecolor{darkblue}{RGB}{120, 165, 215}
\definecolor{lightblue}{RGB}{218, 232, 252}

\usetikzlibrary{matrix,positioning,decorations.pathreplacing,arrows, calc,patterns, calc, 3d,arrows.meta}
\algrenewcommand\algorithmicrequire{\textbf{Input:}}
\algrenewcommand\algorithmicensure{\textbf{Output:}}

\title{Attention Mechanisms Through the Lens of Numerical Methods: Approximation Methods and Alternative Formulations}
\renewcommand{\shorttitle}{Attention Mechanisms Through the Lens of Numerical Methods}

\author[1,*]{Michel Fabrice Serret}
\author[2,*]{Alice Cortinovis}
\author[3]{Yijun Dong}
\author[3]{Diana Halikias}
\author[4]{Anna Ma}
\author[5]{Fabio Matti}
\author[6]{Deanna Needell}
\author[7]{Katherine J. Pearce}
\author[8]{Elizaveta Rebrova}
\author[9]{Disha Shur}
\author[10]{Rudi Smith}
\author[11]{Hai-Xiao Wang}
\author[1,5,$\dagger$]{Laura Grigori}

\affil[1]{Center for Scientific Computing, Theory and Data, Paul Scherrer Institute, Switzerland (\href{mailto:michel.serret@psi.ch}{michel.serret@psi.ch})}
\affil[2]{Department of Computer Science, University of Pisa, Italy and member of INdAM/GNCS}
\affil[3]{Courant Institute, New York University, NY, USA}
\affil[4]{Department of Mathematics, University of California, Irvine, CA, USA}
\affil[5]{Institute of Mathematics, EPFL, Switzerland}
\affil[6]{Department of Mathematics, University of California, Los Angeles, CA, USA}
\affil[7]{Oden Institute, University of Texas, Austin, TX, USA}
\affil[8]{Department of ORFE, Princeton University, NJ, USA}
\affil[9]{Department of Computer Science, Purdue University, IN, USA}
\affil[10]{Department of Mathematics, Virginia Tech, VA, USA}
\affil[11]{Department of Mathematics, University of Wisconsin-Madison, WI, USA.}

\affil[*]{Project group co-leader}
\affil[$\dagger$]{Project group leader}
\begin{document}
\maketitle

\begin{abstract} The attention mechanism is the computational core of modern Transformer architectures, but its quadratic complexity in the input sequence length is the bottleneck for large-scale inference. This has motivated a rapidly growing body of work aimed at accelerating attention through approximation and reformulation. In this survey, we revisit attention mechanisms through the lens of numerical analysis, with a particular emphasis on tools and perspectives from numerical linear algebra.
Our goal is twofold: first, we aim to systematically review and classify fast approximation methods according to the numerical principles they exploit. These include sparsity and clustering approaches, low-rank and subspace projection techniques, randomized sketching methods, and tensor-based decompositions. We also discuss kernel-inspired reformulations of attention and recent architectural variants, such as Latent Attention, that modify the standard softmax formulation to improve efficiency.
Second, by presenting these developments within a unified mathematical framework, we aim to bridge the gap between disciplines and highlight opportunities for further contributions from computational mathematics, particularly numerical linear algebra, to the design of scalable attention mechanisms.
\end{abstract}
\section{Introduction}
In recent years, large language models \cite{openaiGptoss120bGptoss20bModel2025,grattafioriLlama3Herd2024,yangQwen3TechnicalReport2025,teamGemma3Technical2025,guo2025deepseek} based on the attention mechanism \cite{vaswani2017attention} (commonly called LLMs) have become ubiquitous across industry, academia, and society at large, enabling the automation of repetitive and time-consuming tasks.
However, the computational and energy requirements of such models are considerable and represent a nontrivial fraction of total computing resources worldwide.
This paper investigates the mechanism of attention, with a focus on how numerical methods can be leveraged to enhance its efficiency. We discuss different contributions from the  areas of approximation theory and numerical analysis, namely data-sparse representations such as low-rank approximation of matrices and tensors and kernel-based methods used to reduce the computational cost of the attention mechanism. The primary motivation for these methods comes from numerical experiments, across different models,
suggesting that attention matrices can exhibit sparsity or a low-rank structure (see, e.g., \cite{chen2021scatterbrain}), potentially enabling efficient and accurate approximations. Additionally, observed phenomena such as attention sinks \cite{xiaoEfficientStreamingLanguage2024a} and grokking \cite{powerGrokkingGeneralizationOverfitting2022} suggest an underlying structure in attention that could be leveraged to improve efficiency.

Although previous surveys \cite{tayEfficientTransformersSurvey2022,geshkovski2025mathematical} provide a broad view of transformer models, we adopt a more focused approach by studying the attention mechanism and introducing a taxonomy of existing methods for its approximation. We also discuss methods that are fundamentally different from the scaled dot-product attention used in most models today. 

Furthermore, we focus mainly on inference, i.e., the evaluation of a pre-trained model at runtime.  Indeed, we decouple our analysis from training considerations whenever possible, instead referring to works such as Muon \cite{liuMuonScalableLLM2025}, SOAP \cite{vyasSOAPImprovingStabilizing2025}, and related approaches that target improvements in training efficiency. Additionally, we choose to specifically concentrate on the efficiency of the models themselves and not on that of implementations, such as FlashAttention, see~\cite{dao2022flashattentionfastmemoryefficientexact,dao2023flashattention2fasterattentionbetter}. Although not at the core of this work, we provide a brief overview of 
theoretical results concerning the expressivity, learnability, and clustering properties of transformers and of fast attention in Appendix~\ref{sec:tf_theory}.

\paragraph{Outline.} Section \ref{sec:original_notation} introduces the attention mechanism and its mathematical formulations that we will use throughout the paper. 
The remaining sections are organized according to our taxonomy of methods for accelerating attention computation, as summarized in Figure \ref{fig:overview}. We start by reviewing methods based on sparsity and clustering in Section~\ref{sec:clustering}. These strategies leverage the fact that matrices coming from attention computations usually have only a few important entries, employing a variety of techniques to detect and prioritize these entries, with potentially less accurate approximation of the smaller quantities. In Section~\ref{sec:low-rank-approximation} we review methods that use low-rank approximations and projections onto low-dimensional subspaces, which are based on the fact that the matrices involved in the computation of attention often exhibit singular value decay. Section~\ref{sec:kernel-based} discusses the close connection between attention and kernel methods. We focus on (i) how information derived from kernels can be used to approximate the attention and (ii) how the standard formulation of the attention itself can be replaced by a kernel-based formulation that retains the important characteristics of the standard attention but is faster to evaluate. 
We then provide, in Section~\ref{sec:MLA}, a description of a recent modification of the attention mechanism known as \emph{Latent Attention}, that allows for improved efficiency, and consider approximate conversion between regular and latent attention models. Finally, in Section~\ref{sec:tensors} we discuss tensor-based methods, concentrating on three aspects: (i) tensorized matrices across the whole transformer and the application of low-rank tensor decomposition techniques, (ii) incorporation of tensors directly into the attention model, and (iii) the preservation or exploitation of tensor-structured input data. 
\paragraph{Notation.}
Throughout the paper, capital letters generally denote matrices and lowercase letters denote vectors. Indices are denoted by $i,j$, and we write $[N]=\{1,\ldots,N\}$ and $\mathbf{1}_N=(1)_{1\le i\le N}\in\mathbb{R}^N$. Matrix entries are denoted interchangeably by $A_{ij}$ and $A(i,j)$; rows of a matrix $K$ are written as $K(i,:)$, $K_{i,:}$, or using symbols such as $k_i$, while vector entries are denoted by $x_i$, with no ambiguity from context. Functions are denoted by $f$, Greek letters (e.g., $\varphi$), or fraktur symbols (e.g., $\mathfrak{g}$). We reserve $\kappa$ for kernel functions and $\pi$ for permutation vectors. Inner products are written as $\langle \cdot,\cdot\rangle$. The Kronecker, Khatri-Rao, and Hadamard products are denoted by $\otimes$, $*$, and $\odot$, respectively, and $q^{\otimes p}$ denotes the $p$-fold Kronecker product of a vector $q$. The expressions $\|\cdot\|_{\mathrm{F}} $ and $\|\cdot\|_{\textnormal{op}}$ denote the Frobenius and operator (spectral) norm of matrices, respectively, and $\|\cdot\|$ is the Euclidean norm of a vector. We denote by $\indi{\text{event}}$ the indicator function of the event. $\mathbb{N}^*$ is the set of positive integer and $\RR_+$ is the set of non-negative real numbers. Additional notation specific to tensors is introduced in the corresponding section.

\section{Mathematical  model of attention}
\label{sec:original_notation}

\begin{figure}[h!]
\centering
    \includegraphics[width=0.6\textwidth]{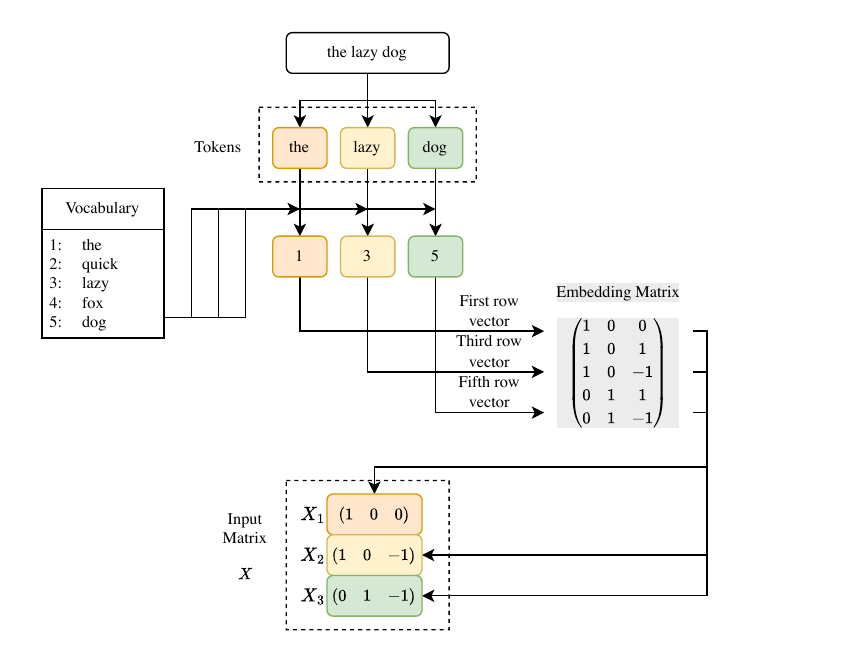}
\caption{Example of sentence embedding for a word-based tokenizer for the phrase 'the lazy dog'.}
\label{fig:tok_embd}
\end{figure}

Attention mechanisms are at the core of modern, Transformer-based \cite{vaswani2017attention}, language models.
These models are able to construct vector spaces encoding the semantic information contained not only in individual words but also in the relationships between words in a text. This encoding of the semantics between elements of a sentence occurs through the attention mechanism.
Before defining exactly how attention works, let us first provide some background on how these methods process textual information into vectors and how the attention mechanism processes these vectors.
\paragraph{Transformers and Natural Language Processing}
Modern language models are not able to directly use the usual string formatting used to encode text.
Instead, they require as input a sequence of vectors corresponding to the text. To create these sequences, the text is cut into small successive chunks of strings called tokens. The choice of the complete set of tokens used by the model, known as the vocabulary, is important.  For example, in the context of a language model, it needs to allow the encoding of enough of the semantic information contained in a text to achieve the task of predicting the next token in a sequence. At the same time, it should be sufficiently small  to avoid unnecessary computational overheads while  efficiently representing the semantic information contained in the text. 
The tokenization procedures tend to be language specific and, ideally, should encode the smallest set of strings with inherent meaning for a given language. 
Second, once the text has been decomposed into a sequence of tokens, a vectorization procedure is required. In what follows, tokens can be taken without loss of generality to be words. 
Given a vocabulary, an embedding step allows to represent each of its element in a vector space through a trained linear application. This is also referred to as vectorization, since each token is thus represented by a vector, a visual example is given in Figure \ref{fig:tok_embd}. 

\begin{figure}[h!]
\centering
    \includegraphics[width=\textwidth]{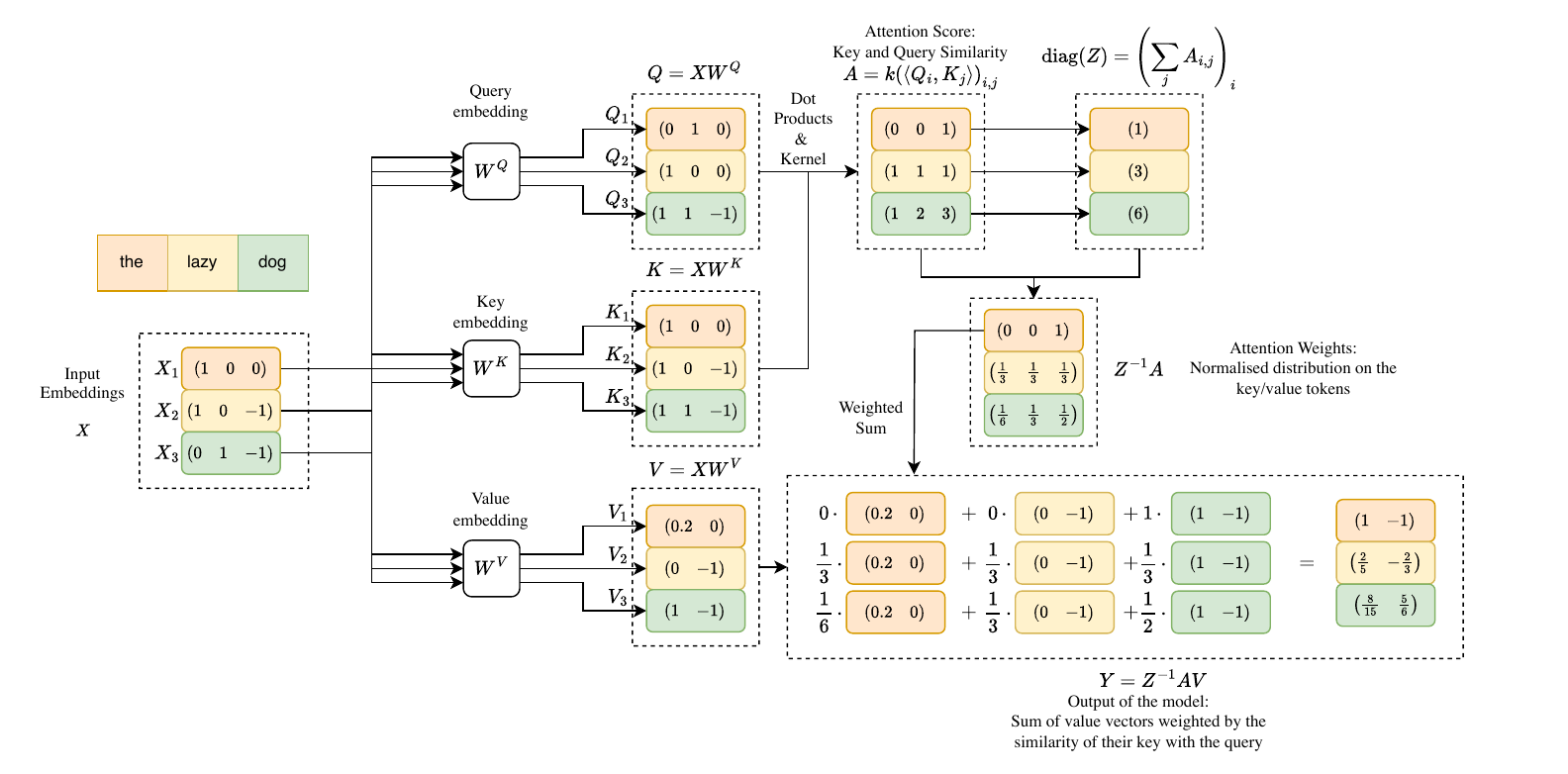}
\caption{Example of an attention layer for the embeddings obtained from the tokenization example in Figure \ref{fig:tok_embd}. It is important to note that in this example we take the kernel to be linear for simplicity.}
\label{fig:attn_meca}
\end{figure}

Attention mechanisms, the building blocks of the Transformer architecture \cite{vaswani2017attention}, allow encoding of semantic information between token embeddings through a database-like structure, wherein, for each query, the database outputs the stored value associated to the key matching the query.
Analogously, given the set of $N$ embedded tokens of dimension $d$, Figure~\ref{fig:attn_meca} illustrates how to obtain, given a set of input vectors $X\in\RR^{N\times d}$, associated to the tokens, the output of the attention mechanism. 
For each token $x\in X$, its associated query, key and value vectors are obtained through a set of linear transformations. Then, output is given as a linear combination of the value embeddings associated with each token $x'\in X$, weighted by a similarity metric relating the query embedding $q$ of $x$ to the 'key' embedding $k'$ of $x'$. Let us now define the attention mechanism more formally.

\begin{tcolorbox}
    In the rest of this section, additional notes will be given inside these gray boxes for those who may be unfamiliar with deep learning models or Natural Language Processing. These notes are not necessary to understand the attention mechanism but provide some additional information about the context in which attention mechanisms are used. For more information on how attention mechanisms are used in transformer-based language models, we refer the interested reader to \cite{serretUnderstandingTransformersAttention2026}, the introduction to transformers and attention mechanisms written for the IPAM workshop at the origin of this review.
\end{tcolorbox}
\subsection{Attention mechanism}
The computation of the standard basic attention mechanism  \cite{vaswani2017attention}, also called the scaled dot-product attention, can be 
expressed as a sequence of linear algebraic operations followed by a nonlinear softmax transformation.
Suppose that the input data consists of $N$ tokens, each of dimension $d$, recorded as a matrix $X\in\RR^{N\times d}$. The dimension $d$ is also called the ``hidden size", as it is the dimension of the hidden state of the neural network, i.e. the size of the vectors sent between successive layers of a model. In most modern LLMs, this size is kept constant throughout the hidden layers of the model for efficiency, and in what follows we will also consider $d$ to be a constant of the model. Three linear maps $\WQ\in\RR^{d\times \dhead}$ (\emph{query weights}), $\WK\in\RR^{d\times \dhead}$  (\emph{key weights}) and $\WV\in\RR^{d\times d}$ (\emph{value weights}) for a model are obtained through separate pre-training and fine-tuning procedures. Interested readers can find more details in the original GPT article \cite{radford2018improving}.  

Then, \emph{Query}, \emph{Key} and \emph{Value matrices} are formed respectively as
\begin{align}
    \Q=X\WQ\in \RR^{N\times \dhead},\quad \K=X\WK\in \RR^{N\times \dhead},\qquad \V=X\WV\in\RR^{N\times d}.\label{eqn:QKV}
\end{align}

Based on these matrices, the \emph{attention scores matrix} $\A\in \RR^{N\times N}$ is defined as
\begin{align}  \A=\exp\left(\frac{\Q\K^{\T}}{\sqrt {\dhead}}\right)=\left(\kappa_{\exp}(q_i,k_j)\right)_{1\le i,j\le N}, \label{eqn:score_matrix}
\end{align}
where the function $\exp$ is applied to matrices entry-wise and, for any $j\in \NN$, we denote the row vectors of the matrices $Q$ and $K$ as $q_j$ and $k_j$ respectively. Furthermore, we note that this can equivalently be seen as the application of the kernel function $\kappa_{\exp}:\RR^{\dhead} \times \RR^{\dhead} \to \RR$ given by $(x,y)\mapsto\exp\left (\frac{\langle x,y\rangle}{\sqrt{\dhead}}\right )$ on key and query vectors.

The diagonal \emph{normalization matrix} $\Z\in\RR^{N\times N}$ is defined entry-wise as 
\begin{align}    \Z(i,i)=\sum\limits_{j=1}^{N}\A(i,j) \quad \text{ for } i = 1, \ldots, N,\label{eqn:diagonal_matrix}
\end{align}
with $Z(i,j) = 0$ for $i \ne j$.
Finally, the  output of the attention mechanism is given by $\Att\in\RR^{N\times d}$ defined by the normalized application of attention scores to the value matrix $V$,
\begin{align}
    \Att=\AttWeights \V = \Z^{-1}\A\V, \label{eqn:att_matrix}
\end{align}
where $\AttWeights = \Z^{-1}\A$ is also known as the \emph{attention weight matrix}.
Intuitively, each token $i \in [N]$ compares its query to all keys, gets a probability distribution $\AttWeights(i,:)$ over the token positions, and based on this, takes a weighted average of the values. Evaluating~\eqref{eqn:att_matrix} via~\eqref{eqn:score_matrix} costs $\mathcal O\left (N^2 (d+\dhead)\right )$ operations in general; we refer to Section~\ref{subsec:cost_inference} for more details about the computational cost of other models. The fact that the cost is quadratic in $N$ constitutes a bottleneck for long input sequences, and this motivates the development of fast algorithms for the computation (or approximation) of~\eqref{eqn:att_matrix}, which are the subject of this overview.

\paragraph{Causality.}
In natural language processing, and in specifically next-token prediction, a causal mask $M$ is applied to the the attention scores $\A$ to prohibit the query vectors from interacting with keys from tokens occurring later in the sentence. This can be formalized by adding a mask $M_{ij}=\delta_{j\le i}$
and \[ A=\exp\left(\frac{\Q\K^{\T}}{\sqrt {\dhead}}\right)\odot M,\] where we denote by $\odot$ the Hadamard, or element-wise, product. This corresponds to removing all but the lower triangular part of the matrix $A$. For simplicity, the mechanisms described in the rest of this section are described without a causal mask.

\subsection{Multi-Headed Attention}\label{sec:multi-head}
One of the most standard generalizations of the basic model, \emph{Multi-Headed Attention (MHA)}, forms several parallel attention heads, whose outputs are further combined through a linear combination. 

Specifically, in MHA, instead of one, a total of $\Nh$ query, key $\head{\WQ}{h}, \head{\WK}{h} \in\RR^{d\times \dhead}$ and value $\head{\WV}{h}\in\RR^{d\times \dhead}$ weights are learned in parallel. Correspondingly, for all $1\le h\le \Nh$, we form
\begin{equation}
    \head{\Q}{h}=X\head{\WQ}{h},\quad \head{\K}{h}=X\head{\WK}{h},\qquad \head{\V}{h}=X\head{\WV}{h},
\end{equation}
and set

\begin{equation}
    \head{\A}{h}=\exp\left(\frac{\head{\Q}{h}\head{\K}{h}^{\T}}{\sqrt {\dhead}}\right), \quad \head{\Z}{h}(i, i)=\sum\limits_{j=1}^{N}\head{\A}{h}(i,j) \quad \text{ and }  \quad   \Atth{h}=\head{\Z}{h}^{-1}\head{\A}{h}\head{\V}{h}.
\end{equation}
We also set $\AttWeights_h=\head{\Z^{-1}}{h}\head{\A}{h}$.
The output of the complete MHA mechanism is denoted by $O\in\RR^{N\times d}$ and is obtained through the linear combination of each head output via another learned weight matrix $\WO\in\RR^{\Nh\dhead\times d}$, i.e.
\begin{equation}
    O=\begin{pmatrix}\Atth{1}|&\cdots&|\Atth{\Nh}\end{pmatrix}\WO. 
\end{equation}
In general, the parameters are set such that $\Nh\dhead = d$, maximizing the potential rank of $W^O$. 

\medskip

\paragraph{Grouped Query Attention.}
A popular and useful in practice generalization of the MHA is \emph{Grouped Query Attention (GQA)} introduced in \cite{ainslieGQATrainingGeneralized2023} which has become all but ubiquitous in modern open source models \cite{grattafioriLlama3Herd2024,teamGemma3Technical2025,yangQwen3TechnicalReport2025,openaiGptoss120bGptoss20bModel2025}. Its idea is to assign a single key-value head to a group of query heads in order to reduce the memory cost of KV caching and reduce inference-time communication overheads. 
\begin{tcolorbox}
    \paragraph{KV Caching}\label{par:KVcaching}
    In the context of causal attention, and most specifically in instances where a stream of tokens needs to be evaluated by the model, the process of KV caching consists in storing in memory the key and value matrices $K$ and $V$ such that when new tokens are added to the current text, the $K$ and $V$ rows corresponding to the previously computed tokens need not be recomputed again. For example, in a discussion with a chatbot, where new tokens are added with each question and each answer, this allows to reduce the computational cost of the language model  by not having to recompute the preceding context and only having to consider the new tokens.  
    
\end{tcolorbox}
Specifically, let $\Ng$ be the number of key/value \emph{groups} and $\Nh$ be the number of query heads ($\Ng \le \Nh$). 
We set, for any $h$ with $1\le h\le \Nh$, its corresponding group index $g_h=\left\lceil{\frac{h\cdot \Ng}{\Nh}}\right\rceil$, the index of the key-value head associated to the query head $h$.
Using the same notation and context as MHA, GQA can then be formulated as
\begin{equation}
    \head{\A}{h}=\exp\left(\frac{\head{\Q}{h}\head{\K}{g_{h}}^{\T}}{\sqrt {\dhead}}\right) \quad \text{ and } \quad \Atth{h}=\head{\Z}{h}^{-1}\head{\A}{h}\head{\V}{g_{h}}.
\end{equation}

Multi-Query Attention (MQA), an important variant of GQA, wherein all query heads share a single key-value group corresponds to the case where $\Ng=1$, minimizing the KV cache size. Also, note that the case where $\Ng=\Nh$ corresponds to regular MHA.

\begin{tcolorbox}
\paragraph{Positional Encoding}
The attention mechanism as described in this section constitutes an idealized model suited to the methods we will consider in the next sections. 
In practice, positional encodings are added to the attention models to provide more information to the model about the position of the tokens in the text. While irrelevant to most of the methods we will consider in this review, the positional encoding is fundamental to understand the Latent Attention mechanism introduced in \cite{deepseekai2024deepseekv2strongeconomicalefficient} used in several recent  models \cite{deepseekai2024deepseekv2strongeconomicalefficient,guo2025deepseek}. Rotary Positional Encoding (RoPE) and Latent Attention are described  in more detail in Section \ref{sec:MLA}.
\end{tcolorbox}

\subsection{Computational cost of inference}\label{subsec:cost_inference}
Let us discuss the computational complexity of  the most general so far GQA model, of which MQA, MHA and regular attention are subcases as we describe in the following table:
\begin{center}
\begin{tabular}{|l|l|l|l|}
\hline
\textbf{Model} & Attention & MHA & MQA  \\
\hline
Number of Query heads  & $\Nh=1$ & $\Nh=\Nh$ & $\Nh=\Nh$ \\
Number of KV heads & $\Ng=1$ & $\Ng=\Nh$ & $\Ng=1$ \\
Head dimension & $\V$: $d$, $\Q/\K$:$ \dhead$ & $\head{\Q}{h}/\head{\K}{h}/\head{\V}{h}$: $\dhead = d/\Nh$ & $\head{\Q}{h}/\head{\K}{g}/\head{\V}{g}$: $\dhead$\\
\hline
\end{tabular}
\end{center}

The computation of the $\Nh$ query matrices and $\Ng$ key/value matrices through projection from the input $X$ requires a total of $(\Nh + 2\Ng)$ matrix-matrix products, incurring a cost of $\mathcal O((\Nh + 2\Ng) N d \dhead)$. For each of the $\Nh$ query heads, the attention score calculation ($\head{\Q}{h} \head{\K}{g_h}^T$) and the subsequent application to the value matrix ($\head{\A}{h} \head{\V}{g_h}$) are both dominated by multiplications involving an $N \times N$ matrix, leading to a cost of $\mathcal O(N^2 \dhead)$ for each operation. The final projection, which combines the head outputs, involves multiplying an $N \times \Nh\dhead$ matrix by an $\Nh\dhead \times d$ matrix, at a cost of $\mathcal O(N \Nh \dhead d)$.

\begin{center}
\begin{tabular}{|l|l|}
\hline
\textbf{Operation} & \textbf{Computational Cost}  \\
\hline
Queries, Keys, Values computations  & $\mathcal O((\Nh + 2\Ng) N d \dhead)$ \\
Attention scores ($\head{\A}{h}$) & $\mathcal O(N^2 \dhead \Nh)$ \\
Attention output ($\Atth{h}$) & $\mathcal O(N^2 \dhead \Nh)$ \\
Final projection ($O$) & $\mathcal O(N \Nh \dhead d)$ \\
\hline
\end{tabular}
\end{center}

The total complexity is dominated by the $N^2$ term, resulting in a $\mathcal O( N^{2}\Nh \dhead)$ complexity. In the usual case, where $\Nh\dhead = d$ with $d$ the hidden dimension (``size'') of the model, this gives us a $\mathcal O(N^2 d)$ complexity. While this complexity is shared by all the above models, the main advantage of GQA resides in ability to provide similar performance to MHA while being more memory-efficient in its KV caching (\ref{par:KVcaching}) as can be seen in the following table.

\begin{center}
\begin{tabular}{|l|l|l|}
\hline
\textbf{Model} & \textbf{KV Cache size} &  \textbf{Computational Complexity}\\
\hline
Single-headed attention & $2Nd$ & $\mathcal O(N^2d)$\\
Multi-headed attention  & $2N\Nh\dhead$ & $\mathcal O(N^2\Nh\dhead)$ \\
Grouped query attention & $2N\Ng\dhead$ & $\mathcal O(N^2\Nh\dhead)$\\
Multi-query attention($\Ng=1$) & $2N\dhead$&$\mathcal O(N^2\Nh\dhead)$\\
\hline
\end{tabular}
\end{center}

The quadratic cost with respect to the number of tokens is a limiting factor to the size of history the model can incorporate from the computational perspective.
In practice, GQA has the same complexity as MHA, however, it allows for reduced memory cost thanks to the reduced number of cached vectors, as well as reduced communications overhead as fewer vectors need to be loaded from memory for inference. Single-headed attention, while having the same complexity and memory requirements as MHA when $\Nh\dhead = d$, suffers from reduced performance as it is not able to attend to different semantic features simultaneously.
While in practice a $N^2$ scaling can be expected given the fact that attention works on all pairs of tokens, it is believed, and experimentally verified, that the information contained in textual data is sparse  and low-rank \cite{chen2021scatterbrain}. The rest of this document is devoted to the review of several methods that attempt to reduce the $\mathcal O(N^2)$ computational cost. 

To assist the reader in navigating the variety of approaches surveyed in this work, we provide in Figure~\ref{fig:overview} a comprehensive summary of the methods discussed throughout the paper. The two main groups of algorithms consist in methods that directly approximate the attention~\eqref{eqn:att_matrix} and methods that use a related, but different, attention mechanism. Within each group, we organize techniques according to their underlying numerical principles, such as sparsity, low-rank approximation, clustering, and tensor decompositions. Beyond serving as a compact overview, this table is intended as a guiding reference: readers may find it useful to return to it when exploring individual sections, in order to contextualize each method within the broader landscape of approximation techniques for fast attention. 

\begin{tcolorbox}
    \paragraph{Neural Network Layers}\label{par:NN}
    In this document, we will consider an abstraction of an artificial neural network described as a composite function made up of blocks of operations called modules or layers. For simplicity, we shall consider the following simplified formulation of a feed-forward neural network.
    Given $\Nl\in\NN$, let $X_0\in\RR^{N\times d_\text{in}}$ be a set of $N$ input vectors to the $L$-layered neural network $f_\text{NN}$, and let, for any $\ell\in [\Nl]$, $f_\ell:\RR^{N\times d_{\ell-1}}\times \RR^{P_
    \ell}\to\RR^{N\times d_{\ell}}$ be the function associated to the $\ell$-th module of the neural network with hidden size $d_\ell\in\NN^*$ and number of parameters $P_\ell$. Then, we set
    $$f_\text{NN}(X_0;\Theta)=f_{L}(\cdot ,\Theta_N)\circ\cdots\circ f_1(\cdot ,\Theta_1)(X_0),$$
    where, for any $\ell\in[\Nl]$, $\Theta_{\ell}$ is the set of weights associated to layer $\ell$ and $\Theta=(\Theta_\ell)_{\ell\in[\Nl]}$. 
    In modern LLMs, each layer of the associated neural network model is composed of an individual attention submodule, identical to the mechanisms described above, as well as additional submodules, illustrated in Figure~\ref{fig:gpt_arch}, which, for simplicity's sake, we will not describe in detail here. Unless specified otherwise, when multi-layer methods are considered, we shall refer to the attention mechanisms associated to the $\ell$-th layer as the $\ell$-th attention layer.
    
    \centering
    \includegraphics[width=0.9\linewidth]{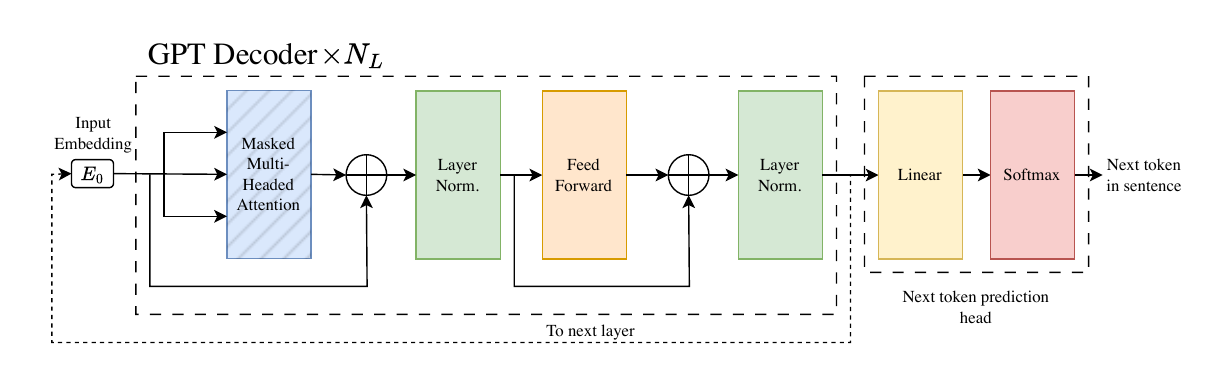}
    \captionof{figure}{GPT architecture \cite{radford2018improving} - each layer of the GPT Decoder contains multiple submodules, only one of which is a masked Multi-Headed Attention submodule.}
    \label{fig:gpt_arch}
\end{tcolorbox}

\begin{figure}
    \centering
    \vspace{-20pt}
    \input{figures/overview_approximations}
    \vspace{-0pt}
    \caption{Overview of approximation techniques for self-attention.}
    \thisfloatpagestyle{headings}
    \label{fig:overview}
\end{figure}

\section{Approximate  attention using clustering/importance sampling}\label{sec:clustering}
\begin{figure}[htb]
    \centering
    \includegraphics[width=\textwidth]{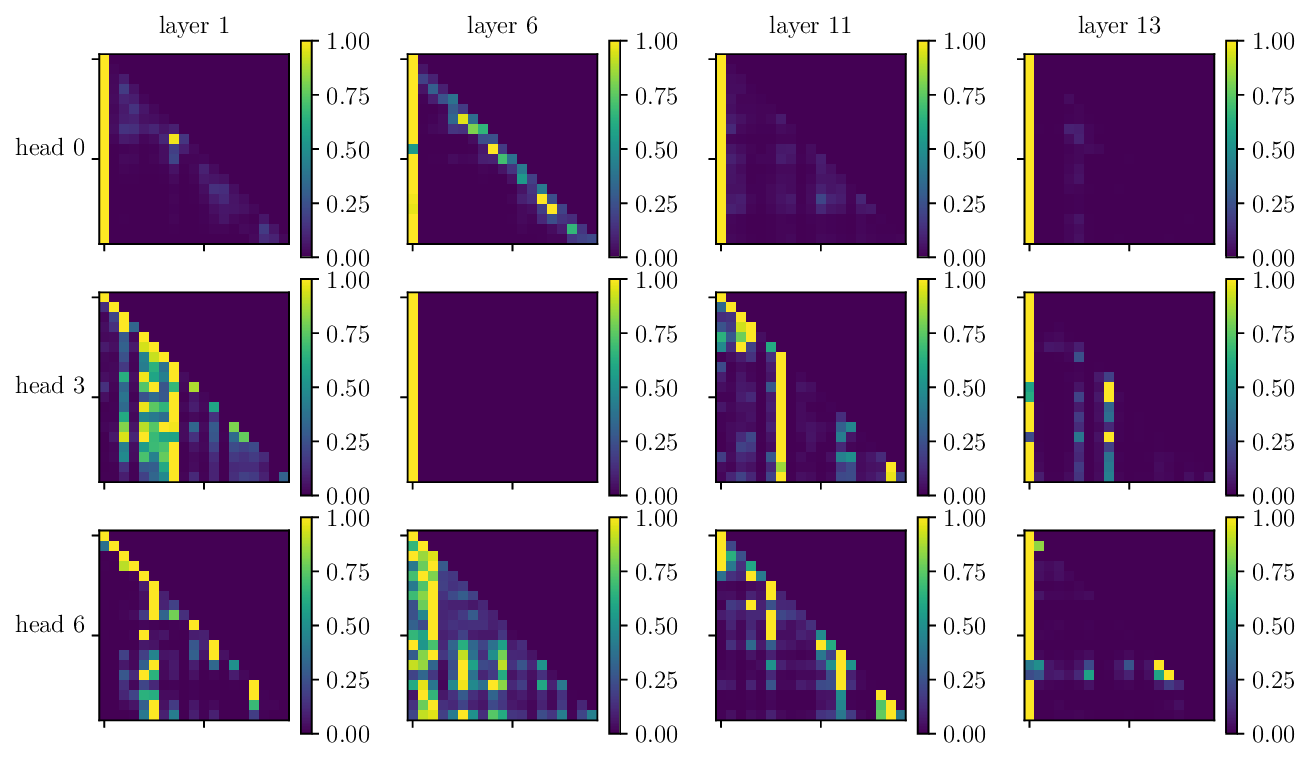}
    \caption{Approximate sparsity pattern of the top $20 \times 20$ block of the masked attention matrix $\Z^{-1}\A \in \mathbb{R}^{309 \times 309}$ corresponding to the $\K$ and $\Q$ matrices given by the Llama 3.2 (1B) model with the HyperAttention \cite{han2023hyperattention} abstract given as input, for different choices of  heads and layers. The colors represent the magnitude of the entries, with lighter colors being larger elements; each row of the matrix has been normalized to have maximum element equal to $1$. It is interesting to note that some heads in specific layers exhibit so-called \emph{attention sinks} \cite{xiaoEfficientStreamingLanguage2024a}, visible as columns with consistently large weights, indicating tokens that attract attention across many query positions.}
    \label{fig:sparsity_pattern}
\end{figure}

The methods discussed in this section are based on the observation that the softmax operation tends to accentuate the difference in magnitude between the entries of a matrix, with the consequence that, in practice, it makes sense to approximate the attention matrix $\AttWeights = \Z^{-1}\A$ by a sparse matrix. Figure~\ref{fig:sparsity_pattern} illustrates this behavior for the Llama 3.2 model. The key idea that encompasses the methods in this section is to compute only the most ``important'' entries of the matrix $\Q\K^T$ (the so-called \emph{heavy hitters}) exactly and cheaply approximate the remaining entries to obtain an approximation of the matrix $$\Z^{-1}\A = \softmax\left (\frac{\Q\K^{\T}}{\sqrt{\dhead}}\right ).$$ Once this is done, one can also cheaply approximate the matrix $\Att = \Z^{-1} A V$. This can reduce $\mathcal O(N^2)$ comparisons to $\mathcal O(N\log N)$ or $\mathcal O(N)$, where $N$ is the number of tokens.

Formally (borrowing notation from~\cite{kitaev2020reformer}\footnote{In contrast to their notation, we explicitly state the $\sqrt{d_{\text{head}}}$ factor and include masking in the definition of $\mathcal P_i$.}),  we can write the attention computation for the $i$-th query  $q_i=Q_{i,:} \in \mathbb R^{\dhead}$ as 
\begin{equation}\label{eq:o_i} o_i := (Z^{-1} A)_{i, :} V = \sum_{j \in\mathcal P_i} \exp\left( \frac{\langle q_i , k_j \rangle}{\sqrt{\dhead}} - \mathfrak{z}(i, \mathcal P_i)\right)v_j, \end{equation}
where we set $\mathcal P_i \subseteq [N]$ to be the set of indices associated to the keys that $q_i$ ``pays attention'' to, i.e. $k_j:=K_{j,:}$ for $j\in\mathcal P_i$, and the associated value vectors $v_j=V_{j,:}$,  and $$\mathfrak{z}(i, \mathcal P_i) = \log \sum_{j\in\mathcal P_i} \exp\left (\frac{\langle q_i, k_j\rangle}{\sqrt{\dhead}}\right ),$$ the logarithm of the $i$-th normalization factor associated to $q_i$, $Z_{ii}$. 

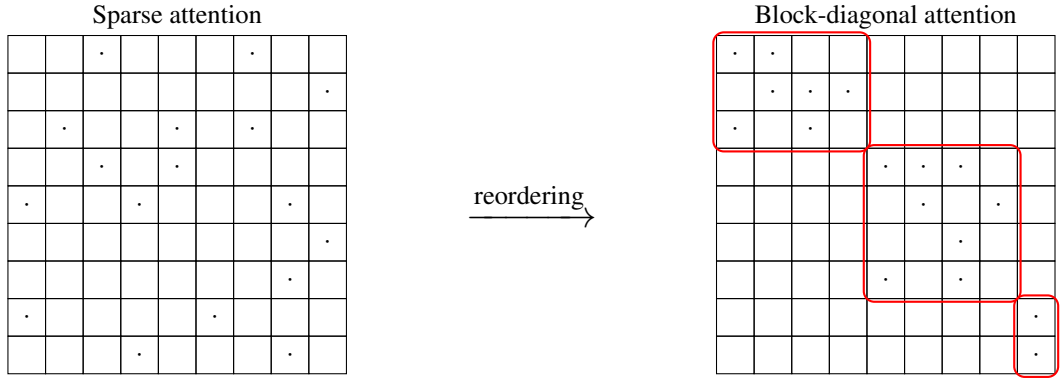
\begin{figure}
    \begin{center}
        \begin{tikzpicture}[every node/.style={inner sep=1pt},>=stealth]

        \tikzset{dot/.style={circle,fill=blue,minimum size=3mm,inner sep=0pt}},

\matrix[matrix of nodes,
    nodes in empty cells,
    nodes={minimum size=5mm, draw, anchor=center},
    column sep=-\pgflinewidth, row sep=-\pgflinewidth,
    label=above:{Sparse attention}
] (A) {
 &   &  $\cdot$ &   &  &   & $\cdot$  &   &   \\
  &   & &   &   &   &   &   & $\cdot$ \\
  & $\cdot$  &   &   &  $\cdot$ &  & $\cdot$  &   &   \\
  &   & $\cdot$  &   &  $\cdot$ &   &  &   &   \\
$\cdot$&   &   &  $\cdot$ &   &   &   &   $\cdot$ &   \\
  &   &   & &   &   &   &   & $\cdot$  \\
  &   &   &   &   &   &   &  $\cdot$ &   \\
  $\cdot$&  & &   &   & $\cdot$  &   &   &   \\
  &   &   &  $\cdot$ &   &   & & $\cdot$  &   \\
};

\node[right=1.5cm of A, align=center] (arrow) {\Large$\xrightarrow{\text{reordering}}$};

\matrix[matrix of nodes,
    nodes in empty cells,
    nodes={minimum size=5mm, draw, anchor=center},
    column sep=-\pgflinewidth, row sep=-\pgflinewidth,
    label=above:{Block-diagonal attention},
    right=1.5cm of arrow
] (B) {
$\cdot$&  $\cdot$ &   &   &   &   &   &   &   \\
  & $\cdot$ &  $\cdot$ & $\cdot$  &   &   &   &   &   \\
  $\cdot$&   & $\cdot$ &   &   &   &   &   &   \\
  &   &   &   & $\cdot$  & $\cdot$  & $\cdot$  &   &   \\
  &   &   &   &  & $\cdot$  &   &  $\cdot$ &   \\
  &   &   &   &   &  & $\cdot$  &   &   \\
  &   &   &   &  $\cdot$ &   &  $\cdot$ &   &   \\
  &   &   &   &   &   &   & &  $\cdot$ \\
  &   &   &   &   &   &   &   &$\cdot$ \\
};

\draw[thick,red,rounded corners] ($(B-1-1.north west)+(-1pt,1pt)$) rectangle ($(B-3-4.south east)+(1pt,-1pt)$);
\draw[thick,red,rounded corners] ($(B-4-5.north west)+(-1pt,1pt)$) rectangle ($(B-7-8.south east)+(1pt,-1pt)$);
\draw[thick,red,rounded corners] ($(B-8-9.north west)+(-1pt,1pt)$) rectangle ($(B-9-9.south east)+(1pt,-1pt)$);

\end{tikzpicture}
\caption{Ideally, the columns and rows of $\A$ can be reordered to produce a matrix that is close to block-diagonal. In practice, it is hard to guarantee that the blocks are comparably sized. }
\label{fig:block_diag}
\end{center}
    \end{figure}

In the normal attention setup, $\mathcal P_i = [N]$ for all $i$. In  masked attention,  $\mathcal P_i = [i] = \{j \in [n] : j \leq i\}$, as a given query only pays attention to preceding keys. Most of the techniques described in this section approximate~\cref{eq:o_i} by computing attention exactly for a much smaller subset $\mathcal P_i$ determined by importance sampling. Once these important entries are found, the keys and queries can be reordered to form an approximately block-diagonal approximation to $\AttWeights$; see Figure~\ref{fig:block_diag} for a visual example. 

The models described in this section develop different techniques for answering the following questions:

\begin{itemize}
    \item How do we efficiently find clusters for important pairs of keys and queries?
    \item How do we ensure these clusters are similarly sized? 
    \item How do we approximate the matrix $\Z$ and the remaining entries of $\Z^{-1}\A$?
\end{itemize}

We describe two categories of models in this section. Some models, such as the Reformer~\cite{kitaev2020reformer}, Routing Transformer~\cite{roy2021efficient}, SMYRF~\cite{daras2020smyrf} (and LevAttention~\cite{kannan2024levattention}, discussed in Section~\ref{subsec:levattention}), compute attention exactly for important pairs of keys and queries, then  set the rest of the attention scores to zero, or approximate them through clustering as in Multipole Attention~\cite{hooper2025multipole},  possibly reordering indices so that $Z^{-1}A$ is sparse and  close to block-diagonal (see~\cref{fig:block_diag}). Other models, such as KDEFormer~\cite{zandieh2023kdeformer} and HyperAttention~\cite{han2023hyperattention, li2025efficient}, approximate the attention scores of the remaining pairs using additional techniques. We summarize the major differences between these models in Table~\ref{table:importance_table}.

\begin{table}[htb]
\centering
\begin{tabularx}{\textwidth}{@{} l >{\raggedright\arraybackslash}X >{\raggedright\arraybackslash}X @{}}
\toprule
\textbf{Model} & \textbf{Importance sampling} & \textbf{Approximation method} \\
\midrule
Reformer
& multi-round LSH as in~\cref{eq:lsh-a} with  $Q = K$; chunking  enforces smaller groups of roughly the same size
& compute \cref{eq:o_i} exactly for  $\mathcal P_i$ defined as the set of keys in $q_i$'s chunk  and the previous chunk \\
\addlinespace
Routing
& $k$-means clustering on keys and queries
 & compute \cref{eq:o_i} exactly for $\mathcal P_i$ defined as the set of keys in $q_i$'s cluster \\
\addlinespace
SMYRF
& $H$ rounds of clustering on keys and queries using LSH as in~\eqref{eq:lsh-c}
& compute \cref{eq:o_i} exactly for $\mathcal P_{i_h}$ for $P_{i_h}$ defined as the set of keys in $q_i$'s cluster at the $h$-th round of LSH; then merge computations by~\cref{eq:merge_rounds}\\
\addlinespace
Multipole
& $k$-means clustering on keys; assign queries and values to corresponding clusters
& compute~\cref{eq:o_i} exactly for $k_j$ in/near $q_i$'s cluster; otherwise, if $k_j$ belongs to cluster $\ell$, approximates $A_{ij} \approx N_j \exp\left (\langle q_i, k_{c_j}\rangle\right ) v_{c_j}$, where $k_{c_\ell}$ and $v_{c_\ell}$ are the key and value centroids of cluster $\ell$ and $N_{\ell}$ is the number of elements of cluster $\ell$\\
\addlinespace
KDEFormer
& Kernel Density Estimation (KDE) and LSH in \eqref{eq:lsh-b}
& $Z^{-1} AV \approx \widetilde{Z}^{-1} A S ^\top \cdot S V$, where $\widetilde{Z}$ is diagonal matrix obtained from a LSH based fast Gaussian KDE and $S$ is a sampling matrix obtained from the KDE approximation \\
\addlinespace
HyperAttention
& Leverage score sampling and LSH as in~\eqref{eq:lsh-b}
& $Z^{-1} AV \approx \widetilde Z^{-1} A S^\top \cdot SV$, where $\widetilde \Z$ is formed from uniformly random rows of $K$ and LSH, and $S$ a sampling matrix based on $V$'s squared row norms \\
\midrule
\end{tabularx}
\caption{Comparison of transformers using importance sampling.}
\label{table:importance_table}
\end{table}

\subsection{Clustering and locality sensitive hashing}

Since classical clustering algorithms, such as $k$-means, are expensive to apply, many of the algorithms for approximating the attention mechanism that we consider in this section use a cheaper strategy for (approximate) clustering based on locality sensitive hashing. 

A hash function $\mathfrak h: \mathbb{R}^{\dhead} \to B$ is said to be \emph{locality sensitive} if nearby vectors, with high probability, are mapped into the same bucket $b \in B$. There are several ways to construct locality sensitive hashing (LSH) maps. Below, we illustrate different approaches, introduced in \cite{andoni2015practical,daras2020smyrf,kitaev2020reformer,zandieh2023kdeformer}, that define a hash function.

\begin{itemize}
    \item[a.] A vector $x \in \RR^{\dhead}$ is mapped into one of the $|B| = 2k$ buckets via
    \begin{align}
        \mathfrak{h}(x) \coloneqq &\, \underset{i\in [2k]}{\argmin} [\mathfrak{H}(x)]_i, \label{eq:lsh-a}\\
        \textnormal{ where} &\, \quad \mathfrak{H}(x) \coloneqq \begin{bmatrix} \langle w_1,x\rangle, & \ldots, & \langle w_{k},x\rangle, & -\langle w_1, x\rangle, & \ldots, & -\langle w_{k}, x\rangle \end{bmatrix}\notag
    \end{align}
    and $w_1, \ldots, w_{k} \in \RR^{\dhead}$ are some fixed (randomly generated) vectors.

    \item[b.] A vector $x \in \RR^{\dhead}$ is mapped into one of the $|B| = 2^{k}$ buckets via
    \begin{equation}
        \mathfrak{h}(x) \coloneqq  \begin{bmatrix} 
        \indi{\langle w_1,x\rangle > 0}, &\, \ldots, &\indi{\langle w_{k},x\rangle > 0} \label{eq:lsh-b}
        \end{bmatrix},
    \end{equation}
 where the indicator function $\indi{A} = 1$ if the event $A$ occurs and $w_1, \ldots, w_{k} \in \RR^{\dhead}$ are some fixed (randomly generated) vectors.
    \item[c.] A vector $x \in \RR^{\dhead}$ is mapped into 
    \begin{equation}\label{eq:lsh-c}
        \mathfrak{h}(x) := \left \lfloor \frac{\langle w,x\rangle + b}{r} \right \rfloor,
    \end{equation}
    where $w \in \RR^{\dhead}$ is a random vector with i.i.d. standard random Gaussian entries, $b \in \RR$ is uniformly chosen in the interval $[0, r]$, and $r \in \RR$ is a scalar parameter that controls the LSH sensitivity.  
\end{itemize}

Many models post-process the output of the LSH function to form reliable clusters. In some cases, multiple rounds of LSH are run in parallel and then combined at the end to ensure better results. For example, the Reformer model~\cite{kitaev2020reformer} performs $\nrounds$ of LSH with distinct hash functions $\{\mathfrak{h}_r(x)\}_{r = 1}^{\nrounds}$. Then, the set $\mathcal P_i$ for each $q_i$ is given by the union over these rounds:
\begin{equation}\label{eq:multi_LSH_reformer}
\mathcal P_i = \bigcup_{r = 1}^{\nrounds} \mathcal P_i^{(r)} \quad \text{where} \quad \mathcal P_i^{(r)} := \{ j : \mathfrak{h}_r(q_i) = \mathfrak{h}_r(k_j)\}.
\end{equation}

\subsection{Clustering-based approximations}

\subsubsection{Reformer }
\label{par:reformer}
The Reformer model~(\cite{kitaev2020reformer}) relies on angular LSH to assign a relatively small number of keys to each query. Reformer is an example of a ``shared-QK transformer,'' where the identical linear layers, $\WQ = \WK$, map embedded input tokens into $Q$ and $K$, enforcing that $Q = K$. The  chosen LSH scheme sorts the key and query vectors $k_i, q_j \in \RR^{\dhead}$ into $b$ buckets based on the hash function $\mathfrak{h}(x)$ defined in~\eqref{eq:lsh-a}.

Ideally, the $i$-th query $q_i$ should only attend to keys in its assigned bucket. However, as written, $\mathfrak{h}(x)$ is not guaranteed to yield hash buckets of approximately equal sizes, or even hash buckets that contain both query and key vectors. Thus, the  $i$-th key vector $k_i$ is modified and set equal to $q_i / \|q_i\|$, and the keys and queries are reordered to group the vectors according to their associated bucket. 
Within buckets, they are ordered sequentially; the resulting matrix is therefore close to block diagonal. Several independent rounds of this process are run in parallel and combined via~\cref{eq:multi_LSH_reformer} to ensure that similar vectors end up in the same bucket with high probability.

Then, we can think of Reformer as computing attention in~\cref{eq:o_i} exactly for keys in the same bucket as $q_i$, i.e.,  in the set $\mathcal P_i$ as defined in~\eqref{eq:multi_LSH_reformer}. In practice, $\mathcal P_i$ is slightly more complicated than this; after the queries and keys are reordered based on their assigned buckets, they are  grouped into smaller contiguous ``chunks'' of size $m$. The $i$-th query $q_i$  pays attention to the keys  that belong to its own chunk and the previous chunk. The chunk length $m$ is set to $2N/b$. This method can be used during training or for compressing the trained model for faster evaluation, or both. The experiments in the paper show for both training and evaluation, 2 or 4 parallel rounds of LSH are enough to get relatively high accuracy compared to full attention. Moreover, 8 parallel  rounds yield close to perfect accuracy.

\subsubsection{Routing Transformer}
\label{par:routing}
The Routing Transformer model~\cite{roy2021efficient} attempts to cluster both keys and queries and ignores all interactions between clusters.

More specifically, in the training phase, both clusters and keys are clustered using mini-batch $k$-means clustering on the same  set of centroid vectors $\mu = (\mu_1, \ldots, \mu_k)$, which are learned as well during the training. Given a fixed input sequence, the rows of the corresponding matrices $\Q$ and $\K$ are clustered in the following way: for each centroid $\mu_i$, the top-$k$ closest rows of $\K$ and the top-$k$ closest rows of $\Q$ are assigned to the $i$-th cluster. Note that, in this way, there might be rows which are not assigned to any cluster and rows which might be assigned to multiple clusters. On the other hand, all clusters have the same number of elements. Then, the attention weight matrix $\AttWeights$ is approximated as $\AttWeights \approx \softmax\left ( W \right )$, where $W \in \mathbb{R}^{N \times N}$ is the sparse matrix that coincides with ${\Q\K^{\top}}/{\sqrt{\dhead}}$ in the entries $(i,j)$ for which $i$ and $j$ are in the same cluster, and is zero in all other entries.

The number of centroids is chosen to be roughly $\sqrt{N}$ in order to balance the cost of the cluster assignments and the cost of the query/key dot products; in this way, the described strategy allows us to reduce the cost for computing the attention from $\mathcal O(N^2 \dhead)$ to $\mathcal O(N^{1.5} \dhead)$.

Essentially, both the Routing Transformer and the Reformer are computing the exact softmax of a ``virtually'' block-diagonal matrix (that is, block-diagonal up to permutation of row and column indices), where the indices in each block are chosen so that the corresponding rows of $\Q$ and $\K$ are in the same cluster. Routing transformer addresses the case in which $\Q \neq \K$ and uses mini-batch $k$-means for clustering, while Reformer uses LSH. 

\subsubsection{SMYRF }
\label{par:smyrf}
In SMYRF~\cite{daras2020smyrf}, the keys and queries are divided into balanced clusters. The strategy is to first define two functions $\varphi, \psi: \mathbb{R}^{\dhead} \to \mathbb{R}^{\dhead+2}$ such that, for any fixed key $k$, the values of $\|\varphi(\cdot) - \psi(k)\|$ maintain the same ordering as $\langle \cdot, k\rangle$, which take the form
\begin{equation*}
\varphi(q_i) := \begin{bmatrix} q_i; 0; \sqrt{M_Q^2 + M_K^2 - \|q_i\|^2} \end{bmatrix}, \quad \psi(k_j) := \begin{bmatrix} k_j; \sqrt{M_Q^2 + M_K^2 - \|k_j\|^2}; 0 \end{bmatrix},
\end{equation*}
where $M_Q := \max_i \|q_i\|^2$ and $M_K := \max_j \|q_j\|^2$ are the maximum squared norms of the queries and keys, respectively. The images $\varphi(q_i)$ and $\psi(k_j)$ are then mapped to $\RR$ using a strategy similar to the LSH function~\eqref{eq:lsh-c}: first, we map each $\varphi(q_i)$ to the real line with the function $x \mapsto \langle w, x \rangle$, for a fixed standard Gaussian random vector $w \in \RR^{\dhead+2}$; then, the resulting images on the real line are divided into $L$ consecutive buckets of equal size, corresponding to the clustering of the keys. The same goes for the keys. For a given query $q_i$, the set $\mathcal{P}_i$ is defined as the set of keys that are in the same cluster as $q_i$; these will be the keys used for the computation of the (approximate) attention.

To improve the recall of the algorithm, $\nrounds$ rounds of hashing are performed, corresponding to random vectors $w_1, \ldots, w_{\nrounds}$, and corresponding to sets of keys $\mathcal{P}_{i_1}, \ldots, \mathcal{P}_{i_{\nrounds}}$, for $i = 1, \ldots, N$. The approximation of the attention $o_i^r$ obtained in the $r$-th round of LSH is done similarly to~\cref{eq:o_i}: 
\begin{equation*}
    o_i^r := \sum_{j \in \mathcal P_{i_r}} \exp \left ( \langle q_i, k_j \rangle - \mathfrak{z}(i, \mathcal P_{i_r})\right ) v_j, \quad \mathfrak{z}(i, \mathcal P_{i_r}) := \sum_{j \in \mathcal P_{i_r}} \exp \left (\langle q_i , k_j\rangle \right).
\end{equation*}
Finally, the $\nrounds$ approximations for each query are merged in the following way:
\begin{equation}\label{eq:merge_rounds}
o_i' := \sum_{r=1}^{\nrounds} a_r o_i^r, \qquad a_r := \frac{\sum_{j \in \mathcal P_{i_r}} \exp \left (\langle q_i, k_j\rangle \right ) }{\sum_{\ell = 1}^{\nrounds} \sum_{j \in \mathcal P_{i_\ell}} \exp\left (\langle q_i , k_j\rangle \right )},
\end{equation}
which means that the final approximation is a weighted sum of the approximations resulting from each round of LSH, where the $r$-th weight is the fraction of the softmax mass that was acquired in the $r$-th round of LSH to the total mass acquired
by all rounds.

\subsubsection{Multipole Attention}
\label{par:multipole}

The Multipole Attention model~\cite{hooper2025multipole} applies $k$-means clustering to  the key vectors, producing $N_{\text{clusters}}$ clusters. For the $i$-th  cluster $c_i$,  the algorithm computes the key centroid $k_{c_i}$ by taking the mean of all the key vectors in that cluster. Then, for a given query $q$, cluster scores are computed for all key centroids, i.e., pairs $\{(q, k_{c_i})\}_{i = 1}^{N_{\text{clusters}}}$ as follows:
\begin{equation}\label{eq:multipole_score}
S_i = \frac{\exp\left (\langle q , k_{c_i}\rangle \right )}{\sum_{j \neq i} N_j \cdot \exp\left (\langle q, k_{c_j}\rangle \right )},
\end{equation}
where $N_j$ denotes the number of keys in cluster $j$. This represents the importance of the $i$-th cluster to the query $q$, relative to the rest of the clusters. Then, the clusters are sorted based on their scores $S_i$, and the exact attention computation is carried out for the keys in the highest-scoring clusters, up to some pre-determined token budget. 

The attention scores for the remaining, less important keys are approximated as follows. For all keys in the $i$-th cluster, the attention contribution is computed as
\begin{equation}
N_i \exp\left (\langle q, k_{c_i}\rangle \right )v_{c_i}, 
\end{equation}
before normalizing for the softmax function. 
Here,  $v_{c_i}$ denotes the mean of the values associated to cluster $i$. Intuitively, this amounts to attributing the same attention score $\exp\left (\langle q, k_{c_i}\rangle \right )$ to each element of the cluster. This process could also be repeated recursively in the generalized hierarchical Multipole attention model. In this case, hierarchical $k$-means clustering is applied at different scales to return a progressively coarser set of clusters and centroids. Then, at each level of the hierarchy, the Multipole attention computation is performed. This is meant to reduce the overhead cost of centroid lookup in the case of one-level Multipole attention.

A key aspect of this model is its efficient method of online clustering. With every new token, the model must update the previously computed clusters, however it is prohibitively expensive to re-cluster the entire set of keys each time. Thus, the authors of~\cite{hooper2025multipole} propose to use a sliding window clustering method, where  with each new query, $k$-means clustering is run on only  a small subset of the dataset. In addition, to ensure that the $k$-means clustering is effective, the model employs the Windowed RoPE strategy~\cite{he20252} over regular RoPE.

\subsection{Combined LSH and sampling methods}

\subsubsection{KDEformer}
\label{subsubsec:kdeformer}
The KDEformer, proposed by \cite{zandieh2023kdeformer}, 
relates the dot-product attention to the Gaussian kernel associated with a matrix $\widetilde \K \in \RR^{N \times (\dhead + 1)}$
\begin{equation*}
\mu_{\widetilde \K} : \RR^{\dhead + 1} \to \RR, \qquad \mu_{\widetilde K}(\widetilde q) := \frac{1}{N} \sum_{i=1}^N \exp \left ( - \frac{\|\widetilde q-\widetilde{k}_i\|^2}{2} \right ),
\end{equation*}
where $\widetilde{k}_i := \widetilde{\K}(i,:)$ for $i \in [N]$, which is amenable to efficient approximation using fast Gaussian Kernel Density Estimation (KDE) methods. This connection is used to construct a diagonal matrix $\widetilde Z$ that approximates $\Z$ from~\eqref{eqn:diagonal_matrix}. In this sense, KDEformer is closely related to the methods we will talk about in Section~\ref{sec:kernel-based}.
 Another key ingredient of KDEformer is the approximation of the matrix-matrix product $\A\V$ via Approximate Matrix Multiplication (AMM); that is, a suitable sampling matrix $S \in \RR^{m\times N}$ with $m = N^{1 - \Omega(1)} \ll N$ is constructed, and the output of the attention mechanism is approximated as
\begin{equation*}
    Y \approx \widetilde\Z^{-1} \A S^\top S \V.
\end{equation*}
The connection with fast KDE methods is used for the construction of the sampling probabilities that are needed to build $S$ as well. Finally, the use of LSH helps to further reduce the sample complexity. 
Below, we will sketch how the KDE is involved in the generation of $\widetilde{\Z}$ and $S$, the theoretical guaranties, and how LSH is used.

\paragraph{The Gaussian KDE trick}    

For any \emph{non-negative} vector $x = (x_1, \ldots, x_N)^\top \in \RR_+^N$, one can write 
\begin{equation}\label{eq:GaussianKDEtrick}
    \kappa(q, K, x) := \sum_{i=1}^N x_i \exp(\langle q, k_i \rangle) = N \mu_{\widetilde{\K}}(\widetilde q) \sum_{i=1}^N x_i \exp \left ( \frac{\| k_i \|^2}{2} \right ),
\end{equation}
for the enlarged vectors and matrices $\widetilde q = \begin{bmatrix} q \\ 0 \end{bmatrix} \in \RR^{\dhead+1}$ and $\widetilde{\K} := \begin{bmatrix} \widetilde k_1^T \\ \vdots \\ \widetilde k_N^T \end{bmatrix}$ with $\widetilde k_i := \begin{bmatrix} k_i \\ y_i \end{bmatrix} \in \RR^{\dhead+1}$, where $y_i \in \RR^+$ solves
\begin{equation*}
    x_i \exp(\|k_i\|^2) \exp\left ( \frac{y_i^2}{2} \right ) = \sum_{j=1}^N x_j \exp \left ( \frac{\|k_j\|^2}{2} \right ).
\end{equation*}
An approximation of~\eqref{eq:GaussianKDEtrick} for all rows of $\Q$ can be obtained by the \emph{Weighted Exponential} KDE function $\textnormal{WExpKDE}(\Q', \K', x, \epsilon)$ proposed in~\cite[]{zandieh2023kdeformer} (inspired by the theoretical guarantee for Gaussian KDE in~\cite[Theorem 2]{charikar2020kernel}), which outputs, for any \emph{non-negative} vector $x \in \RR_+^N$, any matrices $Q', K' \in \RR^{N \times \dhead}$, and any $\varepsilon \in (0,1)$, a vector $w \in \RR_+^N$ such that 
\begin{equation}\label{eq:WExpKDE}
    w_i =\left (1 + O(\varepsilon)\right ) \kappa(q_i', K', x) = \left (1 + O(\varepsilon)\right ) N \mu_{\widetilde{\K}'}(q_i') \sum_{j=1}^N x_j \exp \left ( \frac{\|k_j'\|^2}{2} \right)
\end{equation}
for all rows $q_i'$ of $\Q'$.

Note that the combination of~\eqref{eq:GaussianKDEtrick} with WExpKDE does not directly imply a formula for the product $\A\V$ since $\V$ may have negative entries. However, it can be used to estimate the normalization coefficients $\Z$ and to provide a sampling matrix proportional to the attention weights, $\Z^{-1}\A$, which will allow to apply approximate matrix multiplication methods.

\paragraph{Generating $\widetilde \Z$ and $S$. }
The diagonal elements of the exact matrix $\Z$ can be conveniently written as
\begin{equation*}
    Z_{jj} = \kappa\left (q_j \dhead^{-1/4}, K \dhead^{-1/4}, \mathbf 1_N \right ) 
\end{equation*}
for $j\in[N]$, where $\mathbf 1_N$ is the vector made by all ``ones'', so we can efficiently approximate them using WExpKDE~\eqref{eq:WExpKDE}, which will give us a vector $z \in \RR_+^N$ with non-negative entries and we set $\widetilde \Z = \mathrm{diag}(z) \approx \Z$.

The sampling matrix $S$ is obtained from a probability distribution that ensures samples are proportional to the sum of the squared column norms of $\Z^{-1} \A$ and the squared norms of $\V$. The squared row norms of $\Z^{-1} \A$ can be written as
\begin{equation*}
    \| (\Z^{-1}\A)_{i:} \|_2^2 = \kappa\left ( \sqrt{2}\dhead^{-1/4}q_i,\sqrt{2}\dhead^{-1/4}\K, w  \right ), \qquad w = \begin{bmatrix} z_1^{-2} & \cdots & z_N^{-2} \end{bmatrix}^\top,
\end{equation*}
so they can be again approximated by a vector $\bar{w} \in \RR_+^N$ by WExpKDE. 
At this point, for each $i \in [N]$, we define 
\begin{equation*}
\widetilde{p}_{i} = \bar{w}_i + \frac{\|v_{i}\|^{2}}{\|\V\|^{2}_{\textnormal{op}}}
\end{equation*}
with $v_{i}\in \RR^{d}$ denoting the $i$-th row of $\V$.
Let $[p_{\ell}]_{\ell\in [N]}$ denote the normalized version of $[\widetilde{p}_{\ell}]_{\ell\in [N]}$, creating the probability distribution we were seeking. 
Given $m$ i.i.d. samples $\ell_{1}, \ldots, \ell_{m}\in [N]$ from the distribution $[p_{\ell}]_{\ell\in [N]}$, we define the $i$-th row of $S$ as $\frac{1}{\sqrt{mp_{\ell_{i}}}} e_{\ell_{i}}^{\top}$, where $e_{\ell_{i}}$ is the $\ell_i$-th vector of the canonical basis of $\RR^N$. The number of samples needed is $$m = \Omega\left (\epsilon^{-2}\log(N)\cdot (\mathrm{srank}(\Z^{-1}\A) + \mathrm{srank}(\V))\right ),$$ where $\mathrm{srank}(X) = \|X\|_{\mathrm{F}}^{2}/\|X\|_{\textnormal{op}}^{2}$ denotes the \emph{stable rank} of the matrix $X$.

\paragraph{Theoretical guarantees. }  
It is shown in~\cite[Theorem 3.4]{zandieh2023kdeformer} that is it possible to choose the parameters of KDEformer such that, for any $\varepsilon > 0$, the algorithm has a theoretical computational complexity  of $\mathcal{O}(\epsilon^{-2}d\cdot N^{1.173 + o(1)})$ and produces an approximation such that
\begin{equation}
    \|\Att - \widetilde{\Z}^{-1} \A \, S^{\top} \, S \V \|_{\textnormal{op}} \leq \epsilon \cdot \|\Z^{-1}\A\|_{\textnormal{op}} \cdot \|\V\|_{\textnormal{op}}, \label{eqn:KDE_approxmate_attention}
\end{equation}
with probability at least $1 - 1/\textnormal{poly}(N)$.

\paragraph{Reducing sample complexity via LSH.} 
To reduce the number of samples $m$, \cite{zandieh2023kdeformer} propose a practical technique for reducing the stable rank of $\Z^{-1}\A$ by finding and subtracting its ``heavy'' elements via LSH. Specifically, given the LSH function $\mathcal{H}: \RR^{\dhead} \mapsto [B]$ in~\eqref{eq:lsh-b}, they define the sparse approximation $\A_{\textnormal{spar}}$ below, which contains the dominant entries of $\A$
\[
    [\A_{\textnormal{spar}}]_{ij} = \exp\left ({\dhead}^{-1/2}\< q_{i}, k_{j}\>  \right )\,\cdot \indi{\mathcal{H}(q_{i}) = \mathcal{H}(k_{j})}, \quad \text{ for all } i, j\in [N].
\]
Define the residual attention matrix by $\A_{\textnormal{res}} = \A - \A_{\textnormal{spar}}$. We then apply the previous construction of $\widetilde{Z}$ and $S_\text{res}$ based on $\A_{\textnormal{res}}$. The outcome of the revised algorithm is then given by
\[
\widetilde{\Att} = \widetilde{Z}^{-1} \A_{\textnormal{spar}}V + \widetilde{Z}^{-1} \A_\text{res}S^\sT_\text{res}S_\text{res}V.
\]
This improves the performance of the AMM, as the stable rank of the residual matrix is lower than that of the full attention score matrix.

\subsubsection{HyperAttention}
\label{subsubsec:hyperattention}
HyperAttention, proposed in~\cite{han2023hyperattention}, uses LSH to find the ``most important'' entries of the product $\Q\K^{\top}$ and combines this with random sampling to compute an approximation $\widetilde \Z \approx \Z$. Then, the attention output is approximated by matrix sampling similarly to KDEformer, but with a different sampling distribution. We now explain the main ingredients of HyperAttention in more detail.

First, the LSH method described in~\eqref{eq:lsh-b} is run on the rows of $Q$ and $K$. The buckets, which may contain different numbers of indices, are then arranged in a Hamming-sorted order, and two permutations $\pi_{\Q}$ and $\pi_{\K}$ are defined so that $\pi_{\Q}(i) < \pi_{\Q}(j)$ whenever the bucket corresponding to row $q_i$ comes before the bucket corresponding to row $q_j$, and similarly for $\pi_{\K}$. A block size $b$ is fixed and the clusters are re-defined to have size $b$, that is, rows $q_i$ and $k_j$ are in the same cluster if and only if $\left \lfloor {\pi_{\Q}(i)}/{b} \right \rfloor = \left \lfloor{\pi_{\K}(j)}/{b} \right \rfloor$. Moreover, HyperAttention chooses a small number $\ell$ of indices $j_1, \ldots, j_\ell \in [N]$ uniformly at random. Now, an approximation $\widetilde{\A}$ to $\A$ is defined, entrywise, as
\begin{equation*}
\widetilde{A}_{ij} = \begin{cases}
    \exp\left (\frac{\langle q_i, k_j \rangle}{\sqrt{\dhead}} \right ) & \text{ if } \left \lfloor {\pi_{\Q}(i)}/{b} \right \rfloor = \left \lfloor{\pi_{\K}(j)}/{b} \right \rfloor \text{ or } j \in \{j_1,\ldots,j_{\ell}\};\\
    0 & \text{ otherwise}.
\end{cases}
\end{equation*}
This means that we compute the exact attention scores for the pairs of keys and queries belonging to the same cluster and for the randomly selected keys. 
The motivation for the Hamming-sorted order is that the rows end up in buckets which are hopefully close to their original bucket. 
The diagonal normalization matrix $\Z$ is approximated by summing the rows of $\widetilde{\A}$, which results in a matrix denoted by $\widetilde \Z$. 

Finally, the product $\Z^{-1}\A \V$ is approximated as 
\begin{equation}
\widetilde{\Z}^{-1} \A S^{\top} \cdot S \V,\label{eq:hyperattn}
\end{equation}
where $S$ is a sampling matrix that subselects some rows of $\V$ and the corresponding columns of $\A$; the sampling probabilities are given by the squared row norms of $\V$. 
Note that, to evaluate~\eqref{eq:hyperattn}, it is only necessary to compute the columns of $\A$ corresponding to the selected rows of $\V$; while this means that the information contained in $\widetilde A$ is not enough, it is still computationally efficient if the size of the matrix $S$ is fairly small. 

The main shortcomings of HyperAttention are that (1) LSH may fail to identify \emph{all} blocks where $\Q\K^{\T}$ has large entries, and (2) approximating matrix-vector multiplication as in~\eqref{eq:hyperattn} is pretty slow to converge because it is doing Monte Carlo sampling. 
An attempt to address the first problem was made in Prescoring, recently proposed in \cite{li2025efficient}: here, $s\in \mathbb{N}$ keys (corresponding to rows of $\K$) are pre-selected either using $k$-means clustering or by using leverage scores, inspired by the LevAttention method described in \cite{kannan2024levattention} (which we review in Section~\ref{subsec:levattention}). Then, HyperAttention is applied to the whole matrix $\Q$ and only the selected rows of $\K$ and the corresponding rows of $\V$.

\section{Low-rank techniques}
\label{sec:low-rank-approximation}

\begin{figure}[!ht]
    \centering
    \scalebox{0.9}{\input{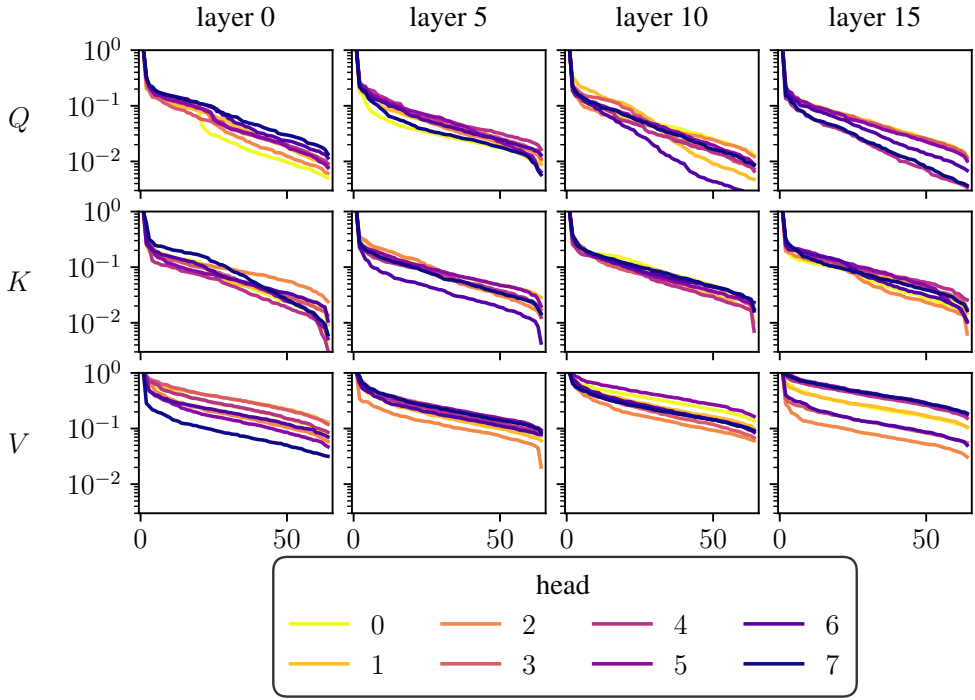}}    
    \caption{Singular values of the query, key, and value matrices involved in the attention computation for a subset of the layers and all attention heads. We use the Llama 3.2 (1B) model in float64 precision with the abstract of \cite{han2023hyperattention} as input to demonstrate the singular value decays which are observed in the literature.}
    \label{fig:singular_values_qkv}
\end{figure}

\begin{figure}[!ht]
    \centering
    \scalebox{0.9}{\input{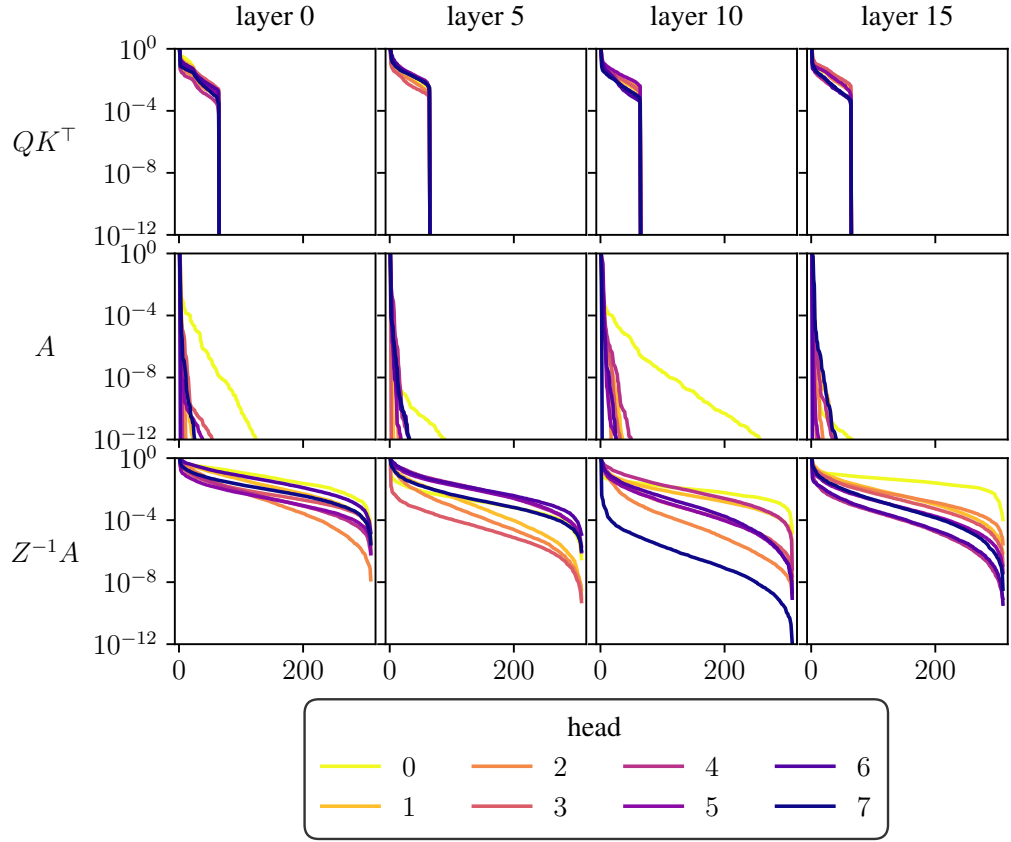}}
    \caption{Singular values of different matrices involved in the attention computation for a subset of the layers and all attention heads. We use the Llama 3.2 (1B) model in float64 precision with the abstract of \cite{han2023hyperattention} as input to demonstrate the singular value decays which are observed in the literature.}
    \label{fig:singular_values_dotproduct}
\end{figure}

The matrices involved in the computation of the attention are often numerically low-rank~\cite{chen2021scatterbrain}. To illustrate this, we show in  Figure~\ref{fig:singular_values_qkv} the singular value decays of $\head\Q h, \head \K h , \head \V h \in \RR^{N\times \dhead}$, and in Figure~\ref{fig:singular_values_dotproduct} the singular value decay of $\head \Q h \head \K h ^\top, \head \A h,$  and $\head \AttWeights h  := \head \Z h ^{-1}\head \A h\in\RR^{N\times N}$, for a selection of heads $h\in\NN^*$ of the Llama 3.2 model\cite{grattafioriLlama3Herd2024} and for a specific choice of input. The experimental settings for these figures are the same as in Figure~\ref{fig:sparsity_pattern}.
While $\head \Q h$ and $\head \K h$, and therefore $\head\Q h\head \K h^\top$, are low-rank by definition, as $\dhead\ll N$, their $\dhead$ singular values do show a moderate decay. More interestingly, a strong decay of the singular values can be observed for all layers and heads of $\head \A h$, and a weaker decay for $\head \Z h ^{-1} \head \A h$. 

This observation opens the door to the approximation of the attention mechanism by certain factorized representations at the benefit of higher computational efficiency. In this section, we review a small collection of methods that either use projections of the matrices $\Q,\K$, and $\V$ before applying the softmax function or directly perform low-rank approximations of the matrix $\A$. It is important to note that some of the methods we consider (Loki in Section~\ref{sec:lowrankK} and Skyformer in Section~\ref{sec:nystromformer})  
aim to compress a pretrained model. In contrast, the others are entirely new mechanisms that must be trained separately. These alternative attention mechanisms are often motivated by their ability to achieve approximately the same result as the standard attention mechanism, while being faster to evaluate. 

\subsection{Methods for compressing self-attention}

\subsubsection{Loki}\label{sec:lowrankK}

Loki~\cite{singhania2406loki} is based on the empirical observation that the matrix $K$ representing the keys usually has an ``effective rank'' which is lower than $d$, its number of  columns. In particular, the metric they consider is the rank at which $90\%$ of the variance is explained by the first $r$ principal components. More precisely, for each layer $\ell$ and each head $h$, they look at the minimum integer $d$ such that $\sum_{i=1}^d \lambda_{i,h,\ell} \ge 0.9$, where $\lambda_{i,h,\ell}$ is the $i$-th normalized eigenvalue of the covariance matrix of the keys at layer $\ell$ and head $h$. \cite[Figure 1]{singhania2406loki} shows that for several models, such as Llama2-7B,
Llama3-70B, Mixtral-8x7B, and Phi3-Mini-4K, the average rank that explains $90\%$ of the variance (averaged over all heads and layers) is around $80$,  while $d = 128$. 
Although the this ``effective'' rank is roughly $60\%$ of the original rank, this still allows for some compression and fast computations.

    The initialization step of Loki consists of  computation of a matrix $P_r \in \mathbb{R}^{\dhead \times r}$, with $r = \dhead/4$ or $r=\dhead/8$ in their experiments, 
    that stores the principal components of the key matrix generated from a calibration dataset, such as BookCorpus, C4, and Wikitext (which contain millions, or even trillions, of tokens). 
    The influence of the particular dataset used seems to be limited, as low-rank behavior has been observed across different calibration datasets and the generalizability of the principal matrix has been explored in~\cite[Section 6.3]{singhania2406loki}. Once the matrix $P_r$ has been computed, Loki proceeds in three steps: 
\begin{enumerate}
    \item The attention scores for the query $q_i$ are approximated as 
    \begin{equation*}
    \softmax\left ( \frac{q_i^T P_r (\K P_r)^T}{\sqrt{\dhead}}\right ) V.
    \end{equation*}

    \item The top-$k$ keys are selected as those with the highest approximate attention scores. In their numerical experiments, they choose $k = N/8$ or $k = N/4$. 

    \item The attention scores are recomputed exactly for the top-$k$ keys.
\end{enumerate}

Assuming that a suitable matrix $P_r$ has been already computed offline, and $\K P_r$ is available in $\mathcal O(N \dhead r)$ operations, computing the approximate attention scores (for one value of $i$) costs $\mathcal O(rN)$. The selection of the top-$k$ keys costs $\mathcal O(N \log N)$.
Finally, recomputing the exact attention score for the top-$k$ keys costs $\mathcal O(\dhead k)$. The whole procedure allows us to pass from the $\mathcal O(N^2 \dhead)$ complexity of the standard computation of the attention to $\mathcal O(N^2 r)$; the cost still grows quadratically with the length of the sequence, and the memory requirement of the KV-cache is the same as the original attention, but the multiplicative constant corresponding to the hidden dimension is reduced.  Loki does not require retraining the network and does not require fine-tuning.

\subsubsection{Skyformer and WILDCAT}
\label{par:skyformer}

Empirically, the attention score matrix $\A$ has a relatively fast singular value decay, as can  be observed in \cref{fig:singular_values_dotproduct}. Consequently, it can be well approximated by a low-rank factorization $\widehat{\A} = U W$. Clearly, the low-rank approximated attention $\widehat{\Att} = \diag(\widehat{\A} \boldsymbol{1}_N)^{-1} \widehat{\A} V = \diag(U W \boldsymbol{1}_N)^{-1} U W V $ can be computed faster and using less memory by first computing $W \boldsymbol{1}_N$ and $WV$ before applying the outputs to $U$. However, the element-wise exponential in the definition of $\A$ complicates the computation of an approximation $\widehat{\A}$ significantly.

The authors of the Skyformer paper \cite{chen2021skyformer} notice that the attention scores matrix can be written as
\begin{equation*}
    \A = \exp\left(\frac{\Q \K^{\top}}{\sqrt{\dhead}}\right) = \kappa_{\exp}(\Q, \K)  \in \mathbb{R}^{N \times N},
\end{equation*}
where we extend the notation used in \eqref{eqn:score_matrix} such that $\kappa_{\exp}(\Q, \K) = (\kappa_{\exp}(q_i,k_j))_{1\le i,j\le N}$. Then, the symmetric matrix
\begin{equation*}
    B = \kappa_{\exp}\left(\begin{bmatrix} \Q \\ \K \end{bmatrix}, \begin{bmatrix} \Q \\ \K \end{bmatrix} \right) = \begin{bmatrix} \kappa_{\exp}(\Q, \Q) & \A \\ \A^{\top} & \kappa_{\exp}(\K, \K) \end{bmatrix}  \in \mathbb{R}^{2N \times 2N}
\end{equation*}
is a positive semidefinite (PSD) kernel matrix. 

There exist many effective methods for computing low-rank approximations of PSD kernel matrices from column or row samples. The obvious candidates are column-sampling methods such as the Nyström approximation, which is used by the Skyformer model \cite{chen2021skyformer},  as well as incomplete Cholesky factorizations, for example using randomized pivoting, as they are used by the WILDCAT model \cite{schroeder2026wildcatnearlinearattentiontheory,, chen2025rpcholesky}. Given any approximation $\widehat{B}$ to $B$, it is straightforward to verify that the upper-right $N \times N$ block $\widehat{\A} = (I, 0) \widehat{B} (0, I)^{\top}$ of $\widehat{B}$ satisfies $\lVert \A - \widehat{\A} \rVert \leq \lVert B - \widehat{B} \rVert$ for both the spectral and Frobenius norm. Additionally, \cite[Lemma 1]{schroeder2026wildcatnearlinearattentiontheory} derives an upper bound on the error of the approximated attention $\widehat{\Att}$ in terms of the approximation error of $\widehat{\A}$ to $\A$. 

\begin{figure}[!ht]
    \centering
    \begin{tikzpicture}
    \fill[gray!10] (2.5, 2.25) rectangle (5.5, -0.75);
    \draw[line width=1.3, draw=darkorange, fill=lightorange] (4, 0.75) rectangle (5.5, 2.25) node[pos=0.5] {$\A$};
    \draw[line width=1.3, draw=black!30] (2.5, -0.75) rectangle (4, 0.75) node[pos=0.5] {$\A^{\top}$};
    \draw[line width=1.3] (2.5, 2.25) rectangle (5.5, -0.75);
    \node at (6, 0.75) {$\approx$};
    \fill[gray!10] (6.5, 2.25) rectangle (7, -0.75);
    \draw[line width=1.3, draw=darkgreen, fill=lightgreen] (6.5, 2.25) rectangle (7, 0.75);
    \draw[line width=1.3] (6.5, 2.25) rectangle (7, -0.75);
    \draw[line width=1.3, draw=black, fill=lightgreen] (7.25, 0.5) rectangle (7.75, 1);
    \fill[gray!10] (8., 0.5) rectangle (11, 1);
    \draw[line width=1.3, draw=darkgreen, fill=lightgreen] (11., 0.5) rectangle (9.5, 1);
    \draw[line width=1.3] (8., 0.5) rectangle (11, 1);

    \draw [->] (11.5, 0.75) to (12, 0.75);
    \draw[line width=1.3, draw=darkorange, fill=lightorange] (12.5, 0) rectangle (14, 1.5) node[pos=0.5] {$\A$};
    \node at (14.5, 0.75) {$\approx$};
    \draw[line width=1.3, draw=darkgreen, fill=lightgreen] (15, 0) rectangle (15.5, 1.5);
    \draw[line width=1.3, draw=darkgreen, fill=lightgreen] (15.75, 0.5) rectangle (16.25, 1);
    \draw[line width=1.3, draw=darkgreen, fill=lightgreen] (16.5, 0.5) rectangle (18, 1);
\end{tikzpicture}
    \caption{Sketch of the symmetrization and approximation procedure.}
    \label{fig:symmetrization}
\end{figure}
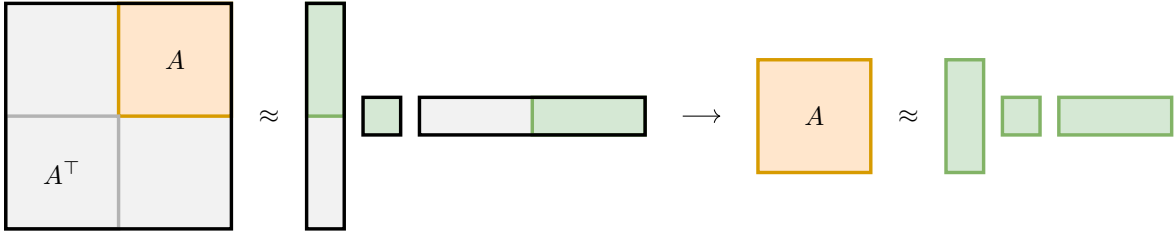

Nevertheless, constructing a good low-rank approximation $\widehat{B}$ to $B$ is more nuanced than it seems. Even when $\A$ exhibits rapid singular value decay -- allowing it to be well approximated by a low-rank factorization -- this decay typically does not carry over to $B$. Thus, achieving a good approximation quality often requires a factorization with a significantly higher rank, which in turn makes the factorization more expensive to compute and store in memory. In \cite{schroeder2026wildcatnearlinearattentiontheory}, the fact that the attention scores $\A$ are invariant under the transformation $\Q \mapsto \tau \Q$ and $\K \mapsto \tau^{-1} \K$ for some $\tau > 0$ is exploited to enforce low-rankness on $\kappa_{\exp}(K, K)$, which is the matrix from which the pivots in their randomly pivoted Nyström method are selected. They provide a closed form value for $\tau$ which balances the trade-off between low-rankness in $\kappa_{\exp}(K, K)$ and large entries in $\kappa_{\exp}(Q, Q)$.

The authors of \cite{chen2021skyformer} notice that in practice, an isolated approximation to $\A$ (without the normalization $\Z^{-1}$) is unstable and prone to floating-point overflow. This is because forming the entries of $\A$ requires exponentiating inner products between rows of $Q$ and $K$, which may already be large. Exponentiating these values quickly leads to numerical overflow. The standard attention mechanism avoids this issue by never explicitly forming $\A$, but instead only ever forming the normalized attention weight matrix $\AttWeights = \Z^{-1} \A$. To prevent this instability, Skyformer instead uses an alternative attention mechanism based on the Gaussian kernel $\kappa_{\mathrm{Gauss}}(x, y) = \exp\left (-\frac{\lVert x - y \rVert^2}{2\sqrt{\dhead}}\right )$, for which the entries in $\A$ are, by design, limited to $[0, 1]$, therefore, making a normalization with $\Z^{-1}$ obsolete \cite{chen2021skyformer}.

\subsection{Alternative self-attention mechanisms}

\subsubsection{Linformer}\label{subsec:linformer}
The \emph{Linformer} attention mechanism was proposed in \cite{wang2020linformer} and was one of the first works to suggest a method for reducing the complexity of self-attention, which typically scales linearly in the number of tokens $N$. In particular, they introduce projection matrices $P_{\K}, P_{\V} \in \mathbb{R}^{N \times k}$ for a projected dimension $k \ll N$. The entries in these matrices are parameters which are learned during training. The projections compress the attention mechanism to
\begin{equation}
   \softmax\left(\frac{\Q (P_{\K} \K)^{\top}}{\sqrt{\dhead}}\right) (P_{\V} \V).
   \label{equ:linformer-attention}
\end{equation}
The smaller the projected dimension $k$ is, the smaller is the memory and time complexity for computing the self-attention. 

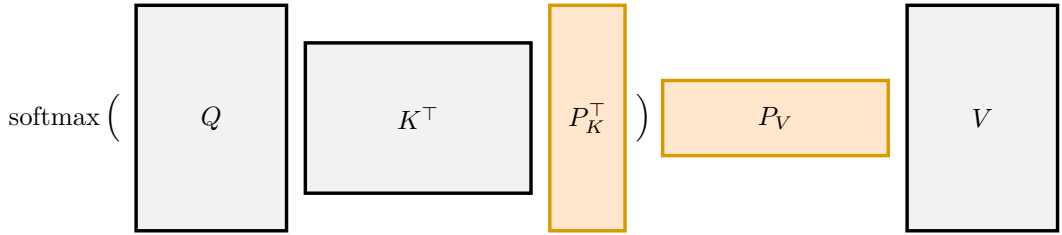
\begin{figure}[!ht]
    \centering
    \begin{tikzpicture}
    \node at (-0.95, 1.5) {$\operatorname{softmax}\Big($};
    \draw[line width=1.3, fill=gray!10] (0, 0) rectangle (2, 3) node[pos=.5] {$\Q$};
    \draw[line width=1.3, fill=gray!10] (2.25, 0.5) rectangle (5.25, 2.5) node[pos=.5] {$\K^{\top}$};
    \draw[line width=1.3, draw=darkorange, fill=lightorange] (5.5, 0) rectangle (6.5, 3) node[pos=.5] {$P_\K^{\top}$};
    \draw[line width=1.3, draw=darkorange, fill=lightorange] (7, 1) rectangle (10, 2) node[pos=.5] {$P_\V$};
    
    \node at (6.75, 1.5) {$\Big)$};
    \draw[line width=1.3, fill=gray!10] (10.25, 0) rectangle (12.25, 3) node[pos=.5] {$\V$};

\end{tikzpicture}
    \caption{Schematic depiction of the Linformer mechanism.}
    \label{fig:linformer}
\end{figure}
While the model has been shown to work empirically, the assumptions for the theoretical justifications are very restrictive. A similar, more recent, randomization-based approximation method with better theoretical guarantees is the Performer architecture, which we explore in section~\ref{subsec:performer}.

\subsubsection{Nyströmformer}\label{sec:nystromformer}
In Nyströmformer \cite{xiong2021nystromformer}, the normalised attention matrix $\AttWeights$ is approximated by considering the projection of both $\Q$ and $\K$ onto a subset of landmark query and key vectors, i.e. the set of vectors from $\Q$ and $\K$ which are the  most important for the reconstruction of $\Z^{-1}\A$. 
In this case, the projection is combined with a Nyström-like (as the softmax function is not a kernel) approximation of $\AttWeights=\Z^{-1}\A$ of the form

\begin{equation}
\softmax\left(\frac{\Q (P_{\K} \K)^{\top}}{\sqrt{\dhead}}\right) \left[\softmax\left(\frac{(P_{\Q} \Q) (P_{\K} \K)^{\top}}{\sqrt{\dhead}}\right)\right]^{\dagger} \softmax\left(\frac{(P_{\Q} \Q) \K^{\top}}{\sqrt{\dhead}}\right) \V,
    \label{equ:nystromformer}
\end{equation}
where $\dagger$ denotes the pseudoinverse of a matrix. The projection matrices $P_\Q$ and $P_\K$ can be chosen in different ways. In the Nyströmformer, they are computed via the segment-means, i.e., by taking the average over a certain number of adjacent rows in $Q$ and $K$. They empirically observe that dividing the rows into $64$ segments is often sufficient to ensure a good approximation. Further, the Nyströmformer uses an iterative procedure for approximating the pseudoinverse of $C=\softmax\left(\frac{(P_{\Q} \Q) (P_{\K} \K)^{\top}}{\sqrt{\dhead}}\right)$. If only a few landmark queries and keys are enough for a good approximation, this factorization speeds up the attention computation and reduces the required memory for storing the weights.

\begin{figure}[!ht]
    \centering
    \begin{tikzpicture}
    \node at (-0.95, 1.5) {$\operatorname{softmax}\Big($};
    \draw[fill=gray!10, line width=1.3] (0, 0) rectangle (2, 3) node[pos=.5] {$\Q$};
    \fill[lightorange!75] (2.25, 0.5) rectangle (2.625, 2.5);
    \fill[lightyellow!75] (2.625, 0.5) rectangle (3, 2.5);
    \fill[lightgreen!75] (3, 0.5) rectangle (3.375, 2.5);
    \fill[lightblue!75] (3.375, 0.5) rectangle (3.75, 2.5);
    \draw[line width=1.3] (2.25, 0.5) rectangle (3.75, 2.5) node[pos=.5] {$(P_\K \K)^{\top}$};
    \fill[lightorange!50] (1.5, 3.5) rectangle (1.75, 5.5);
    \fill[lightorange!75] (1.75, 3.5) rectangle (2, 5.5);
    \fill[lightorange!100] (2, 3.5) rectangle (2.25, 5.5);
    \fill[lightyellow!50] (2.25, 3.5) rectangle (2.5, 5.5);
    \fill[lightyellow!75] (2.5, 3.5) rectangle (2.75, 5.5);
    \fill[lightyellow!100] (2.75, 3.5) rectangle (3, 5.5);
    \fill[lightgreen!50] (3, 3.5) rectangle (3.25, 5.5);
    \fill[lightgreen!75] (3.25, 3.5) rectangle (3.5, 5.5);
    \fill[lightgreen!100] (3.5, 3.5) rectangle (3.75, 5.5);
    \fill[lightblue!50] (3.75, 3.5) rectangle (4, 5.5);
    \fill[lightblue!75] (4, 3.5) rectangle (4.25, 5.5);
    \fill[lightblue!100] (4.25, 3.5) rectangle (4.5, 5.5);
    \draw [decorate, decoration={brace, amplitude=5pt, mirror}, line width=0.5] (1.5, 3.45) to (2.25, 3.45);
    \draw[->] (1.875, 3.3) to (1.875, 3.1) to (2.425, 3.1) to (2.425, 2.6);
    \draw [decorate, decoration={brace, amplitude=5pt, mirror}, line width=0.5] (2.25, 3.45) to (3, 3.45);
    \draw[->] (2.625, 3.3) to (2.625, 3.1) to (2.8, 3.1) to (2.8, 2.6);
    \draw [decorate, decoration={brace, amplitude=5pt, mirror}, line width=0.5] (3, 3.45) to (3.75, 3.45);
    \draw[->] (3.375, 3.3) to (3.375, 3.1) to (3.175, 3.1) to (3.175, 2.6);
    \draw [decorate, decoration={brace, amplitude=5pt, mirror}, line width=0.5] (3.75, 3.45) to (4.5, 3.45);
    \draw[->] (4.125, 3.3) to (4.125, 3.1) to (3.55, 3.1) to (3.55, 2.6);
    \draw[line width=1.3] (1.5, 3.5) rectangle (4.5, 5.5) node[pos=.5] {$\K^{\top}$};
    \node at (3.95, 1.5) {$\Big)$};
    
    \draw[line width=1.3, fill=gray!10] (4.15, 0.75) rectangle (5.6, 2.25) node[pos=.5] {$C^{\dagger}$};
    
    \node at (6.6, 1.5) {$\operatorname{softmax}\Big($};
    \fill[lightorange!75] (8.25, 1.875) rectangle (10.25, 2.25);
    \fill[lightyellow!75] (8.25, 1.5) rectangle (10.25, 1.875);
    \fill[lightgreen!75] (8.25, 1.125) rectangle (10.25, 1.5);
    \fill[lightblue!75] (8.25, 0.75) rectangle (10.25, 1.125);
    \draw[line width=1.3] (8.25, 0.75) rectangle (10.25, 2.25) node[pos=.5] {$P_\Q \Q$};

    \fill[lightorange!50] (8.25, 6) rectangle (10.25, 5.75);
    \fill[lightorange!75] (8.25, 5.75) rectangle (10.25, 5.5);
    \fill[lightorange!100] (8.25, 5.5) rectangle (10.25, 5.25);
    \fill[lightyellow!50] (8.25, 5.25) rectangle (10.25, 5);
    \fill[lightyellow!75] (8.25, 5) rectangle (10.25, 4.75);
    \fill[lightyellow!100] (8.25, 4.75) rectangle (10.25, 4.5);
    \fill[lightgreen!50] (8.25, 4.5) rectangle (10.25, 4.25);
    \fill[lightgreen!75] (8.25, 4.25) rectangle (10.25, 4);
    \fill[lightgreen!100] (8.25, 4) rectangle (10.25, 3.75);
    \fill[lightblue!50] (8.25, 3.75) rectangle (10.25, 3.5);
    \fill[lightblue!75] (8.25, 3.5) rectangle (10.25, 3.25);
    \fill[lightblue!100] (8.25, 3.25) rectangle (10.25, 3);
    \draw[line width=1.3] (8.25, 3) rectangle (10.25, 6) node[pos=.5] {$\Q$};
    \draw [decorate, decoration={brace, amplitude=5pt, mirror}, line width=0.5] (8.2, 6) to (8.2, 5.25);
    \draw[->] (8.05, 5.625) to (7.35, 5.625) to (7.35, 2.05) to (8.2, 2.05);
    \draw [decorate, decoration={brace, amplitude=5pt, mirror}, line width=0.5] (8.2, 5.25) to (8.2, 4.5);
    \draw[->] (8.05, 4.875) to (7.5, 4.875) to (7.5, 1.7) to (8.2, 1.7);
    \draw [decorate, decoration={brace, amplitude=5pt, mirror}, line width=0.5] (8.2, 4.5) to (8.2, 3.75);
    \draw[->] (8.05, 4.125) to (7.65, 4.125) to (7.65, 1.35) to (8.2, 1.35);
    \draw [decorate, decoration={brace, amplitude=5pt, mirror}, line width=0.5] (8.2, 3.75) to (8.2, 3);
    \draw[->] (8.05, 3.375) to (7.8, 3.375) to (7.8, 1) to (8.2, 1);
    
    \draw[line width=1.3, fill=gray!10] (10.5, 0.5) rectangle (13.5, 2.5) node[pos=.5] {$\K^{\top}$};
    \node at (13.7, 1.5) {$\Big)$};
\end{tikzpicture}
    \caption{Schematic depiction of the Nyströmformer approximation.}
    \label{fig:nystromformer}
\end{figure}
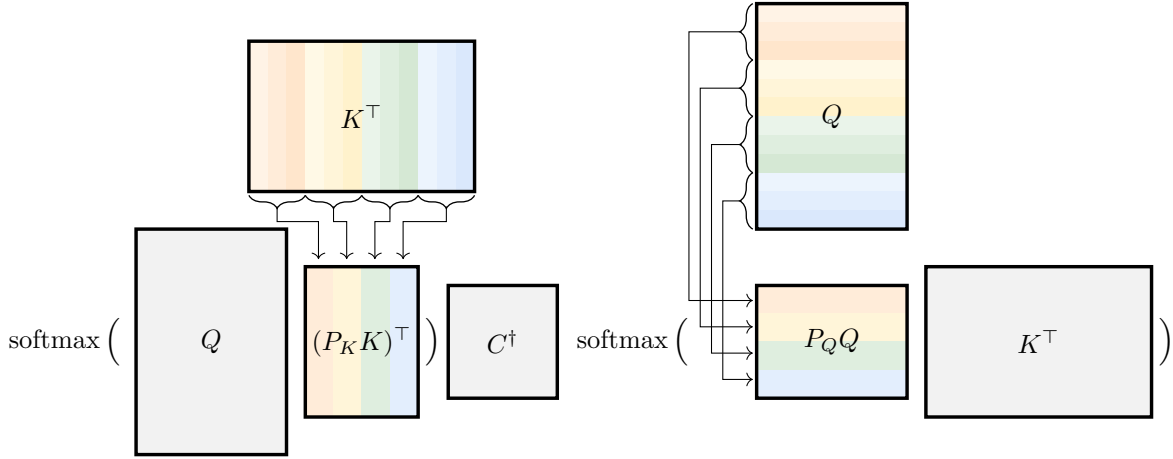

Just like the Nyströmformer, the compressed self-attention with low-rank approximation (CSALR) module proposed in \cite{chen2020compressed} also defines the projections to be the segment-means over a certain number of adjacent queries and keys. It computes the same approximation \eqref{equ:nystromformer} with a small twist: it uses a transformed matrix $C'$ whose entries are defined as 
\begin{equation*}
    C'(i, j) = \begin{cases}
        C(i, j) + 1, & \text{if $P_\Q(i, j) = P_\K(j, k),~\text{for all } k$;} \\
        C(i, j), & \text{else,}
    \end{cases}
\end{equation*}
This guarantees that $(C')^{-1}$ always exists. An identical transformation is also applied to the first softmax matrix in \eqref{equ:nystromformer}.

\section{Kernel-based Methods}\label{sec:kernel-based}

The attention output for a query $q_i$ can be formulated generally as a weighted average of the value vectors $\{v_j\}$, where the weights are determined by a similarity kernel function $\kappa(\cdot,\cdot)$. More specifically, we can rewrite $y_i$, the $i$-th row of $\Att$, as
\begin{equation}\label{eqn:kernel_att}
    y_i = \frac{\sum_{j \in [N]} \kappa(q_i, k_j) v_j^{\T}}{\sum_{j' \in [N]} \kappa(q_i, k_{j'})},
\end{equation}
where $v_1, \ldots, v_N \in \RR^{d}$ are the rows of $\V$. 
The standard self-attention mechanism~\eqref{eqn:att_matrix} relies on the softmax function, which implicitly uses an exponential dot-product kernel $\kappa(q,k) = \exp \left ( \frac{\langle q, k \rangle}{\sqrt{\dhead}} \right )$. In this section, we abuse notation and use the letters $\AttWeights, \A, \Att$ to indicate the attention mechanism corresponding to a more general kernel function $\kappa$. This formulation suggests to view the attention mechanism through the lens of kernel methods~\cite{tsai2019transformer, choromanski2022rethinkingattentionperformers}. When there exists a finite-dimensional feature map $\phi: \RR^{\dhead} \to \RR^M$ such that the kernel can be expressed as an inner product $\kappa(q, k) = \langle\phi(q), \phi(k)\rangle$, it possible to compute quantities of the form~\eqref{eqn:kernel_att} in linear time. Specifically,
\begin{align*}
    y_i &= \frac{\sum_{j \in [N]} (\phi(q_i)^{\T} \phi(k_j)) v_j^{\T}}{\sum_{j' \in [N]} \phi(q_i)^{\T} \phi(k_{j'})}  = \frac{\phi(q_i)^{\T} \left( \sum_{j \in [N]} \phi(k_j) v_j^{\T} \right)}{\phi(q_i)^{\T} \left( \sum_{j' \in [N]} \phi(k_{j'}) \right)}.
\end{align*}

Using this formulation, instead of an $\mathcal{O}(N^2 \dhead)$ complexity, we can compute the whole matrix $\Att$ in time $\mathcal{O}(NM^2)$, linear in the sequence length $N$ albeit with an $M^2$ overhead due to the feature map's dimension. However, the exponential dot-product kernel does not have a finite-dimensional feature map $\phi$, motivating its approximation with a kernel that does. While prior literature has explored alternative similarity measures that inherently possess low-dimensional embeddings~\cite{katharopoulos2020transformers, hua2022transformer}, substituting the exponential kernel with these functions often degrades the model's predictive accuracy or requires architectural adjustments to recover the lost performance. Consequently, the primary objective is to identify a similarity kernel that satisfies two criteria simultaneously: it must maintain an empirical performance profile comparable to standard softmax attention, and it must permit a rigorous approximation via a low-dimensional feature map to ensure computational efficiency. We start by considering the kernel 
\begin{equation}\label{eq:poly_kernel}
\kappa(q,k) = \left ( c + \langle q, k \rangle \right )^p,
\end{equation}
where $c$ is a constant and $p$ is some positive (not necessarily integer) number.

We begin in Section~\ref{subsec:levattention} with LevAttention which tackles the problem through sparsity by utilizing leverage scores to identify a small subset of heavy hitter keys that contribute most significantly to the attention output for polynomial kernels. Building on the utility of polynomial representations, Section~\ref{subsec:polysketch} introduces PolySketchFormer, which substitutes the standard exponential kernel with a high-degree even polynomial, employing Approximate Matrix Multiplication (AMM) and recursive sketching to compute high-dimensional feature maps in linear time. Section~\ref{subsec:tensorsketch} describes Tensor Sketch, which uses an efficient random feature map that avoids explicit high dimensional tensor products by leveraging CountSketch and the fast Fourier transform. Finally, in Section~\ref{subsec:performer} we review Performer, which uses Random Orthogonal Positive Features to approximate the exponential kernel.

\subsection{Polynomial kernels}
\subsubsection{LevAttention}
\label{subsec:levattention}

{ An attempt to address the problem of finding the important tokens making the bulk of the attention scores -- the heavy-hitters -- in the case of polynomial kernels was made in {LevAttention}~\cite{kannan2024levattention}, where the authors attempt to identify a set with the guarantee that all the indices involved in large attention scores are contained in this set.

Let $f: \RR \to \RR$ be a function that takes nonnegative values and consider its elementwise application to $\Q\K^\top$, followed by a normalization of rows so that the sum of each row is $1$. We define the $f$-sensitivity of the $j$-th row $k_j$ of $\K$ as
\begin{equation*}
    \alpha^f_j = \sup_{y \neq 0}\frac{f( \langle k_j, y \rangle)}{\sum_{\ell=1}^N f(\langle k_\ell, y \rangle)}.
\end{equation*} 
By choosing the function $f$ such that, for any $x,y\in\RR^{\dhead}$, we have $\kappa(x,y)=f(\langle x,y\rangle)$, we are able to use this formalism to relate the attention kernel to the properties of $f$ namely through its $f$-sensitivity.

Indeed, for a given $\varepsilon > 0$, the idea of LevAttention is to define the set $\mathcal U = \{ i \in [N] \mid \alpha_i^f > \varepsilon \}$, and approximate the attention matrix $\AttWeights = \Z^{-1}\A$ (corresponding to the function $f$) with a sparse matrix which only has nonzero entries in the columns corresponding to indices $j \in \mathcal U$. Intuitively, the number $\alpha_{j}^f$ can be seen as identifying the maximum value that the normalized attention scores will attribute to the $j$-th key for any input query. By constructing $\mathcal U$, we identify which keys have the possibility to take values greater than $\varepsilon$, motivating the approximation of $\AttWeights$ by restricting attention to the tokens in $\mathcal U$.
More precisely, LevAttention considers the approximation
\begin{equation*}
    \AttWeights(i,j) \approx \frac{f(\langle q_i, k_j \rangle)}{\sum_{\ell \in \mathcal{U}} f(\langle q_i, k_{\ell}\rangle)}.
\end{equation*}

The paper focuses on the case in which $f(x) = |x|^p$, for some $p \ge 1$, which corresponds to the polynomial kernel~\eqref{eq:poly_kernel} when $c=0$ and $p$ is an even positive integer. The quantity
\begin{equation*}
    \psi^f = \sup_{\K' \in \RR^{N \times \dhead}}\sum_{i=1}^{N} \alpha^f_i(\K')
\end{equation*}
is a metric that identifies the degree to which a set of $N$ vectors in the feature map space can be mutually orthogonal and satisfies $\psi^f \le d$ if $p \in [1, 2]$ and $\psi^f \le d^{p/2}$ if $p \ge 2$.

This implies that there exists a universal set $\widetilde{ \mathcal{U}}$ of size bounded by $\psi^f / \varepsilon$ such that, for any choice of $\K$, one has $\widetilde{ \mathcal{U}} \supseteq \mathcal U$. Note that the cardinality of $\widetilde{ \mathcal{U}}$, hence the cardinality of $\mathcal U$, is bounded by a quantity that only depends on $p$ and on $d$, but not the sequence length $N$, as it corresponds to the maximum rank of a linear map in the feature map space.

In practice, the set $\mathcal{U}$ for $p = 2$ can be computed via a QR factorization of $\K$, because the $f$-sensitivities coincide with the leverage scores of $\K$; for $p > 2$, upper bounds on the sensitivity scores can be computed using the so-called $\ell^p$-Lewis sketches, with a cost of $\mathcal O(\mathrm{nnz}(\K) + \textrm{poly}(d/\varepsilon) \textrm{poly}(\log N))$; see~\cite{cohen2014ellprowsamplinglewis}. Optionally, one can estimate the normalization factor $\sum_{\ell=1}^N f\left ( \langle q_i,  k_\ell\rangle \right )$ more precisely; see Theorem 3.2 in~\cite{kannan2024levattention} for details.

\subsubsection{PolySketchFormer}
\label{subsec:polysketch}
 PolySketchFormer \cite{kacham2023polysketchformer} proposes replacing the exponential kernel with a polynomial kernel of the form $\kappa(q, k) = \langle q, k\rangle ^p$, for some large even integer $p$. The associated normalized polynomial attention weight is formulated as
\begin{equation}\label{eq:polynorm}
\AttWeights^{(p)}(i,j) = \frac{\langle q_i, k_j\rangle ^p}{1 + \sum_{j' \in [N]} \langle q_i,k_{j'}\rangle ^p},
\end{equation}
where $q_i$ and $k_j$ are the rows of the matrices $\Q, \K$ which are outputs of a Layer Normalization module (\cite{ba2016layernormalization}), leading to a full attention output of 
\begin{equation*}
\Z^{-1} \left (\Q{\K}^{\T} \right)^{\odot p} \V, \text{ where } \Z = \diag\left (\mathbf{1}_N + (\Q{\K}^{\T})^{\odot p} \mathbf{1}_N\right )
\end{equation*}
and $(\cdot)^{\odot p}$ denotes the element-wise power. 
Note that the ``1'' in the denominator of~\eqref{eq:polynorm} is added to avoid the denominator becoming zero (which does not happen with the softmax kernel but may happen with a polynomial kernel). The choice of even $p$ is an algebraic choice from the authors which allows proper normalization to a probability distribution. However, this means large negative inner products are treated as having high similarity; an extension was proposed in Polaformer ~\cite{meng2025polaformerpolarityawarelinearattention} whereby negative dot products were taken into account in the approximation.

\paragraph{Dependence on the degree parameter.} We take a moment to expand on the capability of the softmax approximation, governed explicitly by the degree parameter $p$, by examining the behavior of the distribution $\AttWeights^{(p)}$ at its lower and upper theoretical limits. When $p=0$, assuming non-zero vectors, the term $\langle q_i, k_j \rangle^0$ is 1 for all entries. Consequently, one gets the uniform distribution
$$ \AttWeights^{(0)}(i,j) = \frac{1}{N}, $$
meaning that the model attends equally to all tokens regardless of relevance. Secondly, let us consider the limit as $p \to \infty$. Let $\mu_i = \max_{j' \in [N]} \langle q_i, k_{j'} \rangle$ be the maximum inner product for the $i$-th query. By dividing the numerator and denominator by $\mu_i^p$, we rewrite the weight as
$$
\AttWeights^{(p)}(i,j) = \frac{\left (\langle q_i, k_j \rangle / \mu_i\right )^p}{\mu_i^{-p} + \sum_{j' \in [N]} \left (\langle q_i, k_{j'} \rangle / \mu_i\right )^p}.
$$
Assuming reasonably that $\mu_i > 1$ (a consequence of the high-dimensional scaling preserved by Layer Normalization), the term $\mu_i^{-p}$ vanishes. For any index $j$ where $\langle q_i, k_j \rangle < \mu_i$, the ratio is strictly less than 1, causing the term to decay to 0. Conversely, for indices in the set of maximisers $\mathcal{S}_i = \{j \mid |\langle q_i, k_j \rangle| = \mu_i\}$, the ratio is 1. Thus, the distribution converges to
$$
\lim_{p \to \infty} \AttWeights^{(p)}(i,j) = \begin{cases} \frac{1}{|\mathcal{S}_i|} & \text{if } j \in \mathcal{S}_i \\ 0 & \text{otherwise} \end{cases}.
$$
This demonstrates that by choosing a sufficiently large parameter $p$, the polynomial formulation transitions from a uniform distribution to a highly peaked distribution, recovering the hardmax operation in the limit.

\paragraph{Attention weight approximations in linear time.} To explain how to quickly compute (or approximate) the weights in linear time, we employ the Kronecker product, denoted by $\otimes$, which satisfies the scalar property $\langle u, v \rangle \langle x, y \rangle = \langle u \otimes x, v \otimes y \rangle$ for any vectors $u,v,x,y$ of compatible sizes. Extending this recursively allows us to express the polynomial kernel as an inner product of $p$-fold tensor products $\langle q, k \rangle^p = \langle q^{\otimes p}, k^{\otimes p} \rangle$, where $q^{\otimes p}$ denotes the Kronecker product of $q$ with itself repeated $p$ times. Then, we can rewrite
\begin{equation*}
    \left (\Q{\K}^{\T}\right )^{\odot p} \V = {\Q}^{\otimes p}\left ({\K}^{\otimes p}\right )^{\T} \V,
\end{equation*}
where, for a matrix $B\in\RR^{d_1\times d_2}$ with $d_1,d_2\in\NN^*$, the rows of the matrix $B^{\otimes p}\in\RR^{d_1\times d_2^p}$ are the Kronecker products, repeated $p$ times, of the rows of $B$ with themselves. By computing 
$(\K^{\otimes p})^\T \V$ first, the overall complexity becomes $\mathcal O\left (N\dhead^{p+1}\right )$. While this complexity is linear with respect to the sequence length $N$, the polynomial dependence on the head dimension $\dhead$ is computationally prohibitive for the large degrees $p$ required to match softmax quality, thus necessitating the use of sketching. 

\paragraph{Attention weight approximation by sketching.} Approximate Matrix Multiplication (AMM), introduced by \cite{woodruff2014sketching}, seeks a sketching matrix $S \in \mathbb{R}^{\dhead^p \times r}$ that satisfies the $(\varepsilon, p)$-AMM property
$$\norm{(\Q^{\otimes p}S)(\K^{\otimes p}S)^\T - \Q^{\otimes p}(\K^{\otimes p})^\T}_{\frob} \le \varepsilon \norm{\Q^{\otimes p}}_{\frob} \norm{\K^{\otimes p}}_{\frob},$$
where $r$ is the sketching dimension. Theoretically, $r$ should scale with $\mathcal{O}(\epsilon^{-2})$, however, typically one can set it to a small multiple of $\dhead$, e.g., $r \in \{32, 64\}$, in practice. The product $\phi(\Q) = \Q^{\otimes p}S$ is computed efficiently using a recursive algorithm from \cite{ahle2020oblivious}, where for $p>1$ the computation is defined as
$$\phi_p(\Q) = \sqrt{1/r} \cdot \left[ (\phi_{p/2}(\Q)G_1) \odot (\phi_{p/2}(\Q)G_2) \right],$$
where $G_1$ and $G_2$ are independent random Gaussian matrices (with dimensions $\dhead \times r$ at the base level $p=2$, and $r \times r$ for recursive levels $p>2$) and $\odot$ denotes the Hadamard product. A critical issue with this standard sketch however is its failure to preserve the non-negativity of the polynomial kernel for even $p$. Therefore, a self-tensoring technique is introduced where an intermediate feature map for degree $p/2$ is computed as $\phi_{p/2}(q) = (q^{\otimes p/2})^\T S_{p/2}$, where $S_{p/2}$ is a sketching matrix of size $\dhead^{p/2} \times r$ 
and the final feature map $\phi'(q)$ is then defined as its Kronecker square, $\phi'(q) = (\phi_{p/2}(q))^{\otimes 2}$. This guarantees a non-negative result since $\langle \phi'(q), \phi'(k)\rangle  = \langle \phi_{p/2}(q), \phi_{p/2}(k)\rangle ^2 \ge 0$. Importantly, if a random sketch $S$ satisfies some necessary JL-moment properties for a degree-$p/2$ kernel, the resulting self-tensored sketch provides a valid approximation for the full degree-$p$ kernel, and the error for the polynomial attention matrix is bounded by 
$$
\norm{ \phi'(\Q)\phi'(\K)^\T - (\Q\K^\T)^{\odot p} }_{\frob} \le \varepsilon\norm{\Q^{\otimes p}}_{\frob}\norm{\K^{\otimes p}}_{\frob}
$$
with high probability.

\subsubsection{Tensor Sketch}
\label{subsec:tensorsketch}
Tensor Sketch \cite{pham2025tensor,pham2013fast} is a random feature map for approximating the polynomial kernel~\eqref{eq:poly_kernel}. Note that a non-zero constant $c$ can be incorporated by appending $\sqrt{c}$ to the input vectors, allowing the method to treat the kernel as homogeneous; for our specific attention formulation, we effectively set $c=0$. As seen previously, the core challenge with polynomial kernels is the computational cost associated with their explicit feature maps. The exponential growth of these feature spaces makes direct computation infeasible, so sketching is a direction to overcome this bottleneck. Tensor Sketch builds upon CountSketch \cite{charikar2002finding}, a technique for dimensionality reduction that approximately preserves inner products. A CountSketch of a vector $x \in \mathbb{R}^{\dhead}$ is a linear projection into $\mathbb{R}^r$ defined by two hash functions, a binning function $\mathfrak{h}: [\dhead] \to [r]$ from a 2-wise independent family and a sign function $\mathfrak{s}: [\dhead] \to \{-1, 1\}$ from a 4-wise independent family, where the $b$-th component of the sketch $Cx$ is given by
$$ (Cx)_b = \sum_{j\in\{i \in [\dhead]|\mathfrak h(i)=b\}} \mathfrak s(j)x_j, $$
where the indexing on $x_j$ denotes the $j$-th scalar component of the vector $x$. The inner product of two sketched vectors, $\langle Cx, Cy\rangle $, serves as an unbiased estimator for $\langle x, y\rangle $ with variance bounded by 

$$ \mathbb{V}[\langle Cx, Cy\rangle ] = \frac{1}{r}\left(\sum_{i\ne j}x_{i}^{2}y_{j}^{2}+\sum_{i\ne j}x_{i}y_{i}x_{j}y_{j}\right) \le \frac{2}{r}\|x\|^2 \|y\|^2. $$

To efficiently compute a CountSketch of the high-dimensional tensor product $x^{\otimes p}$, TensorSketch avoids constructing the massive vector directly. To formalize this, let $\{C_t\}_{t=1}^{p}$ be $p$ independent CountSketches of the original vector $x$, each represented as a polynomial
$$
\mathfrak{P}_t(\omega) = \sum_{i=1}^{\dhead} \mathfrak{s}_t(i) x_i \omega^{\mathfrak{h}_t(i)}.
$$
Here, $\omega$ is a formal indeterminate in the polynomial ring $\mathbb{R}[\omega]$ whose exponent tracks the mapped bucket index. A key property of this representation is that polynomial multiplication natively generates the cross-terms of the tensor product, yielding
$$
\prod_{t=1}^p \mathfrak{P}_t(\omega) = \sum_{i_1, \dots, i_p} \left(\prod_{t=1}^p \mathfrak{s}_t(i_t) x_{i_t}\right) \omega^{\sum_t \mathfrak{h}_t(i_t)}.
$$
Notice how the algebra perfectly simulates the sketching process: multiplying the coefficients constructs the tensor product elements, while adding the exponents naturally computes the combined hash bucket for those elements. To finalize the sketch, we reduce the polynomial's exponents modulo the target dimension $r$, summing together all coefficients whose exponents leave the same remainder. The $r$ coefficients of this compacted polynomial correspond exactly to the CountSketch of $x^{\otimes p}$. Crucially, this entire mathematical sequence is algebraically equivalent to the circular convolution of the $p$ original sketches. The procedure runs in $\mathcal O(\dhead + r \log r)$ time and is given by
$$ \phi_{\mathrm{TS}}(x) = \text{FFT}^{-1} \left( \text{FFT}(C_1x) \odot \dots \odot \text{FFT}(C_px) \right)$$
where $\odot$ denotes the element-wise product between vectors and $r$ is the Tensor Sketch dimension. Applying this to our attention formulation, we use this feature map on specific query vectors $q_i$ and key vectors $k_j$. Given their sketches $\phi_{\mathrm{TS}}(q_i), \phi_{\mathrm{TS}}(k_j) \in \RR^r$, we have that
\begin{equation*}
    \mathbb{E}\left[\left\langle \phi_{\mathrm{TS}}(q_i), \phi_{\mathrm{TS}}(k_j)\right\rangle \right] = \left\langle q_i^{\otimes p}, k_j^{\otimes p}\right\rangle  = \langle q_i, k_j\rangle ^p
\end{equation*}
and the variance satisfies
\begin{equation*}
    \mathbb{V}\left[\left\langle \phi_{\mathrm{TS}}(q_i), \phi_{\mathrm{TS}}(k_j)\right\rangle \right] \leq \frac{3^p-1}{r}\|q_i\|^{2p}\|k_j\|^{2p},
\end{equation*}
allowing for the efficient estimation of the polynomial attention weights without explicitly forming the high-dimensional tensor products.

\subsection{Performer method: Randomized kernel features}\label{subsec:performer}
\newcommand{\NR}{N_{R}}
\newcommand{\Nf}{N_{f}}

Similar  to the methods described above, the main principle behind the Performer architecture proposed in \cite{choromanski2022rethinkingattentionperformers}, and generalized in \cite{likhosherstovChefsRandomTables2022,likhosherstovFAVORSharpAttention2023}, is to reduce the computational
 cost of calculating attention to $\mathcal O(N)$ by reformulating the calculation of $A$ to allow reordering the order of the matrix multiplications in the attention mechanism. Instead of using kernels with finite-dimensional feature maps to represent the key/value vectors, Performer uses  their projection on random features such that, in expectation, the scalar product between the projections approximates the kernel. Unlike polynomial methods, Performer's architecture is not restricted by the dimension of the feature map space and can be applied to regular attention.

More formally, by constructing stand-in matrices, $\Q'$ and $\K'$, for the query and keys respectively, such that $\mathbb E[\langle q_i', k_j' \rangle]= \kappa(q_i,k_j)$, Performer's architecture allows a linear complexity by changing the order of computation to calculate ${\K'}^\T \V$ first and $(\Z')^{-1} \Q' \left ({\K'}^\T \V\right )$ second which has cost $\mathcal O(Nrd)$, where $Z':=Q'\left ({K'}^{\T}\mathbf 1_N\right )$.

 Given an ordered set $\omega = \left( \omega_{i}\right)_{1\le i\le \NR}$ of $\NR\in\NN^*$ random vectors in $\RR^{\dhead}$,
 an ordered set of $\Nf\in\NN^*$ deterministic functions, $f=\left(f_j:\RR\to\RR\right)_{1\le j\le \Nf}$, and a function $\mathfrak{h}:\RR^{\dhead}\to \RR$, the authors introduce
the Random Orthogonal Positive features function \(\phi : \RR^{\dhead} \to \RR^r_{+} \) as 
\begin{align}
    \phi(x) = \frac{\mathfrak{h}(x)}{\sqrt{\NR}}\left (f_1(\omega_1^\T x), \ldots,f_1(\omega_{\NR}^\T x), \ldots, f_{\Nf}(\omega_1^\T x), \ldots,  f_{\Nf}(\omega_{\NR}^\T x) \right )
\end{align}
 where the rank $r$ is thus given by $r=\NR\Nf$.
These randomized features can be engineered, through the selection of $f$, $\omega$ and $\mathfrak{h}$, such that for any $q, k\in\RR^{\dhead}$,
\begin{equation*}
\mathbb E [\langle\phi(q),\phi(k)\rangle]= \exp\left (\frac{\langle q,k\rangle}{\sqrt{\dhead}} \right ).
\end{equation*}

The authors concentrate on three different  Random Features, the first is based on Random Fourier features \cite{rahimiWeightedSumsRandom2008} and the other two are novel and ensure the positivity of the scalar product and a reduced variance of the estimators:
\begin{itemize}
\item  $\phi_{\text{trig}}$, with parameters $\mathfrak h:x \mapsto \exp\left ({\|x\|^2}/2 \right )$, and functions $f_1: u\mapsto \sin(u)$, $f_2:u \mapsto \cos(u)$, and $\omega_p \sim \mathcal{D} = \mathcal{N}(0, I_d)$;
    
    \item $\phi_{+}$, with parameters $\mathfrak h:x \mapsto \exp\left (-{\|x\|^2}/2\right )$, $f_1: u\mapsto \exp(u)$, and $\omega_p \sim \mathcal{D} = \mathcal{N}(0, I_d)$;
    \item  $\phi_{\text{hyp}+}$, based on the hyperbolic cosine, with parameters $\mathfrak h:x \mapsto \frac1{\sqrt{2}}\exp\left (-{\|x\|^2}/2\right )$, $f_1: u\mapsto \exp(u)$, $f_2:u \mapsto \exp(-u)$, and $\omega_p \sim \mathcal{D} = \mathcal{N}(0, I_d)$.
\end{itemize} Note that the above PRF are written without the rescaling factor $1/{\sqrt {\dhead}}$, for simplicity.

The output dimension \(r\) can be further reduced (which also translates to requiring fewer random features) by selecting the random vectors $\left( \omega_{p}\right)_{1\le p\le m}$ to be orthogonal. This involves orthogonalizing general \(\omega_{p} \sim \mathcal{D}\). In practice, if more than $d$ vectors are needed to achieve a target performance, the authors propose using the same procedure to construct multiple blocks composed of $d$ orthonormal vectors.

The above approximation procedure leads to low-variance and low-rank attention computation that is linear in time and space complexity. Furthermore, in practice, noticing that the matrices ${K'}^\T$ and $(V\mid 1_N)$ can be combined into ${K'}^\T(V\mid 1_N)\in \RR^{\NR\Nf\times (d+1)}$, the KV cache size can be reduced to $\mathcal O(\NR \Nf (d+1))$, albeit at the cost of keeping the set of random vectors $\omega$, of size $\NR d$, in memory; this provides a memory benefit in the event that $\NR(\Nf+1)\le 2N$.
It is important to note, however, that while such methods approximate attention with low variance, they fail to accurately represent spiked distributions in attention, due to the projection onto random vectors blurring the output. 
The Scatterbrain article~\cite{chen2021scatterbrain}, which provides the experimental results highlighting the sparsity and low-rankness of attention mechanisms, attempts to solve this by combining one of these randomized features, namely $\phi_+$, to leverage low-rankness, along with an LSH-based method to make use of the sparsity of the attention matrices.
Several extensions to this method were proposed in \cite{likhosherstovChefsRandomTables2022,likhosherstovFAVORSharpAttention2023}. Performer can be reformulated to treat $\omega$ and $x$ differently in each function $f_j$, for $1\le j\le\Nf$, by integrating $\mathfrak h$ into the $f_j$ functions, and setting 
$$f_j(\omega, x)= D\exp\left(\omega^\T A\omega + \omega^\T B_j x +x^\T C_j x\right),$$ where $A\in\RR^{\dhead\times \dhead}$, $D\in\RR$ and for any $1\le j\le \Nf$, $B_j,C_j\in\RR^{\dhead\times \dhead}$. Using the above formulation, the PRF described in \cite{choromanski2022rethinkingattentionperformers} can be obtained by setting $A=0$, $B_1=I_{\dhead}$, $C_j=-I_{\dhead}$ and $D=1$ for $\phi_+$ and $A=0$ $B_j=(-1)^j I_{\dhead}$ and $D=1/\sqrt 2$ in the case of $\phi_{\mathrm{hyp}+}$. 
Furthermore, in \cite{choromanski2022hybrid}, the Hybrid Random Features method was proposed, which allows to combine estimators by adaptively choosing a subset of estimators best suited for a given input.

\section{Latent Attention: An important variant of the attention mechanisms}\label{sec:MLA}
While the most prevalent form of attention in recent large language models makes use of the MHA mechanisms described in Section \ref{sec:original_notation}, recent alternative attention mechanisms, such as Latent Attention, introduced by DeepSeek in \cite{deepseekai2024deepseekv2strongeconomicalefficient}, have shown increased efficiency through the use of a shared latent space for the key/value vectors. 
In this section, we highlight the lifting of equivalence between the attention models described previously and Latent Attention which occurs when positional encodings are applied.
To do so, we first define the Multi-headed Latent Attention (MLA) mechanism, then we describe the ubiquitous Rotary Positional Embeddings (RoPE) \cite{suRoFormerEnhancedTransformer2023} method and its MLA counterpart, highlighting the fundamental differences in the way the models apply the positional embeddings.
Finally, we describe the TransMLA \cite{meng2025transmla} method, which consists of approximately converting a GQA model (including MHA/MQA) into an MLA model. This provides increased efficiency through the use of the DeepSeek MLA pipeline at a reasonable cost in terms of performance after finetuning.

\subsection{Latent Attention}
\label{subsec:latent}
Multi-headed Latent Attention (MLA) consists of using a latent space, of dimension $d_L\in\NN^*$, from which the $\Nh$ head embeddings, of dimension $\dhead\in\NN^*$, are reconstructed.
For the key and value heads, a shared latent space, obtained through $\WLKV\in\RR^{d\times d_{L}}$, is used, while a separate latent space, obtained through $\WLQ\in\RR^{d\times d_{L}},$ is used for the query embeddings.
The embeddings are given by \begin{equation}
\LQ=X\WLQ,\quad  \LKV=X\WLKV  .
\end{equation}
For each head, we set \[
\head{\WQ}{h}=\WLQ{\head{\WUQ}{h}},\quad \head{\WK}{h}=\WLKV{\head{\WUK}{h}},\quad \text{ and }\head{\WV}{h}=\WLKV{\head{\WUV}{h}} \]  where, for the $h$-th head, ${\head{\WUQ}{h}}\in\RR^{d_{L}\times \dhead}, {\head{\WUK}{h}}\in\RR^{d_{L}\times \dhead}$ and ${\head{\WUV}{h}}\in\RR^{d_{L}\times \dhead}$ are the query, key and value ``up-projection'' weight matrices from the latent subspaces. 
Hence, we obtain the factorizations
\[ 
\head{\Q}{h}=\LQ{\head{\WUQ}{h}},\quad \head{\K}{h}=\LKV{\head{\WUK}{h}}\quad \text{ and }\head{\V}{h}=\LKV{\head{\WUV}{h}}.
\]
The attention scores are then obtained as
\begin{equation}   
    \head{\A}{h}= \exp \left ( \frac{\head{\Q}{h} \head{\K}{h}^\top}{\sqrt{\dhead}} \right ) = \exp\left (\frac{\LQ{\head{\WUQ}{h}}\left({\LKV}\head{\WUK}{h}\right)^{\T}}{\sqrt {\dhead}} \right ).\label{eq:attn_score_latent}
\end{equation}
This means that, at inference time, the weight matrices can be merged into $\head{\WQK}{h}={\head{\WUQ}{h}}\left({\head{\WUK}{h}}\right)^\T\in\RR^{d_{L}\times d_{L}}$. Setting ${\Atth{h}^L}=\head{\Z}{h}^{-1}\head{\A}{h}\LKV$ and defining ${\head{\WUV}{h}}\in\RR^{d_L\times \dhead}$ such that $\head{\V}{h}=\LKV{\head{\WUV}{h}}$, the result of each head is then given by
\begin{equation}
    \Atth{h}=\head{\Z}{h}^{-1}\head{\A}{h}\underbrace{\LKV{\head{\WUV}{h}}}_{=\head\V h}=\Atth{h}^L{\head{\WUV}{h}}.
\end{equation}
The output of the MLA mechanism then becomes
\begin{equation}
    O=\begin{pmatrix}\Atth{1}&\cdots&\Atth{\Nh}\end{pmatrix}\WO =\concat{h\in\Nh}{\head{\Att}{h}^L}\WLVO
\end{equation}
where $\WLVO\in \RR^{\Nh d_L\times d}$ corresponds to the inference-time absorption of the value weight matrices into the output weight matrix $\WO\in\RR^{\Nh\dhead\times d}$.

\subsubsection{Computational Cost (Inference)}
The initial projection of the input sequence $X$ into the query and key-value latent spaces of dimension $d_L$ costs $\mathcal O(N d d_L)$. For each of the $\Nh$ heads, the attention scores are computed via two consecutive matrix multiplications ($\LQ \head{\WQK}{h}$ and the result with $(\LKV)^\T$), which have a combined cost of $\mathcal O(N d_L^2 + N^2 d_L)$. The application of the attention scores to the latent key-value matrix $\LKV$ requires a further multiplication costing $\mathcal O(N^2 d_L)$. The final output projection aggregates the head results and multiplies by an output matrix, incurring a cost of $\mathcal O(N \Nh d_L d)$.

\begin{center}
\begin{tabular}{|l|l|}
\hline
\textbf{Operation} & \textbf{Computational Cost} \\
\hline
Latent Projections ($\LQ, \LKV$) & $\mathcal O(N d d_L)$ \\
Attention Scores ($\head{\A}{h}$) (per head) & $\mathcal O(N d_L^2 + N^2 d_L)$ \\
Attention Output ($\head{\Att}{h}^L$) (per head) & $\mathcal O(N^2 d_L)$ \\
Final Projection ($O$) & $\mathcal O(N \Nh d_L d)$ \\
\hline
\end{tabular}
\end{center}

The total complexity is 
\begin{equation*}
    \mathcal O\left (N d d_L + \Nh(N d_L^2 + 2 N^2 d_L) + N \Nh d_L d\right )
\end{equation*}
and it is dominated by $\mathcal O( N^{2}\Nh d_L)$. In the case where $\Nh d_{L}\approx d$, this is comparable to a MHA model with hidden size $d=\Nh\dhead$.
However, in terms of memory cost, note that in MLA we only need to keep the latent vectors $\LKV$, of size $\RR^{N\times d_{L}}$, instead of all key and value vectors for every individual $KV$ head, i.e. a tensor of size ${N\times 2\Ng\times \dhead}$.

\subsection{Rotary Position Embeddings}\label{sec:RoPE}
In the previous sections, we have considered attention mechanisms with no positional encodings.
In practice, positional encodings are added to the attention models to provide more information to the model about the position of the tokens in the text. Indeed, apart from the implicit positional dependence that appears in masked attention, the sequential nature of a sentence is lost in the attention mechanism, inducing the need for such positional encodings.
For GQA and MHA, 
the most common strategy is the Rotary Position Embeddings (RoPE) method: it consists in using rotation matrices, applied to each successive pairs of dimensions $p\in \{ (2\ell+1,2\ell+2) \mid {0\le \ell\le d/2-1} \}$, to induce positional encoding based on the relative distance between tokens.
Indeed, given the $i$-th and $j$-th token embeddings, RoPE encodes the relative distance $i-j$ through rotations of angles $\theta_\ell=b^{-2\ell/d}$, for $0\le \ell\le d/2 - 1$  and where $b\in\RR$ is known as the base, applied to the pair of dimensions indexed by $(2\ell+1,2\ell+2)$. This corresponds to a block diagonal matrix with $2\times2$ blocks corresponding to the rotations

\begin{equation}\label{eq:rope_mat}
    R_j^d = \begin{bmatrix}
    R_j(\theta_0) & & & & \\
    & R_j(\theta_1) & & & \\
    & & \ddots & & \\
    & & & R_j(\theta_{d/2-1}) & \\
    \end{bmatrix}
\end{equation}
where $R_j(\theta_
\ell) = \begin{bmatrix} \cos(j\theta_\ell) & -\sin(j\theta_\ell) \\ \sin(j\theta_\ell) & \cos(j\theta_\ell) \end{bmatrix}$.

The angles $\theta_\ell$, for each pairs of dimensions, and the base $b\in\RR$ are set so as to allow the model to encode more positional information through the different frequencies.
This positional embedding is applied to both the query and key embeddings in the attention score calculation, giving:
\begin{equation}
    A(i,j)=\exp\left(\frac{\Q(i,:)R^{d}_{i}{R^{d}_j}^\T\K(j,:)^\T}{\sqrt{\dhead}}\right)=\exp\left(\frac{\Q(i,:)R^{d}_{i-j}\K(j,:)^\T}{\sqrt{\dhead}}\right),
\end{equation}
where we omit the scaling factor for simplicity. 
We define the rotated query and key vectors as \[\Q^{R}(i,:)=\Q(i,:)R^{d}_i\quad \text{ and } \K^{R}(j,:)=\K(j,:)R^{d}_j.\] For KV caching (\ref{par:KVcaching}), the $\K^{R}$ matrix, instead of $\K$, is stored in memory so that the rotations are computed only once for each token.
In the case of GQA or MHA, the same rotation matrices are applied to all heads, we also denote by $\head \Q h^R=\head \Q h (i,:)R^{\dhead}_i$ and $\head \K g^R=\head \K g (i,:)R^{\dhead}_i$ the associated query and key matrices for each query head $h\in[\Nh]$ and key head $g\in[\Ng]$.
We note that as $R^{d}$ possesses a $2\times 2$ block diagonal structure, the effect of the positional embedding is negligible on the computational cost of the attention mechanism as its complexity is of the order of $\mathcal O(N\dhead)$.
\subsubsection{Rotary Positional Embeddings for Latent Attention}
When trying to apply RoPE to Latent Attention, a problem arises: indeed, one of the performance improvements of the architecture resides in the fact that we can merge ${\head{\WUQ}{h}}\left({\head{\WUK}{h}}\right)^\T$ from eq. \eqref{eq:attn_score_latent} into a single matrix $\head{\WQK}{h}$. This improvement is not compatible with the RoPE method above as it would consist in inserting a different rotation matrix for each input vector in between ${\head{\WUQ}{h}}$ and $\left(\head{\WUK}{h}\right)^\T$.
To remedy this while maintaining the increase in performance from the use of positional encoding, \cite{deepseekai2024deepseekv2strongeconomicalefficient} introduces a new method which consists in enlarging  each vector of each head by appending a vector containing the positional information. 

To do so, letting $\dr\in\NN^*$ be the dimension of the rotary embedding, new weight matrices are introduced.
For all key heads, a single shared weight matrix  $\WRK\in\RR^{d \times \dr}$ is introduced.
For each query head, a weight matrix $\head{\WRQ}{h}\in\RR^{\dl \times \dr}$, with $1\le h\le \Nh$, is introduced such that, for any token $i\in\NN^*$,  we set
\begin{equation}
    \KRL(i,:)=X\WRK R^{\dr}_{i},\qquad \head{Q^R}{h}(i,:)=\LQ\head{\WRQ}{h} R^{\dr}_{i},
\end{equation} where the rotation matrices, $R^{\dr}_n$, are defined as in eq. \eqref{eq:rope_mat}.
Finally, the attention scores are modified such that 
\begin{equation}
    \head{\A}{h}=\exp\left(\head{Q^{\MLA}}{h}{\head{\K^\MLA}{h}}^\T\right),
\end{equation}
where $\head{\Q^\MLA}{h}=\begin{bmatrix}\head{\Q}{h}&\head{Q^R}{h}\end{bmatrix}\in\RR^{N\times (\dhead+\dr)}$,  $\head{\K^\MLA}{h}=\begin{bmatrix}\head{\K}{h}&\KRL\end{bmatrix}\in\RR^{N\times (\dhead+\dr)}$ and we omit the scaling factor for simplicity.
 Notice that, as $\KRL$ does not vary from head to head, only $\LKV$ and $\K^{R,\mathrm{MLA}}$ need to be cached for MLA with RoPE embeddings. 
 We note that $\dr$ is generally a small fraction of $\dhead$, the computational overhead incurred is even smaller than in the case of RoPE for regular attention.
 
 \subsection{TransMLA - From GQA to MLA}\label{subsec:transmla_short}
Due to the prevalence of well-performing already trained GQA models and the optimized open-source pipelines of the DeepSeek architectures, there has been growing interest in converting from one architecture to the other in order to improve the practical performance of GQA models.
As it turns out, without rotary embeddings,
MLA and MHA/GQA can be shown to be completely equivalent to each other through algebraic rewriting.
Indeed, if we consider models where no rotation matrices are applied, by setting $$\head{\WQ}{h}= \WLQ \begin{bmatrix}
\head{\WUQ}{h} & \head{\WRQ}{h}
\end{bmatrix},\,\head{\WK}{h}= \WLKV \begin{bmatrix}
\head{\WUK}{h} & \WRK
\end{bmatrix}\, \text{ and } \,\head{\WV}{h}= \WLKV \head{\WUV}{h},$$
one can construct perfectly equivalent models.
However, their different RoPE mechanisms destroy this equivalence and no procedure allows for a direct conversion from one form to the other as the RoPE operation acts on different parts of the query/key embedding space, i.e. the complete subspace for GQA/MQA, and the last $d_R$ dimensions for MLA corresponding to $\WRK$ and $\head{\WRQ}{h}$.
While exact conversion is impossible, the TransMLA \cite{meng2025transmla} method allows this conversion with reduced performance loss compared to previous methods. 

Given a dataset $X\in\RR^{N\times d}$, TransMLA allows the concentration of the dataset-specific, most important positional features of the key heads of a GQA model into its first few heads. It does so through the application of a matrix $U^\downarrow$. This matrix is constructed by extracting, for each pair of ``key'' output dimensions $p_\ell=(2\ell+1,2\ell+2)$ indexed by $\ell\in\{0,\dhead/2-1\}$, a unitary matrix $U_\ell\in \RR^{\Ng\times\Ng}$ acting on the ``head'' dimension, mixing the information contained in each head for the $\ell$-th par of dimensions.
The authors look at pairs of dimensions as these share the same RoPE frequency $\theta_\ell$), and are mixed through the RoPE procedure\ref{eq:rope_mat}.
The $U_\ell$ matrices are chosen to maximise the following cost function
\begin{equation*}
U_\ell=\underset{\tilde U_\ell\in\RR^{\Ng\times\Ng}}{\argmax}\Tr\left([{\tilde U}_\ell \cov_\ell[K^R] {\tilde U}_\ell] (:m,:m)\right),
\end{equation*}

 where we denote, for any matrix $M\in\RR^{d_1\times d_2}$ with $d_1,d_2\in\NN^*$ and for any integer $1\le m <d_1$ and $1\le n <d_2$ and, by $M(:m,:n)$ the submatrix composed of the first $m$ rows and first $n$ columns of $M$, and we set $$\cov_\ell[\K^R](g,g')=\sum_{i=1}^{N}\K^R_{g}(i,2\ell)\K^R_{g'}(i,2\ell) + \K^R_g(i,2\ell+1)\K_{g'}^R(i,2\ell+1).$$
The solution to the above maximization problem turns out to simply be given by the $m$ first eigenvectors of the covariance matrix.

After the application of $ U^\downarrow$, the features containing the positional information are hence concentrated in the first few heads.
By also applying $U^\downarrow$ to the query matrices, we obtain a reformulation of the exact same model, thanks to the unitarity of $U_\ell$ for each $0\le \ell\le d/2-1$.
Finally, an approximation restricting the positional embedding to the first new head allows to obtain an approximate conversion to MLA. 
However this conversion is extremely inefficient in terms of matrix sizes and so the weight matrices are then compressed by considering a rank-$r$ singular value decomposition of the column-concatenated key and value matrices.
A more thorough description of the transMLA method is given in Appendix \ref{subsec:transmla}.

\section{Tensor-based attention mechanisms}\label{sec:tensors}

Tensor-based attention mechanisms are higher-order generalizations of standard matrix-based attention, developed to address computational bottlenecks in the transformer architecture.
In this section, we describe how tensor structure has been imposed and exploited for improved computational efficiency and greater expressivity. 
We consider three approaches to ``tensor-based attention'' in our work:
\begin{itemize}
    \item matrix-based quantities, -- e.g. weights, queries, keys, and values -- are folded across heads (and potentially layers), and tensor decompositions are employed to factorize and compress resulting  tensors (Section~\ref{subsec:tensorized_attention});
    \item tensor products and models are incorporated directly into the attention model, without explicit use of folding or stacking matrix-based quantities (Section~\ref{subsect:tensorial_models});
    \item  the input data itself is represented as a tensor, and tensor-based attention computations are needed to preserve or exploit its structure (Section~\ref{subsec:tensorial_inputs}).
\end{itemize}
We begin with a brief overview of tensor preliminaries relevant to our discussion. Readers who are familiar with standard tensor notation and decompositions (CP and Tucker) are recommended to proceed to Section~\ref{subsec:tensorized_attention} for further reading.

\subsection{Preliminaries and Notation}
\label{subsec:tensor_prelims}

We refer to multi-dimensional arrays as \textit{tensors} and the number of dimensions as the \textit{order}. For example, vectors are first-order tensors, matrices are second-order tensors, and so forth, as visualized in Figure~\ref{fig:tensor_basic}. When a matrix is only indexed by one subscript, this subscript refers to a sequence index, for example $\{M_{h}\}_{h=1}^n$ is the collection of $n$ matrices $M_h \in \mathbb{R}^{d_1 \times d_2}$.

We will reserve calligraphic letters for tensors of order three or more. Here, we provide a brief introduction to tensor rank, specifically the CANDECOMP-PARAFAC (CP) and Tucker rank. For simplification, we introduce the CP and Tucker decompositions for third-order tensors. We will use $*$ to denote the Khatri-Rao product and $\otimes$ for the outer product. The product $\times_k$ denotes the $k$-mode product between a tensor $\mathcal{T} \in \mathbb{R}^{d_1 \times d_2 \times ...\times d_n}$ and a matrix $M\in\mathbb{R}^{d \times d_k}$ whose output is a tensor of dimension $d_1 \times \ldots \times d_{k-1} \times d \times d_{k+1} \times \ldots \times d_n$, defined element-wise as 
$$\left(\mathcal{T} \times_k M\right)_{i_1, \ldots, i_{k-1}, j, i_{k+1}, \ldots, i_n} = \sum_{i_k = 1}^{d_k} \mathcal{T}_{i_1, \ldots, i_k,\ldots i_n} M_{ji_k}, \qquad j = 1, \ldots, d.$$

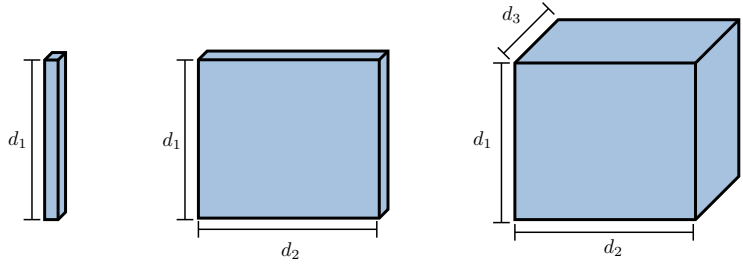
\begin{figure}
    \centering
    \tikzset{every picture/.style={line width=0.75pt}} 
\resizebox{0.6\textwidth}{!}{
\begin{tikzpicture}[x=0.75pt,y=0.75pt,yscale=-1,xscale=1]

\draw  [fill={rgb, 255:red, 166; green, 194; blue, 222 }  ,fill opacity=1 ][line width=1.5]  (439,51) -- (468,22) -- (589,22) -- (589,127) -- (560,156) -- (439,156) -- cycle ; \draw  [line width=1.5]  (589,22) -- (560,51) -- (439,51) ; \draw  [line width=1.5]  (560,51) -- (560,156) ;
\draw [fill={rgb, 255:red, 166; green, 194; blue, 222 }  ,fill opacity=1 ]   (430,51) -- (430,158) ;
\draw [shift={(430,158)}, rotate = 270] [color={rgb, 255:red, 0; green, 0; blue, 0 }  ][line width=0.75]    (0,5.59) -- (0,-5.59)   ;
\draw [shift={(430,51)}, rotate = 270] [color={rgb, 255:red, 0; green, 0; blue, 0 }  ][line width=0.75]    (0,5.59) -- (0,-5.59)   ;
\draw [fill={rgb, 255:red, 166; green, 194; blue, 222 }  ,fill opacity=1 ]   (558.5,164.5) -- (439,164.5) ;
\draw [shift={(439,164.5)}, rotate = 360] [color={rgb, 255:red, 0; green, 0; blue, 0 }  ][line width=0.75]    (0,5.59) -- (0,-5.59)   ;
\draw [shift={(558.5,164.5)}, rotate = 360] [color={rgb, 255:red, 0; green, 0; blue, 0 }  ][line width=0.75]    (0,5.59) -- (0,-5.59)   ;
\draw [fill={rgb, 255:red, 166; green, 194; blue, 222 }  ,fill opacity=1 ]   (431.5,44.5) -- (461.5,14.5) ;
\draw [shift={(461.5,14.5)}, rotate = 135] [color={rgb, 255:red, 0; green, 0; blue, 0 }  ][line width=0.75]    (0,5.59) -- (0,-5.59)   ;
\draw [shift={(431.5,44.5)}, rotate = 135] [color={rgb, 255:red, 0; green, 0; blue, 0 }  ][line width=0.75]    (0,5.59) -- (0,-5.59)   ;

\draw  [fill={rgb, 255:red, 166; green, 194; blue, 222 }  ,fill opacity=1 ][line width=1.5]  (227,48.85) -- (232.85,43) -- (354,43) -- (354,149.15) -- (348.15,155) -- (227,155) -- cycle ; \draw  [line width=1.5]  (354,43) -- (348.15,48.85) -- (227,48.85) ; \draw  [line width=1.5]  (348.15,48.85) -- (348.15,155) ;
\draw    (218,49) -- (218,156) ;
\draw [shift={(218,156)}, rotate = 270] [color={rgb, 255:red, 0; green, 0; blue, 0 }  ][line width=0.75]    (0,5.59) -- (0,-5.59)   ;
\draw [shift={(218,49)}, rotate = 270] [color={rgb, 255:red, 0; green, 0; blue, 0 }  ][line width=0.75]    (0,5.59) -- (0,-5.59)   ;
\draw    (346.5,162.5) -- (227,162.5) ;
\draw [shift={(227,162.5)}, rotate = 360] [color={rgb, 255:red, 0; green, 0; blue, 0 }  ][line width=0.75]    (0,5.59) -- (0,-5.59)   ;
\draw [shift={(346.5,162.5)}, rotate = 360] [color={rgb, 255:red, 0; green, 0; blue, 0 }  ][line width=0.75]    (0,5.59) -- (0,-5.59)   ;
\draw  [fill={rgb, 255:red, 166; green, 194; blue, 222 }  ,fill opacity=1 ][line width=1.5]  (124,49) -- (129,44) -- (138,44) -- (138,151) -- (133,156) -- (124,156) -- cycle ; \draw  [line width=1.5]  (138,44) -- (133,49) -- (124,49) ; \draw  [line width=1.5]  (133,49) -- (133,156) ;
\draw    (115.5,49) -- (115.5,156) ;
\draw [shift={(115.5,156)}, rotate = 270] [color={rgb, 255:red, 0; green, 0; blue, 0 }  ][line width=0.75]    (0,5.59) -- (0,-5.59)   ;
\draw [shift={(115.5,49)}, rotate = 270] [color={rgb, 255:red, 0; green, 0; blue, 0 }  ][line width=0.75]    (0,5.59) -- (0,-5.59)   ;

\draw (409.5,95.4) node [anchor=north west][inner sep=0.75pt]    {$d_{1}$};
\draw (497,167.4) node [anchor=north west][inner sep=0.75pt]    {$d_{2}$};
\draw (429,11.4) node [anchor=north west][inner sep=0.75pt]    {$d_{3}$};
\draw (98,95.4) node [anchor=north west][inner sep=0.75pt]    {$d_{1}$};
\draw (201.5,95.9) node [anchor=north west][inner sep=0.75pt]    {$d_{1}$};
\draw (281,169.4) node [anchor=north west][inner sep=0.75pt]    {$d_{2}$};

\end{tikzpicture}
}
    \caption{Visualization of a vector $v \in \mathbb{R}^{d_1}$ (left), a matrix $M \in \mathbb{R}^{d_1 \times d_2}$ (middle), and a tensor $\mathcal{T} \in \mathbb{R}^{d_1 \times d_2 \times d_3}$ (right).}
    \label{fig:tensor_basic}
\end{figure}

To generate third-order tensors from a sequence of matrices $M_i \in \mathbb{R}^{d_2 \times d_3}$ for $i=1,...,d_1$, we define the function $\text{fold} \left( \{ M_i \}_{i = 1}^{d_1} \right) \in \mathbb{R}^{d_1 \times d_2 \times d_3}$, which horizontally stacks matrices from top to bottom, as visualized in Figure~\ref{fig:tensor_operations}. It is often useful to unfold tensors as well. Unfolding is defined mode-wise and we denote $T_{(i)} \in \mathbb{R}^{d_i \times \prod_{j=1, j\neq i}^n d_j}$ the mode-$i$ unfolding of $\mathcal{T}$ whose columns are comprised of vectors obtained by fixing all but the $i^{th}$ index. 

\begin{figure}[h]
    \centering
    \tikzset{every picture/.style={line width=0.75pt}} 

\begin{tikzpicture}[x=0.75pt,y=0.75pt,yscale=-1,xscale=1]

\draw  [fill={rgb, 255:red, 186; green, 202; blue, 218 }  ,fill opacity=1 ][line width=1.5]  (533.04,96.4) -- (554.81,74.62) -- (626.36,74.62) -- (626.36,83.96) -- (604.58,105.73) -- (533.04,105.73) -- cycle ; \draw  [line width=1.5]  (626.36,74.62) -- (604.58,96.4) -- (533.04,96.4) ; \draw  [line width=1.5]  (604.58,96.4) -- (604.58,105.73) ;
\draw [fill={rgb, 255:red, 166; green, 194; blue, 222 }  ,fill opacity=1 ]   (527,50.58) -- (527,106.58) ;
\draw [shift={(527,106.58)}, rotate = 270] [color={rgb, 255:red, 0; green, 0; blue, 0 }  ][line width=0.75]    (0,5.59) -- (0,-5.59)   ;
\draw [shift={(527,50.58)}, rotate = 270] [color={rgb, 255:red, 0; green, 0; blue, 0 }  ][line width=0.75]    (0,5.59) -- (0,-5.59)   ;
\draw [fill={rgb, 255:red, 166; green, 194; blue, 222 }  ,fill opacity=1 ]   (625,19.5) -- (555.04,19.26) ;
\draw [shift={(555.04,19.26)}, rotate = 0.19] [color={rgb, 255:red, 0; green, 0; blue, 0 }  ][line width=0.75]    (0,5.59) -- (0,-5.59)   ;
\draw [shift={(625,19.5)}, rotate = 0.19] [color={rgb, 255:red, 0; green, 0; blue, 0 }  ][line width=0.75]    (0,5.59) -- (0,-5.59)   ;
\draw  [fill={rgb, 255:red, 144; green, 185; blue, 225 }  ,fill opacity=1 ][line width=1.5]  (533.04,87.07) -- (554.81,65.29) -- (626.36,65.29) -- (626.36,74.62) -- (604.58,96.4) -- (533.04,96.4) -- cycle ; \draw  [line width=1.5]  (626.36,65.29) -- (604.58,87.07) -- (533.04,87.07) ; \draw  [line width=1.5]  (604.58,87.07) -- (604.58,96.4) ;
\draw  [fill={rgb, 255:red, 129; green, 173; blue, 222 }  ,fill opacity=1 ][line width=1.5]  (533.04,77.73) -- (554.81,55.96) -- (626.36,55.96) -- (626.36,65.29) -- (604.58,87.07) -- (533.04,87.07) -- cycle ; \draw  [line width=1.5]  (626.36,55.96) -- (604.58,77.73) -- (533.04,77.73) ; \draw  [line width=1.5]  (604.58,77.73) -- (604.58,87.07) ;
\draw  [fill={rgb, 255:red, 108; green, 165; blue, 215 }  ,fill opacity=1 ][line width=1.5]  (533.04,68.4) -- (554.81,46.63) -- (626.36,46.63) -- (626.36,55.96) -- (604.58,77.73) -- (533.04,77.73) -- cycle ; \draw  [line width=1.5]  (626.36,46.63) -- (604.58,68.4) -- (533.04,68.4) ; \draw  [line width=1.5]  (604.58,68.4) -- (604.58,77.73) ;
\draw  [fill={rgb, 255:red, 74; green, 147; blue, 220 }  ,fill opacity=1 ][line width=1.5]  (533.04,59.07) -- (554.81,37.3) -- (626.36,37.3) -- (626.36,46.63) -- (604.58,68.4) -- (533.04,68.4) -- cycle ; \draw  [line width=1.5]  (626.36,37.3) -- (604.58,59.07) -- (533.04,59.07) ; \draw  [line width=1.5]  (604.58,59.07) -- (604.58,68.4) ;
\draw  [fill={rgb, 255:red, 24; green, 121; blue, 218 }  ,fill opacity=1 ][line width=1.5]  (533.04,49.74) -- (554.81,27.96) -- (626.36,27.96) -- (626.36,37.3) -- (604.58,59.07) -- (533.04,59.07) -- cycle ; \draw  [line width=1.5]  (626.36,27.96) -- (604.58,49.74) -- (533.04,49.74) ; \draw  [line width=1.5]  (604.58,49.74) -- (604.58,59.07) ;
\draw  [fill={rgb, 255:red, 24; green, 121; blue, 218 }  ,fill opacity=1 ][line width=1.5]  (46,45.26) -- (49.26,42) -- (116.87,42) -- (116.87,101.24) -- (113.61,104.5) -- (46,104.5) -- cycle ; \draw  [line width=1.5]  (116.87,42) -- (113.61,45.26) -- (46,45.26) ; \draw  [line width=1.5]  (113.61,45.26) -- (113.61,104.5) ;
\draw  [fill={rgb, 255:red, 74; green, 147; blue, 220 }  ,fill opacity=1 ][line width=1.5]  (126,45.26) -- (129.26,42) -- (196.87,42) -- (196.87,101.24) -- (193.61,104.5) -- (126,104.5) -- cycle ; \draw  [line width=1.5]  (196.87,42) -- (193.61,45.26) -- (126,45.26) ; \draw  [line width=1.5]  (193.61,45.26) -- (193.61,104.5) ;
\draw  [fill={rgb, 255:red, 144; green, 185; blue, 225 }  ,fill opacity=1 ][line width=1.5]  (252,43.26) -- (255.26,40) -- (322.87,40) -- (322.87,99.24) -- (319.61,102.5) -- (252,102.5) -- cycle ; \draw  [line width=1.5]  (322.87,40) -- (319.61,43.26) -- (252,43.26) ; \draw  [line width=1.5]  (319.61,43.26) -- (319.61,102.5) ;
\draw  [fill={rgb, 255:red, 186; green, 202; blue, 218 }  ,fill opacity=1 ][line width=1.5]  (332,43.26) -- (335.26,40) -- (402.87,40) -- (402.87,99.24) -- (399.61,102.5) -- (332,102.5) -- cycle ; \draw  [line width=1.5]  (402.87,40) -- (399.61,43.26) -- (332,43.26) ; \draw  [line width=1.5]  (399.61,43.26) -- (399.61,102.5) ;
\draw [fill={rgb, 255:red, 166; green, 194; blue, 222 }  ,fill opacity=1 ]   (528,42.58) -- (548,22.58) ;
\draw [shift={(548,22.58)}, rotate = 135] [color={rgb, 255:red, 0; green, 0; blue, 0 }  ][line width=0.75]    (0,5.59) -- (0,-5.59)   ;
\draw [shift={(528,42.58)}, rotate = 135] [color={rgb, 255:red, 0; green, 0; blue, 0 }  ][line width=0.75]    (0,5.59) -- (0,-5.59)   ;
\draw  [fill={rgb, 255:red, 0; green, 0; blue, 0 }  ,fill opacity=1 ] (436,70.64) -- (467.2,70.64) -- (467.2,65) -- (488,76.29) -- (467.2,87.58) -- (467.2,81.93) -- (436,81.93) -- cycle ;

\draw (509.29,69.02) node [anchor=north west][inner sep=0.75pt]  [font=\footnotesize]  {$d_{1}$};
\draw (581.72,2.82) node [anchor=north west][inner sep=0.75pt]  [font=\footnotesize]  {$d_{2}$};
\draw (519.42,14.77) node [anchor=north west][inner sep=0.75pt]  [font=\footnotesize]  {$d_{3}$};
\draw (210,80.4) node [anchor=north west][inner sep=0.75pt]  [font=\huge]  {$...$};
\draw (68,66.4) node [anchor=north west][inner sep=0.75pt]    {$M_{1}$};
\draw (147,67.4) node [anchor=north west][inner sep=0.75pt]    {$M_{2}$};
\draw (268,65.4) node [anchor=north west][inner sep=0.75pt]    {$M_{d_{1} -1}$};
\draw (355,62.4) node [anchor=north west][inner sep=0.75pt]    {$M_{d_{1}}$};
\draw (566.52,45.06) node [anchor=north west][inner sep=0.75pt]  [font=\scriptsize,rotate=-309.53,xslant=-0.84]  {$M_{1}$};

\end{tikzpicture}
    \caption{Visualization of the $\text{fold}(\cdot)$ operation applied to matrices $M_1$, ..., $M_{d_1}$.}
    \label{fig:tensor_operations}
\end{figure}
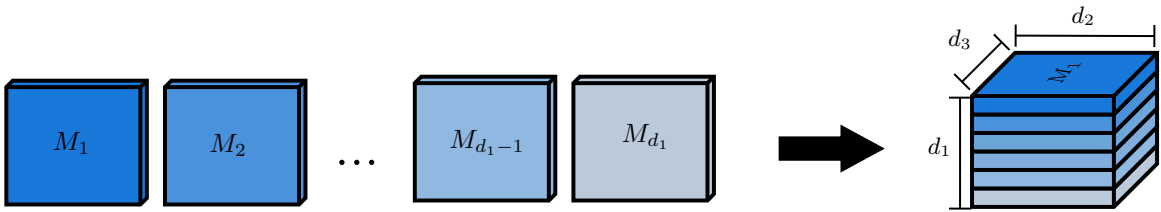

A third order tensor $\mathcal{T} \in \mathbb{R}^{ d_1 \times d_2 \times d_3}$ is said to have CANDECOMP-PARAFAC rank (CP-rank) $R$ if $R$ is the smallest integer such that $\mathcal{T}$ can be decomposed into the sum of $R$ rank-one tensors:
\begin{equation} 
\mathcal{T} = \sum_{i=1}^R \lambda_i (a_i \otimes b_i \otimes c_i),
\label{eq:cp}
\end{equation}
where $a_i \in \mathbb{R}^{d_1}$, $b_i \in \mathbb{R}^{d_2}$, and $c_i \in \mathbb{R}^{d_3}$ for $i=1,...,R$ and $\lambda = \begin{bmatrix} \lambda_1 & \cdots & \lambda_R \end{bmatrix}^T \in \mathbb{R}^R$; see, e.g.,~\cite{kiers2000towards,kolda2009tensor}. The column-wise concatenation of the vectors $a_i$, $b_i$, and $c_i$ define the \emph{factor matrices} $A \in \mathbb{R}^{d_1 \times R}$, $B \in \mathbb{R}^{d_2 \times R}$ and $C \in \mathbb{R}^{d_3 \times R}$, and we represent the CP-decomposition in short-hand as $\mathcal{T} = [[\lambda; A, B, C]]$ for notational convenience. If all $\lambda_i$'s are equal to $1$, we suppress the notation further to  $\mathcal{T} = [[A, B, C]]$.

\begin{figure}[h]
    \centering
    \tikzset{every picture/.style={line width=0.75pt}} 

\begin{tikzpicture}[x=0.75pt,y=0.75pt,yscale=-1,xscale=1]

\draw  [fill={rgb, 255:red, 166; green, 194; blue, 222 }  ,fill opacity=1 ][line width=0.75]  (11,55.05) -- (18.75,47.3) -- (61.8,47.3) -- (61.8,97.2) -- (54.05,104.95) -- (11,104.95) -- cycle ; \draw  [line width=0.75]  (61.8,47.3) -- (54.05,55.05) -- (11,55.05) ; \draw  [line width=0.75]  (54.05,55.05) -- (54.05,104.95) ;
\draw  [fill={rgb, 255:red, 166; green, 194; blue, 222 }  ,fill opacity=1 ][line width=0.75]  (111.72,57.18) -- (113.86,55.05) -- (117.7,55.05) -- (117.7,102.81) -- (115.56,104.95) -- (111.72,104.95) -- cycle ; \draw  [line width=0.75]  (117.7,55.05) -- (115.56,57.18) -- (111.72,57.18) ; \draw  [line width=0.75]  (115.56,57.18) -- (115.56,104.95) ;
\draw  [fill={rgb, 255:red, 166; green, 194; blue, 222 }  ,fill opacity=1 ][line width=0.75]  (162.65,54.8) -- (164.84,52.61) -- (164.84,49.27) -- (123.33,49.27) -- (121.13,51.47) -- (121.13,54.8) -- cycle ; \draw  [line width=0.75]  (164.84,49.27) -- (162.65,51.47) -- (162.65,54.8) ; \draw  [line width=0.75]  (162.65,51.47) -- (121.13,51.47) ;
\draw  [fill={rgb, 255:red, 166; green, 194; blue, 222 }  ,fill opacity=1 ][line width=0.75]  (118.95,48.64) -- (126.07,41.38) -- (126.04,38.44) -- (121.93,38.48) -- (114.81,45.74) -- (114.84,48.68) -- cycle ; \draw  [line width=0.75]  (126.04,38.44) -- (118.92,45.7) -- (118.95,48.64) ; \draw  [line width=0.75]  (118.92,45.7) -- (114.81,45.74) ;

\draw  [fill={rgb, 255:red, 166; green, 194; blue, 222 }  ,fill opacity=1 ][line width=0.75]  (210.01,57.18) -- (212.15,55.05) -- (215.99,55.05) -- (215.99,102.81) -- (213.86,104.95) -- (210.01,104.95) -- cycle ; \draw  [line width=0.75]  (215.99,55.05) -- (213.86,57.18) -- (210.01,57.18) ; \draw  [line width=0.75]  (213.86,57.18) -- (213.86,104.95) ;
\draw  [fill={rgb, 255:red, 166; green, 194; blue, 222 }  ,fill opacity=1 ][line width=0.75]  (260.94,54.8) -- (263.14,52.61) -- (263.14,49.27) -- (221.62,49.27) -- (219.43,51.47) -- (219.43,54.8) -- cycle ; \draw  [line width=0.75]  (263.14,49.27) -- (260.94,51.47) -- (260.94,54.8) ; \draw  [line width=0.75]  (260.94,51.47) -- (219.43,51.47) ;
\draw  [fill={rgb, 255:red, 166; green, 194; blue, 222 }  ,fill opacity=1 ][line width=0.75]  (217.25,48.64) -- (224.36,41.38) -- (224.33,38.44) -- (220.22,38.48) -- (213.1,45.74) -- (213.13,48.68) -- cycle ; \draw  [line width=0.75]  (224.33,38.44) -- (217.22,45.7) -- (217.25,48.64) ; \draw  [line width=0.75]  (217.22,45.7) -- (213.1,45.74) ;

\draw  [fill={rgb, 255:red, 166; green, 194; blue, 222 }  ,fill opacity=1 ][line width=0.75]  (346.85,57.18) -- (348.98,55.05) -- (352.83,55.05) -- (352.83,102.81) -- (350.69,104.95) -- (346.85,104.95) -- cycle ; \draw  [line width=0.75]  (352.83,55.05) -- (350.69,57.18) -- (346.85,57.18) ; \draw  [line width=0.75]  (350.69,57.18) -- (350.69,104.95) ;
\draw  [fill={rgb, 255:red, 166; green, 194; blue, 222 }  ,fill opacity=1 ][line width=0.75]  (397.78,54.8) -- (399.97,52.61) -- (399.97,49.27) -- (358.45,49.27) -- (356.26,51.47) -- (356.26,54.8) -- cycle ; \draw  [line width=0.75]  (399.97,49.27) -- (397.78,51.47) -- (397.78,54.8) ; \draw  [line width=0.75]  (397.78,51.47) -- (356.26,51.47) ;
\draw  [fill={rgb, 255:red, 166; green, 194; blue, 222 }  ,fill opacity=1 ][line width=0.75]  (354.08,48.64) -- (361.19,41.38) -- (361.17,38.44) -- (357.05,38.48) -- (349.94,45.74) -- (349.97,48.68) -- cycle ; \draw  [line width=0.75]  (361.17,38.44) -- (354.05,45.7) -- (354.08,48.64) ; \draw  [line width=0.75]  (354.05,45.7) -- (349.94,45.74) ;

\draw  [fill={rgb, 255:red, 166; green, 194; blue, 222 }  ,fill opacity=1 ][line width=0.75]  (553.34,55.8) -- (562.09,47.05) -- (583,47.05) -- (583,68.92) -- (574.25,77.67) -- (553.34,77.67) -- cycle ; \draw  [line width=0.75]  (583,47.05) -- (574.25,55.8) -- (553.34,55.8) ; \draw  [line width=0.75]  (574.25,55.8) -- (574.25,77.67) ;
\draw  [fill={rgb, 255:red, 166; green, 194; blue, 222 }  ,fill opacity=1 ][line width=0.75]  (585.57,41.05) -- (607.11,19.5) -- (607.11,15.33) -- (587.88,15.33) -- (566.33,36.88) -- (566.33,41.05) -- cycle ; \draw  [line width=0.75]  (607.11,15.33) -- (585.57,36.88) -- (585.57,41.05) ; \draw  [line width=0.75]  (585.57,36.88) -- (566.33,36.88) ;
\draw  [fill={rgb, 255:red, 166; green, 194; blue, 222 }  ,fill opacity=1 ][line width=0.75]  (632.45,68.96) -- (634.31,67.05) -- (634.11,47.05) -- (589.41,47.48) -- (587.55,49.38) -- (587.75,69.39) -- cycle ; \draw  [line width=0.75]  (634.11,47.05) -- (632.25,48.95) -- (632.45,68.96) ; \draw  [line width=0.75]  (632.25,48.95) -- (587.55,49.38) ;
\draw  [fill={rgb, 255:red, 166; green, 194; blue, 222 }  ,fill opacity=1 ][line width=0.75]  (525,49.86) -- (527.82,47.05) -- (547.62,47.05) -- (547.62,102.13) -- (544.8,104.95) -- (525,104.95) -- cycle ; \draw  [line width=0.75]  (547.62,47.05) -- (544.8,49.86) -- (525,49.86) ; \draw  [line width=0.75]  (544.8,49.86) -- (544.8,104.95) ;
\draw  [fill={rgb, 255:red, 166; green, 194; blue, 222 }  ,fill opacity=1 ][line width=0.75]  (440.67,55.05) -- (448.42,47.3) -- (491.46,47.3) -- (491.46,97.2) -- (483.71,104.95) -- (440.67,104.95) -- cycle ; \draw  [line width=0.75]  (491.46,47.3) -- (483.71,55.05) -- (440.67,55.05) ; \draw  [line width=0.75]  (483.71,55.05) -- (483.71,104.95) ;

\draw (69.16,73.57) node [anchor=north west][inner sep=0.75pt]  [font=\normalsize]  {$\approx $};
\draw (26.79,70.57) node [anchor=north west][inner sep=0.75pt]  [font=\normalsize]  {$\mathcal{T}$};
\draw (225.21,29.39) node [anchor=north west][inner sep=0.75pt]  [font=\footnotesize]  {$c_{2}$};
\draw (264.15,44.74) node [anchor=north west][inner sep=0.75pt]  [font=\footnotesize]  {$b_{2}$};
\draw (205.1,106.15) node [anchor=north west][inner sep=0.75pt]  [font=\footnotesize]  {$a_{2}$};
\draw (361.8,29.39) node [anchor=north west][inner sep=0.75pt]  [font=\footnotesize]  {$c_{R}$};
\draw (400.6,44.74) node [anchor=north west][inner sep=0.75pt]  [font=\footnotesize]  {$b_{R}$};
\draw (341.55,106.15) node [anchor=north west][inner sep=0.75pt]  [font=\footnotesize]  {$a_{R}$};
\draw (162.09,70.57) node [anchor=north west][inner sep=0.75pt]  [font=\normalsize]  {$+$};
\draw (267.69,70.57) node [anchor=north west][inner sep=0.75pt]  [font=\normalsize]  {${\displaystyle +\ ...+}$};
\draw (184.61,69.57) node [anchor=north west][inner sep=0.75pt]  [font=\normalsize]  {$\lambda _{2}$};
\draw (318.86,69.57) node [anchor=north west][inner sep=0.75pt]  [font=\normalsize]  {$\lambda _{R}$};
\draw (455.65,70.57) node [anchor=north west][inner sep=0.75pt]  [font=\normalsize]  {$\mathcal{T}$};
\draw (497.8,69.37) node [anchor=north west][inner sep=0.75pt]  [font=\normalsize]  {$\approx $};
\draw (559.34,60.07) node [anchor=north west][inner sep=0.75pt]  [font=\footnotesize]  {$\mathcal{G}$};
\draw (527.14,71.57) node [anchor=north west][inner sep=0.75pt]  [font=\footnotesize]  {$A$};
\draw (603.34,53.45) node [anchor=north west][inner sep=0.75pt]  [font=\footnotesize]  {$B$};
\draw (584.26,18.07) node [anchor=north west][inner sep=0.75pt]  [font=\footnotesize,rotate=-21.19]  {$C$};
\draw (91.72,69.57) node [anchor=north west][inner sep=0.75pt]  [font=\normalsize]  {$\lambda _{1}$};
\draw (106.8,106.15) node [anchor=north west][inner sep=0.75pt]  [font=\footnotesize]  {$a_{1}$};
\draw (165.86,44.74) node [anchor=north west][inner sep=0.75pt]  [font=\footnotesize]  {$b_{1}$};
\draw (126.91,29.39) node [anchor=north west][inner sep=0.75pt]  [font=\footnotesize]  {$c_{1}$};

\end{tikzpicture}
    \caption{Visualization of the CP decomposition (left) and Tucker decomposition (right) of the tensor $\mathcal{T}$.}
    \label{fig:decomp}
\end{figure}
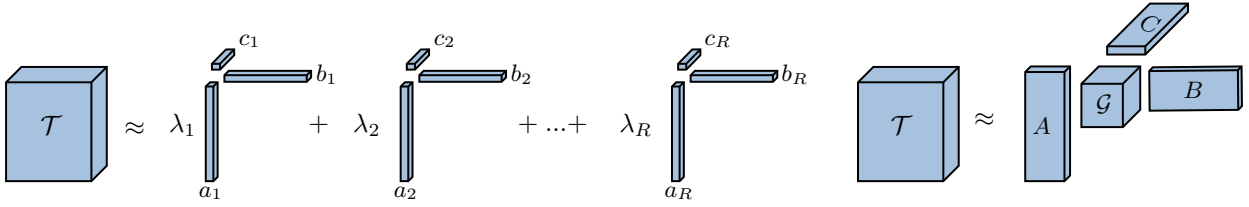

The tuple rank $(R_1, R_2, R_3)$ Tucker Decomposition~\cite{tucker1966some,kolda2009tensor} is defined as 
\begin{equation}
    \mathcal{T} = \sum_{i=1}^{R_1}\sum_{j=1}^{R_2} \sum_{k=1}^{R_3} g_{ijk} (a_i \otimes b_j \otimes c_k)
    \label{eq:tucker}
\end{equation}
and can also be written using mode-wise products:
\begin{equation}
    \mathcal{T} = \mathcal{G} \times_1 A \times_2 B \times_3 C,
\end{equation}
where $\mathcal{G} \in \mathbb{R}^{R_1 \times R_2 \times R_3}$, $A \in \mathbb{R}^{d_1 \times R_1}$,  $B \in \mathbb{R}^{d_2 \times R_2}$, and  $C \in \mathbb{R}^{d_3 \times R_3}$. The shorthand notation for the Tucker decomposition will be written as $\mathcal{T} = [[\mathcal{G}; A, B, C]]$. Note that we often refer to the tuple rank $(R,R,R)$ Tucker decomposition simply as a rank-$R$ Tucker decomposition. Visualizations of the CP-decomposition and Tucker decomposition are shown in Figure~\ref{fig:decomp}.

\subsection{Tensorized Parameters}
\label{subsec:tensorized_attention}

We first explore how the matrix-valued quantities and models of Section~\ref{sec:original_notation} can be represented as (higher-order) tensors, and why doing so can be advantageous.
We recall that standard attention involves embeddings, or weight matrices, that map the inputs into query, key, and value matrices in each attention head within a given layer of the transformer architecture. 
In MHA, for example, the final output of the attention model is a concatenation of the attention outputs across all heads, mapped back to the original model dimensions (see Section~\ref{sec:multi-head}). 
Rather than matrix-based concatenations, we can instead represent attention computations in terms of higher-order tensors.

By imposing tensor structure on standard attention matrix quantities, we aim to alleviate several computational challenges.
For example, with limited context windows (cf.~causality in Section~\ref{sec:original_notation}), standard attention models can fail to capture long-range token dependencies, particularly in long sequences; tensorized attention models may reduce these to short-range interactions within context windows along another dimension, cf. \cite{feng2024long}.
Additionally, standard attention is known to struggle with KV cache overload
and is limited in expressivity, as only two tokens are considered at a time, vs. higher-order relationships (e.g., triplets of tokens).
Tensor representations, on the other hand, may be used for compressing weights or queries, keys, and values, offering better performance and reduced storage requirements while enhancing expressivity and perplexity.

In this section, we summarize and highlight examples of works that have employed tensorization within pre-trained transformer models for the purposes of compression. We begin with a discussion on tensorization of weight matrices~\cite{cordonnier2020multi, ren2022exploring, luo2025trawl} and then segue into the tensorization of query, key, and value matrices~\cite{zhang2025tensor}.

Tensor-based representations of the weight matrices in standard attention are motivated by the observation that the weights of a model can be seen as having more than two dimensions, e.g., 
$$\tensor{W}_{\ell}^Q = \fold{ \{\head{\WQ}{h} \}_{h \in \Nh}  } \in \RR^{\Nh \times \dhead \times d}$$ 
corresponding to the query weight matrix in the $h$-th head of the $\ell$-th layer.
We can also further tensorize, e.g.,
\begin{align*}
    \tensor{W}^Q \in \RR^{L \times \Nh \times \dhead \times d} \textup{ given by } \tensor{W}^Q(\ell, :, :, :) &= \tensor{W}_{\ell}^Q, \ \ell = 1,\ldots,L,
\end{align*}
or $\tensor{W} \in \RR^{4 \times L \times \Nh \times \dhead \times d}$ corresponding to the tensor comprised of $\tensor{W}^Q, \tensor{W}^K, \tensor{W}^V,$ and $\tensor{W}^O$.
Methods that tensorize the weight matrices attempt to find higher-order correlations among some subset of these dimensions to achieve compression or speed-up of the associated model, based on a reasonable assumption that the matrices between different heads of the same layer of the same type ($Q, K, V$, or $O$) will be highly correlated.

We discuss three works in more detail. We begin with \cite{cordonnier2020multi}, one of the first works to propose stacking matrices from the attention computation (specifically, across heads) into a third-order tensor and learning this tensor via classical tensor decompositions, thereby capturing potential similarities across heads. We then turn to \cite{ren2022exploring}, which extends this idea by stacking a larger collection of related weight matrices -- including multiple attention matrices within a layer, matrices across different layers, and feed-forward parameters -- into higher-order tensors and applying a variety of tensor decompositions. 
Finally, we conclude with \cite{zhang2025tensor}, a recent and influential paper that does not directly employ classical tensor decompositions but instead builds on the idea of approximating attention with low-rank structure across multiple dimensions.

\subsubsection{Concatenation vs.\ Collaboration in Multi-Head Attention} 
\label{subsubsec:collaboration}

The work in \cite{cordonnier2020multi} aims to better detect and quantify attention head redundancy by asking whether independent heads learn overlapping or distinct concepts. This is accomplished via a rank-$R$ tensor decomposition of 
  \begin{equation}
      \mathcal{X} := \fold{\{\head{\WQ}{h}(\head{\WK}{h})^T\}_{h \in \Nh}} \in \mathbb{R}^{\Nh \times d \times d}.
      \label{eq:tensor_cordonnier}
  \end{equation}
 In particular, the CP-decomposition of \eqref{eq:tensor_cordonnier}
 $$\mathcal{X} \approx \llbracket M, \tilde{W}_Q, \tilde{W}_K\rrbracket$$ 
 is used, where $M \in \mathbb{R}^{\Nh \times R}$, $\tilde{W}_Q \in \mathbb{R}^{d \times R}$, and $\tilde{W}_K^T \in \mathbb{R}^{d \times R}$. These factors define joint weight matrices $\tilde{W}_Q, \tilde{W}_K$ and a ``mixing matrix'' $M$ which are used to define the collaborative attention head. 
Intuitively, this approach aims to obtain a low-rank approximation of $\head{\WQ}{h}(\head{\WK}{h})^T$ for each head, using weights across all heads simultaneously. That is, for each $h \in [\Nh]$, we have
\begin{equation}
\head{\WQ}{h}(\head{\WK}{h})^T \approx \tilde{W}_Q M_h \tilde{W}_K^T,
    \label{eq:eq:cordonnier_interp}
\end{equation}
where the factors $\tilde{W}_Q$, $\tilde{W}_K^T $, and the diagonal matrix $M_h := \text{diag}({M}(h,:)) \in \mathbb{R}^{R \times R}$ built from the $h^{th}$ row of $M$, are all obtained via CP-decomposition.

In practice, after the CP-decomposition is computed, the matrices $\head{\WK}{h}$ and $\head{\WQ}{h}$ are replaced with the joint matrix $\head{\WK}{h} \to \tilde{W}_K$ and by correlated matrices $\head{\WQ}{h} \to \tilde{W}_Q M_h$. 
As a result, instead of forming scalar scores via dot products based on $Q$ and $K$, to obtain scores $S(i,j)=\langle Q(i,:),K(j,:)\rangle$,
one introduces a feature-wise score tensor
\[
\tilde{\mathcal{S}}(i,j,\ell)=\tilde{Q}(i,\ell)\,\tilde{K}(j,\ell),
\]
where $\tilde{\mathcal{S}} \in \RR^{N\times N\times \dhead}$ encodes token--token interactions separately for each feature dimension~\cite{shen2019tensorized}.
Using this approach, numerical results demonstrate that the size of query and key projections can be significantly reduced, while preserving speed and accuracy~\cite{cordonnier2020multi}. Note that this reduces the computational cost in computing the attention matrix from $\mathcal{O}(N^2 \dhead \Nh)$ to $\mathcal{O}(R N \dhead \Nh)$ where $R$ is the CP-rank.

\subsubsection{Exploring extreme parameter compression for pre-trained language models}
\label{subsubsec:extreme-compression}

The idea of learning a joint representation over heads was further developed in several works, including \cite{ren2022exploring} and \cite{luo2025trawl}, where the proposed approach was to stack all trainable weight matrices across all layers into a third-order tensor. Although the model weights are generally not low-rank, cf. \cite{YuWu2023CompressingTransformers}, these methods demonstrate that significant performance improvements can be achieved by compressing weight matrices without adversely affecting model inference, especially when the low-rank representation is informed by the higher-order tensor rank.

A more general version of this can be found in \cite{ren2022exploring}, which takes the approach of, assuming $\dhead = d$,
$$\mathcal{X} := \text{fold}(\{W_i \}_{i=1}^{m})  \in \mathbb{R}^{m \times d \times d},$$
where $m = 12 L$ and $L$ is the number of layers, with $12$ matrices per layer. The 12 matrices correspond to the standard query, key, value, and output weight matrices, $\head{\WQ}{h},\head{\WK}{h}, \head{\WV}{h}, W^O_h$ in addition to 4 input weight matrices and 4 output weight matrices from the feedforward network. In~\cite{luo2025trawl}, $\mathcal{X} \in \mathbb{R}^{4\ell \times d \times d}$, using only the weight query, key, value, and output weight matrices ($\ell = 1$ if compressing layer-wise and $\ell=L$ if compressing across all layers). Given the tensor $\mathcal{X}$, three different low-rank compression models are considered in~\cite{ren2022exploring}: (i) matrix-based decomposition, (ii) the shared factor model, and (iii) the Tucker decomposition, for which we give some more details in this section.

\begin{itemize}
\item[(i)] The matrix-based decomposition approach uses the rank-$R$ truncated SVD of each weight matrix $\mathcal{X}(i,:,:) =: X_i \in \mathbb{R}^{d \times d}$:
\begin{equation}
    X_i \approx U_i\Sigma_i V_i^T,
    \label{eq:ren_model1}
\end{equation}
where $U_i, V_i \in \mathbb{R}^{d \times R}$ and $\Sigma_i \in \mathbb{R}^{R \times R}$. This approach is equivalent to what has been done in other works, where cross-matrix redundancy is not exploited (see Section~\ref{sec:low-rank-approximation}). 

\item[(ii)] The shared factor approach assumes there are common low-rank embedding matrices $B, C \in \mathbb{R}^{d \times R}$ between all weight matrices in $\mathcal{X}$. In other words, for $i \in [m]$:
\begin{equation}
    X_i \approx B G_i C^T,
    \label{eq:ren_model2}
\end{equation}
where $G_i \in \mathbb{R}^{R \times R}$. It is important to emphasize that $B$ and $C$ are learned and shared across all weight matrices of which $\mathcal{X}$ is comprised. This can also be seen as a special case of the Tucker decomposition~\eqref{eq:tucker},  
\begin{equation}
    \mathcal{X} \approx \mathcal{G} \times_1 I \times_2 B \times_3 C,
    \label{eq:ren_model2_tensor}
\end{equation}
where $I$ is an $m\times m$ identity matrix, $\mathcal{G} \in \mathbb{R}^{m \times R \times R}$ and each horizontal slice of $\mathcal{G}$ corresponds to $\mathcal{G}(i,:,:) = G_i$ of~\eqref{eq:ren_model2}. 

\item[(iii)] The shared factor model is generalized by the full Tucker decomposition,
\begin{equation}
\mathcal{X} \approx \mathcal{G} \times_1 A \times_2 B \times_3 C.
\label{eq:ren_model3_tensor}
\end{equation}
Here, $A \in\RR^{m\times \ell}$ contains block-specific coefficients, and $\mathcal{G}\in\RR^{\ell\times R\times R}$ is a shared core tensor whose slices encode a bank of $\ell$ basis matrices in $\RR^{R\times R}$. Equivalently, for each block index $i\in\{1,\dots,m\}$, the corresponding weight matrix is approximated via low-rank representation:
\begin{equation}
X_i \approx  B\left(\sum_{k=1}^{\ell} A(i,k)\,G_i\right)C^{\top},
\tag{F}
\label{eq:ren_model3}
\end{equation}
where again $\mathcal{G}(i,:,:) = G_i$.
Thus, each weight matrix is expressed as a bilinear map with shared left and right factors, while its intermediate representation is a linear combination of a small number of shared core matrices.

\end{itemize}

Using these different low-rank approximation approaches demonstrated experimental success~\cite {ren2022exploring}. For example, in comparison to BERT-base~\cite{devlin2019bert}, which originally used $86$M parameters, using the (ii) shared model factor~\ref{eq:ren_model2_tensor} and (iii) Tucker Decomposition~\ref{eq:ren_model3_tensor}, only $1.8$M and $1.9$M parameters are needed, respectively. Despite the significant decrease in memory storage and thus inference computational cost, the empirical average accuracy across test sets in GLUE~\cite{wang2019glue} is $83\%$ for the baseline, $80\%$ for (ii), and $80.8\%$ and (iii). For more details on the numerical benefits of tensorized weights, we refer the reader to~\cite{luo2025trawl,ren2022exploring}.

\subsubsection{Tensor Product Attention} 

As an alternative to tensorizing weight matrices, the work of Tensor Product Attention (TPA)~\cite{zhang2025tensor} considers the tensorization of query, key, and value matrices directly. 
In TPA, we now consider the tensors
\begin{align}
\label{eq:TPAtensors}
    \mathcal{Q} &= \fold{ \{ Q_n \}_{n \in N} }  \in \RR^{N \times \Nh \times \dhead} \\
    \mathcal{K} &= \fold{ \{ K_n \}_{n \in N} }  \in \RR^{N \times \Nh \times \dhead} \\
    \mathcal{V} &= \fold{ \{ V_n \}_{n \in N} }  \in \RR^{N \times \Nh \times \dhead},
\end{align}
where the (per-token) slices $Q_n, K_n, V_n \in \RR^{\Nh \times \dhead}$ are factorized in terms of learned latent weight matrices.
Namely, we define
\begin{align}
    \label{eq:TPAmats}
    \Q_n &= \frac{1}{R_Q} \sum_{r = 1}^{R_Q} a_r^Q(x_n) \otimes b_r^Q(x_t), \\
    \K_n &= \frac{1}{R_K} \sum_{r = 1}^{R_K} a_r^K(x_n) \otimes b_r^K(x_t), \\
    \V_n &= \frac{1}{R_V} \sum_{r = 1}^{R_V} a_r^V(x_n) \otimes b_r^V(x_t),
\end{align}
where each pair of vectors, e.g. $a_r^Q(x_n) \in \RR^{\Nh}$ and $b_r^Q(x_n) \in \RR^{\dhead}$, depend on learned weight matrices, e.g. $\W^{a^Q}_r \in \RR^{\Nh \times d}$ and $\W^{b^Q}_r \in \RR^{\dhead \times d}$,
\begin{align*}
    a_r^Q(x_n) &= W^{a^Q}_{r} x_n, \\
    b_r^Q(x_n) &= W^{b^Q}_{r} x_n.
\end{align*}
Equivalently, we can express $Q_n$, $K_n$, and $V_n$ in matrix form, e.g.
\begin{align}
    \Q_n = \frac{1}{R_Q} A_Q(x_n)^\top B_Q(x_n),
\end{align}
where the rows of $A_Q(x_n) \in \RR^{R_Q \times \Nh}$ are given by $a_r^Q(x_n)^\top$, and the rows of $B_Q(x_n) \in \RR^{R_Q \times \dhead}$ are given by $b_r^Q(x_n)^\top$.
Repeating for all tokens, we arrive at the tensors $\mathcal{Q}, \mathcal{K}$ and $\mathcal{V}$ in (\ref{eq:TPAtensors}).
Once the tensors $\mathcal{Q}, \mathcal{K},$ and $\mathcal{V}$ are constructed in this way, slices are taken for each head, e.g. for $h=1,\ldots,\Nh$
\begin{align}
    Q_h = \mathcal{Q}(:,h,:) \in \RR^{N \times \dhead}, \\
    K_h = \mathcal{K}(:,h,:) \in \RR^{N \times \dhead}, \\
    V_h = \mathcal{Q}(:,h,:) \in \RR^{N \times \dhead}.
\end{align}
We then compute $Y_h = \textup{softmax} \left ( \frac{Q_h K_h^\top}{\sqrt{\dhead}} \right ) V_h \in \RR^{N \times {\dhead}}$ and form the output
$$ O = \left ( Y_1 | \cdots | Y_{\Nh}   \right ) \W^O \in \RR^{N \times d},$$
the usual MHA output.

In other words, TPA factorizes each token's query, key, and value matrix as a \textit{contextual} tensor product, concisely assembling information from each head and each head's dimension, as opposed to standard MHA.
Namely, in (\ref{eq:TPAmats}), the token-indexed matrices $\Q_n, \K_n, \V_n$ can have different ranks, and their representations incorporate contextual information from latent factors, e.g., $a_r^Q(x_n) \in \R^{\Nh}, b_r^Q(x_n) \in \R^{\dhead}$. 
This is particularly useful for efficient KV caching; standard attention requires $K_n$ and $V_n$ to be stored for every previous token $n$, leading to a storage cost of $2 \Nh \dhead$.
Because TPA only needs to store the latent factors, the per-token memory cost of TPA is $(R_K + R_V)(\Nh +\dhead),$ where typically $R_K$ and $R_V$ are taken to be very small, e.g. $R_K = R_V = 2$. 
The result is that TPA can handle much longer input sequences under given memory constraints than competing methods (i.e. GQA, MLA, MHA).

Additionally, within each attention layer, TPA typically involves fewer parameters than GQA or MLA. 
The complexity of TPA can be broken down into the cost of forming tensors $\mathcal{Q}, \mathcal{K}, \mathcal{V}$, given by $\Theta(N \Nh \dhead (R_Q + R_K + R_V))$, followed by the cost of MHA to compute attention.
In \cite{zhang2025tensor}, it is actually observed that MHA, MLA, and GQA may be thought of as non-contextual versions of TPA, though TPA consistently outperforms each of these competing methods in the numerical experiments of \cite{zhang2025tensor}, particularly with longer input sequences.

We note there are several modifications that can be made to the basic TPA format.
For example, RoPE can be directly incorporated  into TPA by applying standard RoPE to each query token-slice $Q_n = \frac{1}{R_Q}A_Q(x_n)^\top B_Q(x_n)$, cf. \cite[Theorem 3.1]{zhang2025tensor}.
Additionally, higher-order versions of TPA are explored, e.g. using latent factors $a^Q_r(x_n) \in \RR^{\Nh}, b^Q_r(x_n) \in \RR^{b}$ and $c^Q_r(x_n) \in \RR^{c}$, where $bc = \dhead.$
We refer the reader to \cite[Appendix C]{zhang2025tensor} for more details.

\subsection{Tensorized attention models}
\label{subsect:tensorial_models}

While attention computation captures pairwise information between words in a sentence, it may not be able to capture higher-order correlations, such as those between triplets of words. Rather than using tensorization as a tool to compress stacked weight matrices, this section shifts the focus to works that change the \textit{attention mechanism} directly, incorporating tensor structure, either through tensor products and/or decompositions to compute the attention scores~\eqref{eqn:att_matrix}.

We begin with one of the earliest and most influential works in this direction \cite{ma2019tensorized}, which proposes to replace the standard softmax-based attention entirely by learning a low-rank tensor-structured parameterization. We then discuss the approach of \cite{feng2024long}, which tensorizes the input data itself and applies attention mode-wise on the resulting tensor, yielding a memory-efficient, Kronecker-structured attention operator that imposes structured restrictions on token interactions. Finally, the work of \cite{sanford2023representational} establishes a theoretical basis for tensorized attention, showing that the usual attention layers cannot succeed at tasks that require detection of higher order correlations. They also introduce multilinear generalizations of the attention matrix designed to capture higher-order correlations in the data.

 \subsubsection{Tensorized (Multi-linear) Attention} 
\label{subsubsec:tensorized-MLA}
In~\cite{ma2019tensorized}, the authors proposed a tensorized transformer with Block-Term Tensor Decomposition (BTD) for compression. A \textit{tensorized attention mechanism} replaces the softmax-based attention scores by a learned low-rank multilinear interaction among queries, keys, and values. The new tensor attention is defined as follows. Let $\dhead\in\NN$ denote the embedding dimension used for queries and keys, and let $R\in\NN$ be a prescribed rank parameter. Given the matrices $\Q,\K\in\RR^{N\times \dhead}$ and $\V\in\RR^{N\times d}$ defined in \eqref{eqn:QKV}, the model introduces a trainable weight vector $g=(g_1,\dots,g_R)\in\RR^R$ and defines a third-order tensor $\mathcal{A}\in\RR^{N\times N\times N}$ with entries
\[
\mathcal{A}(i,j,m)
=
\sum_{r=1}^R g_r\, \Q(i,r)\,\K(j,r)\,\V(m,r),
\qquad
1\le i,j,m\le N,
\]
In other words, this representation corresponds to a rank-$R$ CP decomposition of $\mathcal{A} = \llbracket g; \Q,\K,\V \rrbracket$, with factor matrices $\Q,\K,\V$ and weights $g$.

To produce an output compatible with standard Transformer layers, the tensor $\mathcal{A}$ is mapped back to a matrix through a mode-1 unfolding. Recall ${A}_{(1)}\in\RR^{N\times N^2}$ denotes the mode-$1$ unfolding of the tensor $\mathcal{A}$. Then, up to a fixed permutation of columns, one can write
\[
{A}_{(1)} = \Q\,\mathrm{diag}(g)\,(\V * \K)^{\T},
\]
where $*$ denotes the Khatri-Rao (columnwise Kronecker) product. 

 The output of the tensorized attention block is then given by
\[
O = {A}_{(1)}\,\WO \in \RR^{N\times d},
\]
where $\WO\in\RR^{N^2\times d}$ is a learned output projection. In the multi-headed setting, several diagonal cores $g^{(1)},\dots,g^{(\Nh)}$ are used in parallel with shared matrices $\Q,\K,\V$, the corresponding tensors are averaged, and the same SplitConcat and output projection are applied.

The tensorized attention mechanism replaces the data-dependent softmax normalization by a learned multilinear interaction of fixed rank. As a result, the mapping from queries to outputs is linear for fixed keys and values, and the resulting representation does not enforce positivity or simplex constraints on attention weights. This removes the natural interpretation of attention as a selector or convex combination of values, while introducing explicit low-rank structure. 
Consequently, the expressivity of the model is governed by the rank parameter $R$ and the output projection $\WO$.

\subsubsection{Representational Strengths and Limitations of Transformers}
\label{par:tensor-attention} 

While attention computation captures pairwise information between words in a sentence, it cannot capture higher-order correlations, such as those between triplets of words. This is captured more concisely in~\cite{sanford2023representational}. In this work, the authors aim to understand approximation-theoretic properties of self-attention. More specifically, they ask whether self-attention can efficiently (with respect to parameter size) represent decomposable functions into pairwise interactions, and then take it a step further by asking the same question for triple-wise interactions. Informally, their results demonstrate that while self-attention can efficiently accomplish pairwise tasks, it cannot do the same for triple-wise tasks. 

To address this drawback, the authors propose a ``third-order tensor self-attention", which can accomplish triple-wise tasks~\cite{sanford2023representational}. The tensor self-attention is defined as follows.
Given query, key, and value matrices $\Q, \K_1, \K_2, \V_1, \V_2 \in \mathbb{R}^{N \times \dhead}$, where the additional key and value matrices can be obtained from different views or modalities~\cite{cao2024training}, the proposed generalization of the attention model is defined as 
\begin{equation}
D^{-1} A (\V_1 * \V_2) \in \mathbb{R}^{N \times \dhead},
    \label{eq:tensor_attn_comp}
\end{equation}
where 
\begin{equation}
A = \exp \left( \frac{Q(K_1 * K_2)^\top}{\sqrt{\dhead}}\right) \in \mathbb{R}^{N \times N^2},
    \label{eq:tensor_attn_mat}
\end{equation}
where $D = \diag(A \mathbf{1}_{N^2}) \in \mathbb{R}^{ N \times N}$, $\mathbf{1}_{N^2}$ is the $N^2$-dimensional all-ones vector, and $\exp()$ is applied element-wise. Intuitively, one can interpret $A \in \mathbb{R}^{N \times N^2}$ as an unfolded $N \times N \times N$ tensor. 

Under this new tensor-based attention model, it is shown that triple-wise interactions can be efficiently represented. However, the computational cost of such an attention mechanism becomes a major bottleneck. Instead of explicitly computing this tensor attention, the main contribution of~\cite{alman2023capture} focuses on the approximate computation of~\eqref{eq:tensor_attn_comp} in near linear time and they show that such a method exists if elements of the key, value, and query matrices to be bounded by $o(\sqrt[3]{\log N})$. The work of {\cite{liang2024tensor} also builds on the work of \cite{sanford2023representational}, providing an accelerated algorithm for tensor attention computation in near-linear time, and without restrictive bounds on the norm of the elements.

\subsection{Tensorial inputs}
\label{subsec:tensorial_inputs}

Tensorial or multiway data has become increasingly prevalent in modern scientific computing and data science applications.
For example, tensorial data arises frequently in applications such as weather prediction, climate modeling,  economic and financial data, discretizations of multivariate functions, simulations of partial differential equations, or 3D image reconstruction from, e.g., MRI data. 
However, standard dot-product attention often requires flattening of tensorial input data, which is often cost-prohibitive, destroys inherent structure, and obscures multi-way correlations or dependencies. In this section, we discuss a method that has been developed to compute attention within transformer architectures that preserves the tensorial structure of the input data. This is closely related to Section~\ref{subsect:tensorial_models}, where we discuss tensorized attention models for vector/matrix input data, but here, resizing the input/output to vector/matrix valued objects is not required.

   In \cite{omranpour2024higher}, the authors present an approach to higher-order attention for tensorial input data. Given input tensor $\mathcal{X} \in \RR^{N_1 \times N_2 \times d}$ for some hidden dimension $d$, they compute the query, key, and value tensors for each head $h$:
    \begin{align}
        \mathcal{Q}^{h} = \mathcal{X} \times_{3} (W^h_Q)^\top \in \RR^{N_1 \times N_2 \times \dhead}, \\
        \mathcal{K}^{h} = \mathcal{X} \times_{3} (W^h_K)^\top \in \RR^{N_1 \times N_2 \times \dhead}, \\
        \mathcal{V}^{h} = \mathcal{X} \times_{3} (W^h_V)^\top \in \RR^{N_1 \times N_2 \times \dhead},
    \end{align}
    where $\times_{3}$ denotes multiplication along the 3rd mode and $W_{Q/K/V/O}^h \in \mathbb{R}^{d \times \dhead}$.
    The attention scores $S^h \in \RR^{N_1 N_2  \times N_1 N_2}$ are then given by
    \begin{align}
        S^h = \textup{softmax} \left ( \frac{(Q^h_{(3)})^\top K^h_{(3)})}{\sqrt{\dhead}} \right ),
    \end{align}
    where $Q^h_{(3)}, K^h_{(3)} \in \mathbb{R}^{\dhead \times N_1 N_2}$ represent matricizations of $\mathcal{Q}^{h}$ and $\mathcal{K}^{h}$, respectively, along the 3rd mode.
    The output of the higher-order attention mechanism:
    $$O_{(3)} = \sum_h {\W_O^h} V_{(3)}^{h} S^h \in \mathbb{R}^{\dhead \times N_1 N_2 }, $$
    is then refolded back to the original tensor shape $\mathcal{O} \in \mathbb{R}^{N_1 \times N_2 \times d}$.
    
 Note that if naively implemented, the computational cost would be $\mathcal O\left (\dhead(N_1 N_2)^2\right )$, so the authors \cite{omranpour2024higher} propose the strategy of representing the attention matrix $S^h$ as a Kronecker decomposition: $S = \sum_{h} S_h^{(1)} \otimes S_h^{(2)}$, where $S_h^{(i)} \in \RR^{N_i \times N_i}$ is a factor matrix corresponding to the attention weights in the $i$th mode for head $h$. 
    Properties of Kronecker products are then exploited to apply $S^h$ to $\V^h_{(3)}$ without explicitly forming $S^h$, which reduces the complexity to $\mathcal{O}\left (\dhead N_i (N_1 N_2) \right )$.
    Lastly, low-rank tensor decompositions can be used to further drive down computational complexity.

       We end with additional references for tensorial attention models that assume tensor-structured input data. 
       Tensor-Augmented Transformers (TEAFormers)~\cite{kong2025teaformers} uses tensor input data but reformulates attention for (assumably, learnable) tensor-based weights, with applications to time-series data. Their method includes a compression step for the input tensor followed by a ``tensor-augmented multi-head attention'' computation, which uses the trained weight tensors $\mathcal{W}_{Q/K/V/O}$ to compute the layer's output tensor $\mathcal{O}$. Axial attention~\cite{ho2019axial} also targets inputs that are naturally represented as higher-order tensors (e.g., images or videos) and avoids flattening them into long token sequences. In that work, instead of forming a full $N\times N$ attention matrix, attention is applied sequentially along individual tensor modes (axes), such as rows and columns, yielding a factorized, mode-wise attention mechanism. We refer the reader to both of these works for more details on implementation and numerical results.
}

\section{Acknowledgments}
Part of this research was performed while the authors were visiting the Institute for Pure and Applied Mathematics (IPAM), which is supported by the National Science Foundation (Grant No. DMS-1925919), for the Research Collaboration Workshop\footnote{A.M. and D.N. were part of the organizing team for this workshop.}, ``Randomized Numerical Linear Algebra” (RNLA) 2025".

\section*{Code and Data Availability}
The code used to generate the QKV matrices used in the figures of this paper is publicly available at \href{https://github.com/rnla-transformers/qkv_extractor}{https://github.com/rnla-transformers/qkv\_extractor}. 
The repository includes the implementation, instructions for extracting the $\Q$, $\K$ and $\V$ matrices of an input text, and some examples. 
The Huggingface transformers \cite{wolf-etal-2020-transformers} python library is used in the process.
The input text we used for obtaining Figures~\ref{fig:sparsity_pattern}, \ref{fig:singular_values_qkv}, and~\ref{fig:singular_values_dotproduct} is the abstract of~\cite{han2023hyperattention}, in the form available at \href{https://arxiv.org/abs/2310.05869}{https://arxiv.org/abs/2310.05869}, licensed under CC BY 4.0  (\href{https://creativecommons.org/licenses/by/4.0/}{https://creativecommons.org/licenses/by/4.0/}). No endorsement by the original author is implied.  
The text was tokenized, and transformed into $\head Q h$, $\head K h$ and $\head V h$ through the use of the above qkv\_extractor library, using Meta’s   
\href{https://huggingface.co/meta-llama/Llama-3.2-1B}{Llama-3.2-1B}  model~\cite{grattafioriLlama3Herd2024}, released September 25, 2024,  
under the Llama 3.2 Community License Agreement.

\printbibliography
\appendix

\section{Theoretical Perspectives on Transformers}\label{sec:tf_theory}

Transformers are the main engine behind modern large language models, where the theory-practice gap remains wide. 
Fundamental theoretical questions concerning Transformers include what kinds of problems they can represent, what makes attention useful, and when standard training methods can actually discover those useful computations. This section provides a brief, accessible overview of recent explorations into how to address those questions.

\paragraph{Expressivity of transformers.}
A broad picture has emerged in which depth, positional information, and the structure of attention together determine what Transformers can compute.
With appropriate positional encodings, \cite{yun2020transformers} shows that multi‑head self‑attention with feed‑forward blocks is a universal approximator of sequence‑to‑sequence maps on a compact set.
\cite{giannou2023looped} shows that ``looped'' transformers (recycling a fixed block a variable number of times) can emulate programmable computers with a constant number of encoder layers, and are Turing‑complete.
At the other end of the spectrum, constructive results also prove Turing-completeness for idealized variants with \emph{hard} attention (which selects one location almost exactly) or \emph{monotone} attention (which moves through the sequence in order), together with unrealistically high numerical precision; these are useful as existence proofs, though less reflective of practical models~\citep{perez2021attention}.
\cite{sanford2023representational} studies basic representational tradeoffs, identifying tasks that separate attention from feed-forward networks, also called multilayer perceptrons (MLPs), and recurrent models, while clarifying the role of embedding dimension and communication complexity.
On the formal‑language side, \cite{liu2023transformers} shows that shallow decoders can shortcut finite‑state automata on an input sequence of length $T$, giving $O(\log T)$-depth simulators and often $O(1)$‑depth solutions for broad automata subclasses.
Meanwhile, \cite{merrill2024expressive} proves that adding intermediate decoding steps, i.e., chain‑of‑thought (CoT), strictly raises computational power. 
Further, \cite{merrill2025little} shows that even highly uniform Transformers with depth $\Theta(\log n)$ can express regular languages and graph connectivity, tasks that are inexpressible for constant‑depth Transformers under standard conjectures, suggesting that depth scaling can be far more effective than width or CoT steps for these classes.
\cite{roy2025fast} shows that higher‑order attention strengthens the basic dot‑product interaction by proposing 2‑simplicial (tri‑linear, determinant‑based) attention that natively models triple interactions.
Several works give universality/approximation‑rate results for simplified or efficient Transformer families~\citep{jiang2024approximation,de2024positional}, unveiling which efficiency tricks preserve/degrade theoretical expressivity.

\paragraph{Expressivity of fast attention.} 
Recent work connects fast attention mechanisms to parallel computation models, like the massively parallel computation (MPC) model. \citet{liu2025fast} introduces Approximate Nearest Neighbor Attention (ANNA), an LSH-inspired sub-quadratic attention primitive, and proves a sharp equivalence between ANNA-transformers and sublinear-memory MPC. In particular, it is demonstrated that ANNA retains the MPC-level expressivity previously exhibited for standard attention, but unlike standard attention, it can be simulated by MPC using a strongly sub-quadratic (near-linear for a large approximation factor) number of machines. This tighter correspondence transfers MPC round-complexity lower bounds to depth lower bounds for fast-attention transformers, while still permitting explicit fast-attention constructions for reasoning benchmarks such as Match2~\citep{sanford2023representational} and $k$-hop induction heads~\citep{sanford2024transformers}. Moreover, it is shown that constant-depth ANNA's can simulate constant-depth low-rank attention transformers.

\paragraph{Learnability of transformers under gradient descent.}
Beyond what can be represented by a transformer as a function class, another central question is whether and how efficiently common optimization algorithms like gradient descent (GD) find the transformer representation of a target function; this topic has been explored, among others, by~\citep{nichani2024transformers,wang2025learning,yang2025multi,goel2026training}.
Toward answering this question, \cite{nichani2024transformers} analyzes an in-context learning (ICL) task where the data are generated from a latent causal graph and proved that a (simplified two‑layer) transformer trained by GD learns to encode the causal structure in its first attention layer. 
At a high level, the model learns to pay more attention to tokens that are statistically informative about the next prediction. 
\cite{wang2025learning} studies a compositional target family, $k$-fold function compositions expressible by $O(\log k)$-depth transformers, and shows a sharp \emph{statistical–computational gap}: generic statistical query (SQ) learners need exponentially many samples, whereas GD succeeds polynomially under curricula that mix easy and hard instances.
\cite{yang2025multi} studies symbolic multi-step reasoning via path finding on trees, and shows that when the model is trained to generate chain-of-thought intermediate steps, GD can train even a one-layer multi-head transformer to solve reasoning tasks by inducing a multi-phase training trajectory in which different attention heads autonomously specialize and coordinate across subtasks, with generalization guarantees to unseen tree structures.
\cite{goel2026training} analyzes the training of a softmax self-attention layer on a linear regression problem and show that, in the infinite-data limit, the training objective can be rewritten as a weighted matrix factorization problem. Using the geometric understanding of the landscape of such a matrix factorization problem, they design a tailored first-order method (combining spectral initialization, regularization, and preconditioning) that avoids bad stationary points and converges globally at a geometric rate. This yields a clean scaling law in which the excess prediction error splits into a statistical part and an optimization part, with the latter decaying exponentially in the number of gradient steps.

\paragraph{Clustering of token embeddings.} 

Recently, the emergence of clusters in self-attention dynamics has been explored on some simplified attention mechanisms, which partially account for the effectiveness of the clustering methods utilized in Section \ref{sec:clustering}. In \cite{geshkovski2025mathematical} and subsequent works, the authors develop a mathematical framework for analyzing Transformers based on the interpretations among token embeddings as interacting particle systems, revealing the emergence of clusters over long time. Precisely, \cite{geshkovski2023emergence} considers the dynamics of $d$-dimensional token embeddings $\{x_{j}(t)\}_{j=1}^{N}$, where $x_{j}(t)$ evolves with respect to the depth $t$ of the layers. Rather than dealing with the discrete labeling $t\in \NN$ for layers, embedding vectors are treated as particles that evolve continuously, with their dynamics described by an interacting particle system satisfying the ODE below
\begin{align}
    \frac{\mathrm{d}}{\mathrm{d}t}x_{i}(t) =\sum_{j=1}^{N}P_{ij}(t) V x_{j}(t),\quad t\in [0, \infty)
\end{align}
for any $i\in [N]$, where $P_{ij}(t)$ are the entries of an $N\times N$ stochastic matrix $P(t)$, given by
\begin{align}
    P_{ij}(t) \coloneqq \frac{\exp\left (\<Q x_{i}(t), K x_{j}(t) \> \right )}{\sum_{l=1}^{N} \exp\left (\<Q x_{i}(t), K x_{l}(t) \>\right )},
\end{align}
with $Q$, $K$, $V$ representing the \emph{query}, \emph{key}, and \emph{value} matrices, respectively, defined in \eqref{eqn:QKV}. Here, the matrix $P(t)$ is called \emph{self-attention} matrix. The term attention stems precisely from the fact that $P_{ij}(t)$ captures the attention given by token $i$ to token $j$ relative to all tokens $l\in [N]$. In {\cite[Theorem 2.1]{geshkovski2023emergence}}, the authors proved that for any initial sequence of pairwise distinct tokens, $P(t)$ converges to a low-rank boolean matrix. When considering the rescaled token $z_{j}(t) = e^{-tV} x_{j}(t)$ for each $j\in [N]$, {\cite[Theorem 3.1]{geshkovski2023emergence}} proved that there exists a convex polytope $\mathcal{K}\subseteq \R^{d}$ such that for any $j\in [N]$, $z_{j}(t)$ converges to either $0$ or some point on the boundary of $\mathcal{K}$ as $t \to \infty$. Therefore, the clusters emerge. \cite{karagodin2024clustering} presents a modification of the self-attention dynamics to better reflect the practically relevant, causally masked attention used in transformer architectures. They prove the asymptotic convergence of token embeddings to a single cluster for arbitrary key-query matrices and a value matrix equal to the identity. \cite{chen2025quantitative} further investigates the long-term clustering of mean-field transformer models. They establish exponential rates of contraction to a Dirac point mass for any suitably regular initialization of token embeddings. They show that any suitably regular mean-field initialization synchronizes exponentially fast with some quantitative rates. Readers may refer to \cite{geshkovski2025mathematical} for a more comprehensive literature review.

\section{Detailed description of transMLA}
\label{subsec:transmla}

The objective of transMLA \cite{meng2025transmla} is to convert GQA-based models into DeepSeek-like MLA-based models able to take advantage of the optimized pipeline developed by DeepSeek. To do so, several steps are required: (i) first, an intermediate rewriting of GQA which would be equivalent to MLA without positional encoding is provided; (ii) then, a novel method is used to concentrate the features used for the positional information in the first few heads; (iii) finally, an approximation restricting the positional embedding to the first head allows us to obtain an approximate conversion to MLA. In this section we give some further details about the three aforementioned components of this approach.

\subsection{An intermediate encoding of GQA}

The intermediary rewriting of the GQA model consists in simply setting an intermediate shared KV embedding (a {proto} latent space) consisting in the concatenation of the key and value weight matrices for each group head, i.e.
\[
\WIKV=\begin{bmatrix}\head{\WK}{1}|&\cdots&|\head{\WK}{\Ng}|&\head{\WV}{1} |&\cdots &|\head{\WV}{\Ng}\end{bmatrix} =\begin{bmatrix}
    \WKd |& \WVd
\end{bmatrix}\in\RR^{d\times2\Ng\dhead},
\]
 where $\WKd,\WVd\in\RR^{d\times\Ng\dhead}$ are the ``down-projection'' key and value weight matrices into the proto latent subspace.
Note that $\WKd$ and $\WVd$ can equivalently be seen concatenation of the  $\head\WK g$ and $\head \WV g$ for all heads $g\in[\Ng]$. 
We then set, for each head $h$, the ``up-projection'' key and value weight matrices (from the proto latent subspace to the key and value spaces), ${\head{\WKu}{h}}\in\RR^{\Ng\dhead\times\dhead}$ and ${\head{\WVu}{h}}\in\RR^{\Ng\dhead\times\dhead}$, to be the identity matrices for corresponding group head $g_h=\lfloor\frac{h\cdot\Nh}{Ng}\rfloor$, i.e.
\begin{equation*}
{\head{\WKu}{h}}={\head{\WVu}{h}}=\Big[\overbrace{0\mid\cdots\mid\ 0}^{g_{h}-1 \text{ times}} \mid I_{\dhead} \mid  \overbrace{0\mid  \cdots\mid 0}^{{\Ng - g_h \text{ times}}}\Big]^\T.\\ 
\end{equation*}
Hence, the key and value matrices associated with any query head $h$ can be obtained as $$\head{\WK}{g_h} = {\WKd\head{\WKu}{h}},\qquad \head{\WV}{g_h} = {\WVd\head{\WVu}{h}}.$$

Then, we can define an intermediate ``latent'' representation $L^{\interm}$ as
$$L^\interm=
\begin{bmatrix}\head{\K}{1}\mid\cdots\mid\head\K{\Ng}\mid\head{\V}{1} \mid\cdots \mid\head{\V}{\Ng}\end{bmatrix} =X\WIKV\in\RR^{N\times2\Ng\dhead},
$$ with 
\[
\head \K h = L^\interm \begin{bmatrix}
    {\head{\WKu}{h}} \mid0
\end{bmatrix}= X\WKd\head{\WKu}{h},\qquad \head \V h = L^\interm \begin{bmatrix}
    0 \mid {\head{\WVu}{h}} 
\end{bmatrix}= X\WVd\head{\WVu}{h}.
\]
The application of RoPE can then be written as $$\head \Q h^{R}{\head \K h^{R}}^\T=\head \Q h^R{\left(X\WKd {\head{\WKu}{h}}\right)^{R}}^{\T},$$
where we recall from Section \ref{sec:RoPE} that $(\cdot)^R=\RoPE_{\Ng}\left(\cdot\right)$ where $\RoPE_{\Ng}$ consists in the application of the $\dhead$-dimensional RoPE rotations repeated $\Ng$ times, that is, once for each  head.
As the RoPE operation commutes with the identity and the null matrix which are the diagonal blocks of ${{\head{\WKu}{h}}}^\T$, we can move it through the RoPE operation. This allows to rewrite the query and key matrices into an intermediate form wherein only the query matrix depends on the attention head, this reads
\begin{align*}
\head{\Q}{h}^\interm&=\left(\head{\Q}{h}{{{\head{\WKu}{h}}}^\T}\right)^R\in\RR^{N\times\Ng\dhead},\\ \K^\interm&=\left(X\WKd\right)^R\in\RR^{N\times\Ng\dhead}.
\end{align*} For any head, $1\le h\le\Nh$, this formulation indeed gives us $(\head{\Q^{\interm}}{h})(\K^{\interm})^{\T}=(\head{Q^R}{h})(\head{\K^{R}}{h})^\T$.
However, this encoding in addition to being inefficient is also not equivalent to DeepSeek's MLA due to the different RoPE implementations for GQA and MLA and an additional step to separate the rotational embeddings is required.

\subsection{RoRoPE}
The idea behind RoRoPE is to construct a matrix $U^{\downarrow}\in \RR^{\Ng\dhead\times\Ng}$ 
 to concentrate the ``positional features'' in the first output heads which then allows the approximation of the intermediate model by an MLA model by restricting the application of RoPE to the first output head. To do so, for each $0\le \ell \le \dhead/2-1$, TransMLA defines a head-dimension orthogonal matrix $U_\ell$ of size $\CC^{\Ng\times\Ng}$ whose purpose is to extract the positional features contained in each key head and concentrate them into a vector as in Latent Attention's $\KRL$. The method described in \cite{meng2025transmla} extracts this positional information using the $\Q$, $\K$ and $\V$ matrices of an input calibration dataset.
In this subsection we will resort to using tensor notation, i.e. we consider the heads of the matrices to be a third dimension.
\paragraph{Positional information extraction.}
In this paragraph, we consider the tensors for the calibration dataset to be denoted by
$$
\mathcal{Q}(i,j,h)=\head\Q h(i,j),\quad \mathcal{K}(i,j,g)=\head\K g(i,j),\quad \mathcal{V}(i,j,g)=\head\V g(i,j),
$$ and we resort to the same tensor notation for any weight matrix depending on $h$.
To extract the positional information, for simplicity we need to consider the tensors as indexed complex matrices, we thus define, for each $0\le \ell \le \dhead/2-1$, the following matrices:
\[
\widehat{\WK_\ell}=\WK(:,2\ell+1,:)+\mathfrak i\WK(:,2\ell+2,:)\in\CC^{\dhead\times \Ng}
\text{ and }
\widehat{\head{\K}{\ell}}=\K(:,2\ell+1,:)+\mathfrak i\K(:,2\ell+2,:)\in\CC^{N\times\Ng},
\] 
where $\mathfrak i=\sqrt{-1}$.
This consists in embedding the $(2\ell+1)$-th and $(2\ell+2)$-th output dimensions of each head in the complex plane; note that the $\ell$ indices are zero-indexed.

For a given token index $j$ (from the calibration dataset), $1\le n\le N$, the effect of each block, for $0\le \ell\le \dhead/2-1$, of the RoPE matrix on $\mathcal K$ can now be seen as a complex rotation of the same associated angle, and can be written as
\[
    \widehat{\K_\ell^R}(j,:)=\widehat{\K_\ell}(j,:)\exp(\mathfrak i j \theta_\ell ).
\]
The TransMLA method consists in obtaining, for each angle $\theta_\ell$, the unitary matrix $U_\ell\in\RR^{\Ng\times\Ng}$ concentrating the features which are important for the $\ell$-th frequency into the first $m\in\NN^*$ heads. Hence, $U_\ell$ maximizes the norm, 
\[
\|K_\ell^R [U_\ell(:,:m)]\|^2_{\mathrm{F}}=\Tr\left([U_\ell \widehat{\K_\ell^R}^\T\widehat{\K_\ell^R} U_\ell] (:m,:m)\right).
\]
As is readily evident, given the singular value decomposition of $\widehat{\K^R_\ell}$ this maximizer turns out to simply be the $m$ first right singular vectors. 
Equivalently, using the methodology of the article, the maximizer corresponds to the best rank-$m$ approximation of the covariance matrix $\widehat {K_\ell^R}^\T\widehat {K_\ell^R}$, which can be obtained from the eigendecomposition
\[
    \widehat{\K_\ell^R}^\T\widehat{\K_\ell^R}=U_\ell\Lambda_\ell U_\ell^\T\in\CC^{\Ng\times\Ng},
\]
where $\Lambda_\ell$ is the diagonal matrix containing the eigenvalues. 

\paragraph{Rewriting the intermediate model.}
The idea is that now the most important features obtained through RoPE for the $\ell$-th frequency will be contained in the first few heads, allowing the application of the RoPE embedding to the vector contained in the first head only with the minimum error.
The following identity is provided,
\[
\head{\Q^\interm}{h}(\head{\K^\interm}{h})^\T  = \left(\head{\Q}{h}{{{\head{\WKu}{h}}}^\T}U^{\downarrow}\right)^R\left(X\WKd U^{\downarrow}\right)^{R,\T}
\]
where we set, for any $g,g'\in[\Ng]$ and $j\in[\dhead]$, $$U^\downarrow\left(j+\dhead\times (g-1),j+\dhead\times (g'-1)\right) = U_{\lfloor j/2\rfloor}(g,g')$$  and zero elsewhere. This means that the application of $U^\downarrow$ mixes the vectors corresponding to the same dimension $j$ from different heads to concentrate the positional information into the first heads. 
This also shows that we can, without any approximation to the model, combine the $U^\downarrow$ matrices into the $\WKd$ and ${\head{\WKu} h}^\T$ matrices.
The authors also provide an alternative method of calculating the matrices $U_\ell$ by concatenating the vectors corresponding to multiple ``similar'' frequencies and sharing the $U_\ell$ matrices between them: this method is called FreqFold.
It basically consists in finding a matrix $U_G$, for each group $G$ of frequencies, which concentrates the most important shared features in the $m$ first vectors. The matrix $U_G$ is then given by the first $m$ vectors of the eigendecomposition of the sum of the covariance matrices $\sum_{f\in G}\widehat{K_f}^\T\widehat{K_f}$. 

\subsection{MLA Approximation}
By design, the first ``virtual'' head $g'=1$ now contains the part of the vectors which is important for the rotational embeddings. The application of RoPE can thus be (approximately) restricted to $h'=1$ resulting in the following expression:
\begin{multline*}
\head{\Q^\interm}{h}(\head{\K^\interm}{h})^\T  \approx \left(\head{\Q}{h}\left(\left[{\head{\WKu}{h}}^\T U^\downarrow\right](\underbrace{:\dhead}_{ g'=1},:)\right)\right)^R\left(X\left(\left[\WKd U^\downarrow \right](\underbrace{:\dhead}_{ g'=1},:)\right)\right)^{R,\T}\\
 + \head{\Q}{h}\left(\left[{{\head{\WKu}{h}}}^\T U^\downarrow\right](\dhead:,:)\right)\left(X\left(\left[\WKd U^\downarrow\right](\dhead:,:)\right)\right)^{\T}.
\end{multline*}

In order to convert this model to MLA, TransMLA separates the first head, which contains the features most important for the positional encoding, and sets it as the Latent Attention $\KRL$, i.e. 
\[ 
\WRK =\left[\WKd U^\downarrow \right](:\dhead,:).
\] 
The weight matrices for the other ``new'' heads, on which we apply No Positional Embeddings (NoPE), are set to be
\[ \W^{K\downarrow,\NoPE}=\left[{{\head{\WKd}{h}}}^\T U^\downarrow\right](\dhead:,:),\qquad \head{\W^{K\uparrow,\NoPE}}{h}=\left[{{\head{\WKu}{h}}}^\T U^\downarrow\right](\dhead:,:),\] where we denote, for any matrix $M\in\RR^{d_1\times d_2}$ with $d_1,d_2\in\NN^*$, any integer $1\le m <d_1$ and $1\le n <d_2$, by $M(m:,n:)$ the submatrix composed of the last $d_1-m$ rows and last $d_2-n$ columns of $M$.
The intermediate embedding for the remaining Key-Value space is then given by
$$\W^{\mathrm{IKV},\NoPE} = \begin{bmatrix}
    \W^{K\downarrow,\NoPE}| & 
    \WVd
\end{bmatrix}\in\RR^{d \times (2\Ng-1)\dhead}.$$ 
To obtain the latent model weights, a low-rank decomposition for this matrix is used. 
The authors use the same method to that used in the positional encoding case to achieve this decomposition, using a calibration dataset with input matrix $X$.
More specifically, a procedure called ``KV balancing'' is used to balance the norms of the $X\W^{K\downarrow,\NoPE}$ and $X\WVd$ sub-matrices.
A coefficient $\alpha\in\RR$ is introduced and the matrix \[C= \begin{bmatrix}
    \alpha\W^{K\downarrow,\NoPE}| & 
    \WVd
\end{bmatrix},\] is considered. 

The right singular vector matrix $P^\T$ of $C$ is then obtained through the eigendecomposition of the covariance matrix $ C^\T C$.
Given the parameter $r\in\NN^*$, the low-rank approximation is then obtained by setting
\[ \WLKV = \W^{\mathrm{IKV},\NoPE} [P](:, :r)\in\RR^{d\times r},\]
and for all head $1\le h\le \Ng$, \[ \head\WUK h  = \frac{1}{\alpha}\head{\W^{K\uparrow,\NoPE}}{h} [P](:, :r)^\T\in\RR^{r\times\dhead},\qquad
    \head\WUV h = \head{\WVu}{h} [P](:, :r)^\T\in\RR^{r\times\dhead}.\]  
For the latent query subspace, in a similar way, the eigendecomposition of $C_Q=X\begin{bmatrix}\head\WQ 1 |& \cdots &| \head\WQ \Nh\end{bmatrix}$ is used to obtain $P_Q^\T$.
Then, given $r_Q\in\NN^*$, the low-rank decomposition of the model is obtained by setting 
\[
\WLQ = \begin{bmatrix}
    \head\WQ 1 |& \cdots &| \head\WQ \Nh
\end{bmatrix} [P_Q](:, :r_Q)\in\RR^{d\times r_Q},
\] 
and for all head $1\le h\le \Nh$,
\[ \head\WUQ h = [P_Q](:, :r_Q)^\T\in\RR^{r_Q\times \dhead},\qquad\head\WRQ h =\head\WUQ h\left[{\head{\WKu}{h}}^\T U^\downarrow\right](:\dhead,:)\in\RR^{r_Q\times \dhead}.\]
This concludes the conversion of a GQA model into an approximate MLA model.

\end{document}